\documentclass{svjour3}

\smartqed

\usepackage{graphicx}
\usepackage{graphics}
\usepackage{upgreek}
\usepackage{float}
\usepackage{amsmath}
\usepackage{wrapfig}
\usepackage{amssymb}
\usepackage{listings}
\usepackage{cite}
\usepackage{subcaption}
\usepackage{multirow}
\usepackage[percent]{overpic}
\usepackage{multicol}
\usepackage{color}
\usepackage{tikz}
\usepackage{calc}
\usepackage[hidelinks]{hyperref}

\newcommand{\etal}{et~al.}

\newcommand{\hfd}{\hat{f}^{\delta}}

\newcommand{\hfI}{\hat{f}^{\delta I}}

\newcommand{\hfC}{\hat{f}^{\delta C}}
\newcommand{\hfD}{\hat{f}^{\delta D}}

\newcommand{\hgI}{\hat{g}^{\delta I}}

\newcommand{\hgD}{\hat{g}^{\delta D}}

\newcommand{\hud}{\hat{u}^{\delta}}

\newcommand{\bhud}{\hat{\mathbf{u}}^{\delta}}

\newcommand\hb[1]{\hat{\mathbf{#1}}} 
 
\newcommand\hbd[1]{\hat{\mathbf{#1}}^{\delta}} 
\newcommand\bd[1]{\mathbf{#1}^{\delta}} 
\newcommand{\hnabla}{\hat{\nabla}} 

\newcommand\px[2]{\frac{\partial #1}{\partial {#2}}}

\newcommand\dx[2]{\frac{\mathrm{d} #1}{\mathrm{d} #2}}

\usepackage{calc}
\newlength\myheight
\newlength\mydepth
\settototalheight\myheight{Xygp}
\newcommand*\inlinegraphics[1]{%
  \settototalheight\myheight{Xygp}%
  \settodepth\mydepth{Xygp}%
  \raisebox{-\mydepth}{\includegraphics[height=\myheight]{#1}}%
}
\newcommand\orcid[1]{\href{https://orcid.org/#1}{\inlinegraphics{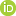}}}

\DeclareRobustCommand\dotted{\tikz[baseline=-0.6ex]\draw[thick,dotted] (0,0)--(0.54,0);}

\usepackage[titletoc]{appendix}
\usepackage[section]{algorithm}

\graphicspath{ {Figs/} }

\begin{document}

\title{Effect of Mesh Quality on Flux Reconstruction in Multi-Dimensions}

\author{Will Trojak \and
        Rob Watson \and
        Ashley Scillitoe \and
        Paul G. Tucker}
 
\institute{Will Trojak (corresponding author)\at
           Department of Engineering, University of Cambridge, Cambridge, UK, CB2 1PZ \\
           \email{wt247@cam.ac.uk} \\
           \and
           Rob Watson \at 
           School of Mechanical and Aerospace Engineering, Queen's University
           Belfast, Belfast, UK, BT9 5AH \\
           \and
           Ashley Scillitoe \at
           Alan Turing Institute, Kings Cross, London, UK, NW1 2DB\\
           \and
           Paul G. Tucker \at
           Department of Engineering, University of Cambridge, Cambridge, UK,
           CB2 1PZ}
           
\date{Received: 18 September 2018}
\maketitle

\begin{abstract}
	Theoretical methods are developed to understand the effect of non-uniform grids on Flux Reconstruction (FR) in multi-dimensions. The analysis reveals that the same effect of expanding and contracting grids is seen in two dimensions as in one dimension. Namely, that expansions cause instability and contractions cause excess dissipation. Subsequent numerical experiments on the Taylor-Green Vortex with jittered elements show the effect of localised regions of expansion and contraction, with an initial increase in the kinetic energy observed on non-uniform meshes. Some comparison is made between second-order FR and second-order finite volume (FV). FR is found to be more resilient to mesh deformation, however, FV is found to be more resolved when operated at second order on the same mesh. In both cases, it is recommended that a kinetic energy preserving/conservation formulation should be used as this can greatly increase resilience to mesh deformation.
\end{abstract}



\section{Introduction}\label{sec:intro}
	Since the inception of Spectral Volume methods by Wang~\cite{Wang2002}, the trajectory of high order methods has trended towards the Flux Reconstruction method of Huynh~\cite{Huynh2007} and Vincent~\etal~\cite{Vincent2010}. This approach draws on the work of Finite Elements, see Brenner and Ridgway-Scott~\cite{Brenner2008}, enabling the high performance of Flux Reconstruction on heterogeneous computing --- as can be seen in the highly efficient use of vast computing resources by Vincent~\etal~\cite{Vincent2017}. However, the move towards high order was not born out of a need for more efficient use of modern HPC environments. For example, Brandvik and Pullan~\cite{Brandvik2011} showed that high throughput could be obtained using second-order Finite Volume (FV) methods. Instead, the main motivating factor has been the increased uptake by industry of turbulence resolving methods, such as Large Eddy Simulation (LES), as this allows for far better exploration of flow physics and moves towards the aim of computational wind tunnels. The main feature of LES is the modelling of the very smallest scales of motion, which removes the need for Direct Numerical Simulation (DNS) levels of resolution. However, Chow and Moin~\cite{Chow2003} and Ghosal~\cite{Ghosal1996} showed that, for LES, the need to keep the truncation error small to enable the sensible use of sub-grid scale models meant that the grid requirements were demanding. A move to higher order would mean that the scaling of the truncation error was to a higher power of grid spacing --- thus lowering the grid requirements and decoupling the scaling of aliasing error and truncation error. Hence, for wall-resolved LES, calculations are often impractical unless the improved mesh resolution requirements of high order methods are considered.    

    The analytical understanding of Flux Reconstruction has been explored to a large extent in the work of Vincent~\etal~\cite{Vincent2011}, Jameson~\etal~\cite{Jameson2012}, and Castonguay~\etal~\cite{Castonguay2014}, where the stability of linear advection, advection-diffusion, and non-linear problems were presented. The key findings were the energy stability of FR on linear problems, and the condition for energy stability on non-linear problems. In addition, by investigating the dispersion and dissipation characteristics of FR, the existence of superconvergence after temporal integration and the corresponding CFL limits were found. This work was limited to one dimension --- although still applicable, the investigation of the exact behaviour of FR in higher dimensions has been limited, such as that of Williams and Jameson~\cite{Williams2014} and Sheshadri and Jameson~\cite{Sheshadri2016}. This work focused primarily on the proof of the Sobolev type energy stability in 2D in a similar manner to that of Hesthaven and Warburton~\cite{Hesthaven2007}, alongside some numerical studies performed for validation.  

    The advantage of FR --- that leads to high performance on heterogeneous and massively parallel architectures --- is its unstructured and sub-domain nature. Unstructured grids also allow far more complex geometries to be considered, but the resulting meshes experience deformation, expansion, and contraction of the elements. We wish to characterise the performance of FR under these conditions, and so far the effect of linear mesh deformation on FR has been considered in one dimension by Trojak~\etal~\cite{Trojak2017a}. Therefore, we make use of the seminal work of Lele~\cite{Lele1992}, in which the dispersion and dissipation of finite difference methods were considered in both one and two dimensions. We wish to repeat this process for FR, but extend it to also consider deformed grids.

    In this paper, we present an extension to the one-dimensional analytical work of Vincent~\etal~\cite{Vincent2011} and Trojak~\etal~\cite{Trojak2017a}. This extension will be shown for a two-dimensional case on quadrilaterals with rectilinear mesh stretching, but could also be performed on higher dimensional hypercubes. From the basis of this more general von Neumann analysis, the behaviour of FR on linearly mapped meshes can be explored. The investigation has been restricted to linear transformation as these are of key importance for complex industrial simulations due to there fundamental nature. For example, they occur in meshes where mesh generators have simply tessellated elements to fill the domain. Therefore, understanding their character is key, however, we should point to some recent work that has numerically investigated curved meshes~\cite{Mengaldo2016,Yu2014}.

     The aim of this work is to understand the effect of moving to higher dimensionality on key metrics governing scheme performance, such as CFL limit, dispersion, and dissipation. Finally, the Taylor-Green vortex will be used to understand how deformed meshes affect full Navier-Stokes calculations, with reference calculations performed by an industrial second order finite volume method.
\section{Flux Reconstruction}
	Flux Reconstruction ~\cite{Castonguay2012,Huynh2007} (FR) applied to the linear advection equation will form the basis of the initial investigation to be carried out, and for the reader's convenience, an overview of the scheme is presented here. However, for a more detailed understanding, the reader should consult Castonguay \cite{Castonguay2012} or Huynh \cite{Huynh2007}. This 1D scheme can be readily converted to two and three dimensions for quadrilaterals and hexahedrals, respectively. First, let us consider the one-dimensional advection equation:
	\begin{equation}
		\px{u}{t} + \px{f}{x} = 0
	\end{equation}
	
The FR method is related to the Discontinuous Galerkin (DG) method~\cite{Reed1973} and makes use of the same subdivision of the domain into discontinuous sub-domains:
	\begin{equation}\label{eq:Domain}
		\mathbf{\Omega}  = \bigcup_{n=1}^{N}{\mathbf{\Omega}_n}
	\end{equation}

\begin{figure}
		\centering
		\begin{subfigure}[b]{0.35\linewidth}
			\centering
			\includegraphics[width=\linewidth,trim= 0mm 0mm 0mm 0mm,clip=true]{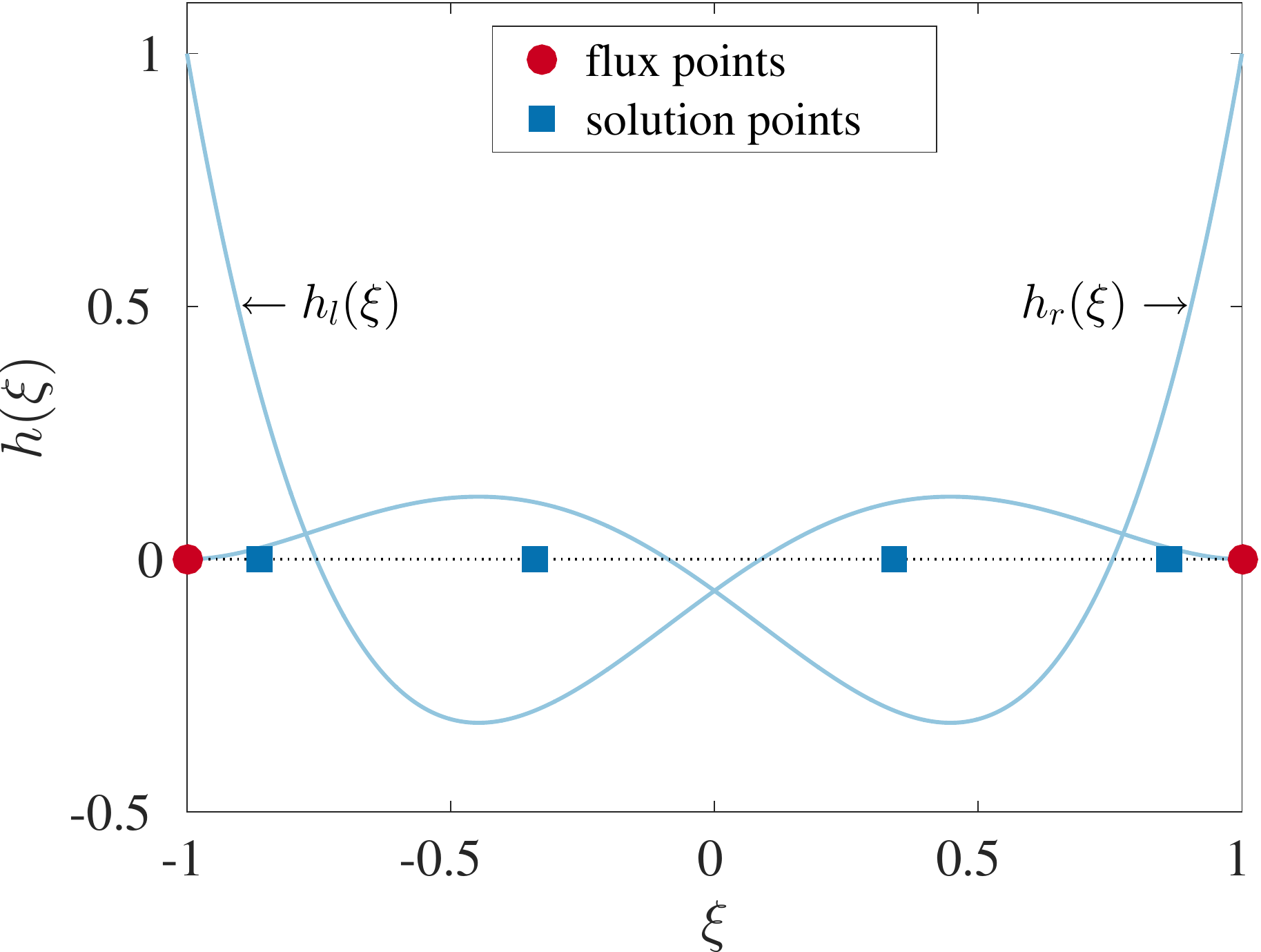}
			\caption{Flux and solution point layout for $p = 3$ in $\hb{\Omega}$, with corresponding left and right Huynh correction functions.}
			\label{fig:1D_layout}
		\end{subfigure}
		~
		\begin{subfigure}[b]{0.5\linewidth}
			\centering
			\includegraphics[width=\linewidth]{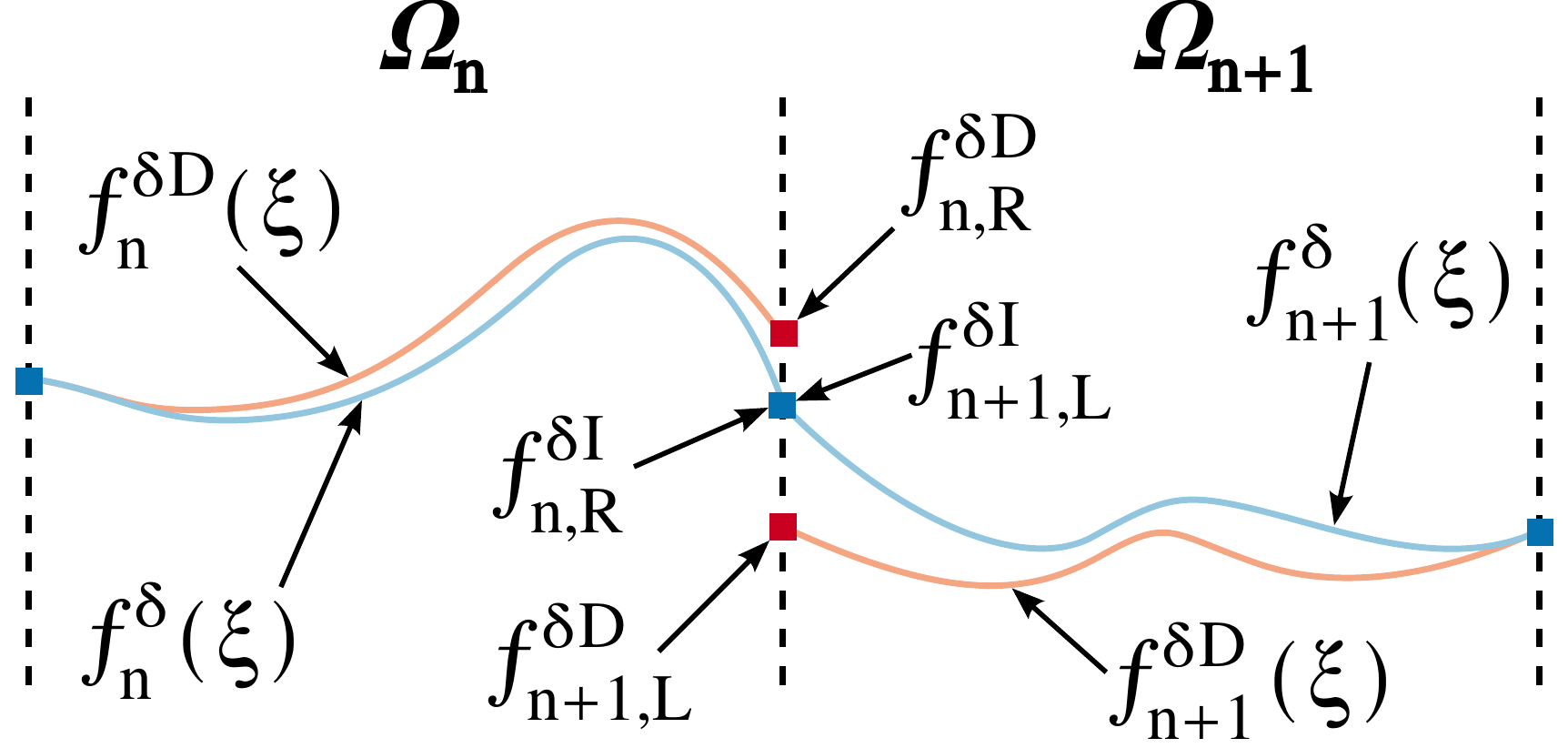}
			\caption{Diagram showing the procedure of correcting one, central, interface.\\ \:\:\: \\}
			\label{fig:1D_interface_u}
		\end{subfigure}
		
		\caption{Point layout in $\hb{\Omega}$ for $p=3$ and cell interface topology.}
		\label{fig:fr_interface}
	\end{figure}

	Within the standardised sub-domain, $\hb{\Omega} \in \mathbb{R}^d$, computational spatial variables are defined. When $d=1$,  $\hb{\Omega} = [-1,1]$, using $\xi$ to denote the value taken. This computational space is discretised with $(p+1)^d$ solution points, and $2d(p+1)^{d-1}$ flux points, placed at the edges of the sub-domain. The solution and flux point locations are determined using a tensor grid of a 1D quadrature. Fig.~\ref{fig:1D_layout} shows a 1D example of this. To transform from $\mathbf{\Omega}_n \rightarrow \hb{\Omega}$, a Jacobian $J_n$ is defined such that:
	\begin{equation}
		\hud = \hud(\xi;t) = J_n u^{\delta}(x;t) 
	\end{equation}	
	With this domain set up, we now proceed with defining the steps to construct a continuous solution from the discontinuous segments. The first stage is to define a local solution polynomial in $\hat{\mathbf{\Omega}}$ using Lagrange interpolation.
	\begin{align}
		l_k(\xi) &= \prod_{i=0,i\ne k}^{p}{\frac{\xi - \xi_i}{\xi_k - \xi_i}}\\
		\hud(\xi) &= \sum_{i=0}^{p}{\hud_il_i(\xi)}
	\end{align}
	
	Repeating the interpolation for the discontinuous flux in $\hat{\mathbf{\Omega}}$:
	\begin{equation}
		\hfD = \hfD(\xi,t) = \sum_{i=0}^{p}{\hfD_il_i(\xi)}
	\end{equation}
	Here we define $\hfD$ as the transformed discontinuous flux polynomial. Now using the Jacobian and the solution polynomials, the primitive values can be calculated in the physical domain $\mathbf{\Omega}_n$:
	\begin{equation}
		\overline{u}^{\delta}(x) = \frac{\hud(\xi)}{J_n} = \sum_{i=0}^{p}{u^{\delta}_i l_i(\xi)}
	\end{equation} 
		
	The primitive polynomial can then be interpolated to the interface and define as $\hud_l = \hud(-1)$ and $\hud_r = \hud(1)$. The values at the interface, $I$, then allow for a common interface flux, $f^{\delta I}_{I}$, to be calculated in the physical domain. The is shown graphically in Fig.~\ref{fig:1D_interface_u}. For a general case, this is done using a Riemann solver on the primitives at the interface, such as: Roe~\cite{Roe1981}; flux vector splitting~\cite{vanLeer1982}; or HLL~\cite{Harten1983}. In order to get a spatially continuous solution over $\mathbf{\Omega}$, the common interface flux has to be incorporated into the solution. For FR this is done by using a correction function to propagate the corrected flux into $\mathbf{\Omega}_n$. The definition of the correction function was shown to be important in the determination of the characteristics of FR, Vincent~\etal~\cite{Vincent2011}. Primarily the correction function is a polynomial which, in one dimension, has the boundary conditions:
	\begin{align}
		h_L(-1) = h_R(1) = 1& \\
		h_R(-1) = h_L(1) = 0 	
	\end{align}
	Beyond this, several sets of stable correction functions have been defined, firstly unified by Vincent~\etal~\cite{Vincent2010} and later expanded in \cite{Vincent2015,Trojak2018, Trojak2018b,Trojak2018f}. However, in this paper, we will focus on the correction function defined by Huynh~\cite{Huynh2007}, in particular, the Huynh $g_2$ correction function, which is shown in Fig.~\ref{fig:1D_layout}. We will also consider the correction that can recover Nodal Discontinuous Galerkin (NDG)~\cite{Hesthaven2008} as it provides a good point of comparison, due to DG relevance and maturity. It should be noted that the NDG correction function will only recover NDG in FR for homogeneous linear flux functions, due to the different mechanisms of aliasing. At times we will also explore some of the effects of correction function in two dimensions and for this, we will restrict ourselves to the Original Stable FR (OSFR) family of \cite{Vincent2010}, for brevity. This is a one-parameter family of the correction functions for which we will call the control parameter $\iota$.
	
	The correction to the flux function is then calculated using the difference between the discontinuous and common interface values and a correction function. The correction is then defined as:
	\begin{equation}
		\hfC = (\hfI_L - \hfD_L)h_L(\xi) + (\hfI_R - \hfD_R)h_R(\xi)
	\end{equation}
	and, hence the corrected continuous gradient of the flux is then:
	\begin{align}\label{eq:discrete_1D_full}
		\px{\hfd}{\xi} &= \dx{\hfD}{\xi} + \dx{\hfC}{\xi} \\
		& =  \sum_{j=0}^{p}\hfD_j\dx{l_j(\xi)}{\xi} + (\hfI_L - \hfD_L)\dx{h_L(\xi)}{\xi} + (\hfI_R - \hfD_R)\dx{h_R(\xi)}{\xi}
	\end{align}

	Finally, the solution is advanced in time following Eq.~(\ref{eq:discrete_1D_ad}) --- this can be done via a sensible choice of temporal integration. 
	
	\begin{equation} \label{eq:discrete_1D_ad}
		\px{\hud}{t} = - \px{\hfd}{\xi}
	\end{equation}
	The method detailed here was shown for simplicity in one dimension, but in subsequent sections, the method for extending this to a high dimension will be detailed. We will briefly state that to increase the dimensionality of the method a tensor product is used, which is this same method as is used in analysis and in the formal implementation of a solver for hypercube elements.
\section{Two-Dimensional Von Neumann Analysis}\label{sec:vn_method}
\begin{figure}
		\centering
		\begin{subfigure}[b]{0.46\linewidth}
			\centering
			\includegraphics[width=\linewidth,trim= 0mm 0mm 0mm 0mm,clip=true]{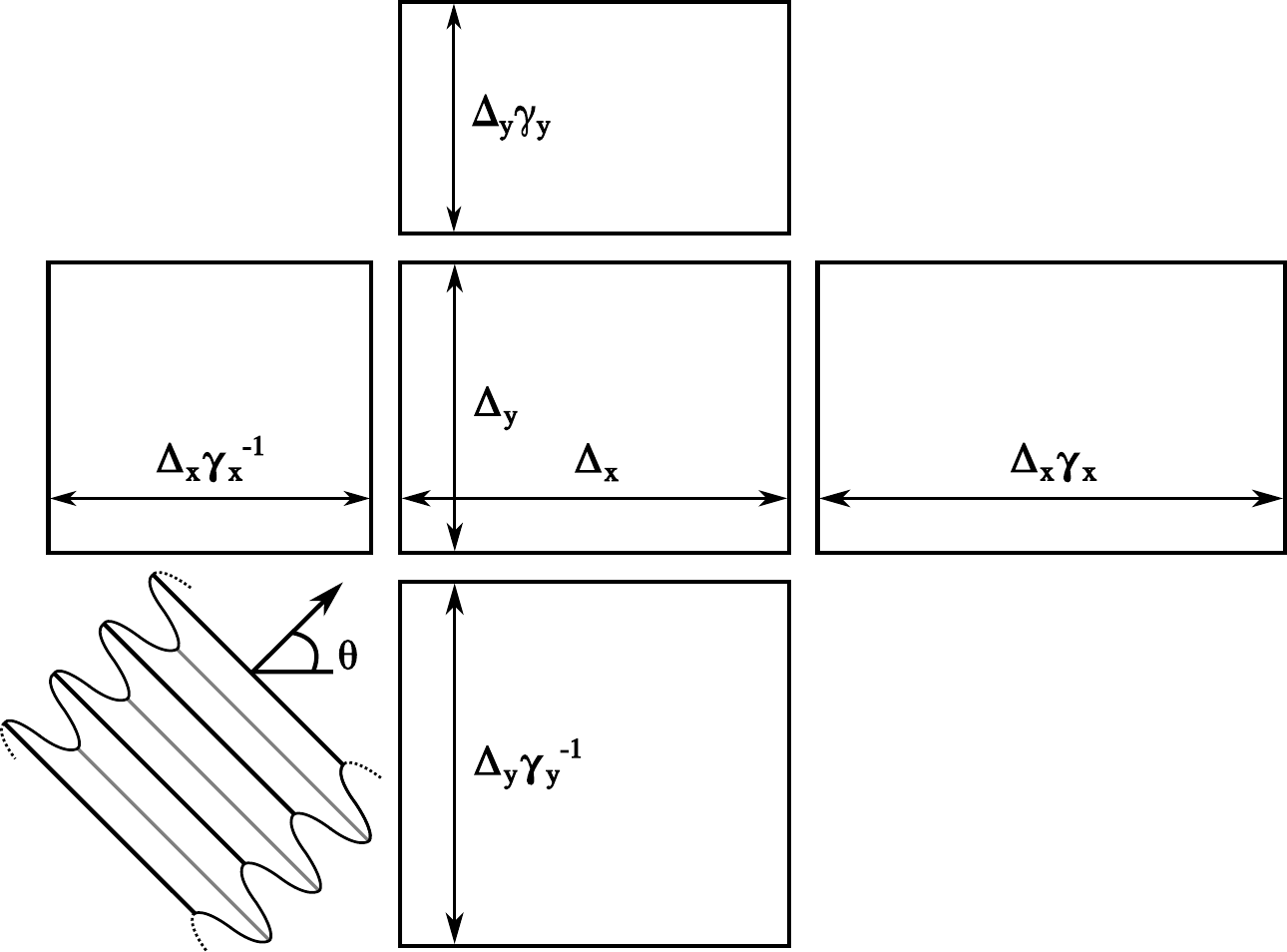}
			\caption{Schematic showing inclined plane wave passing through a cell with a geometrically transformed rectilinear stencil.}
			\label{fig:2D_schematic}
		\end{subfigure}
		~
		\begin{subfigure}[b]{0.38\linewidth}
			\centering
			\includegraphics[width=\linewidth,trim= 0mm 0mm 0mm 0mm,clip=true]{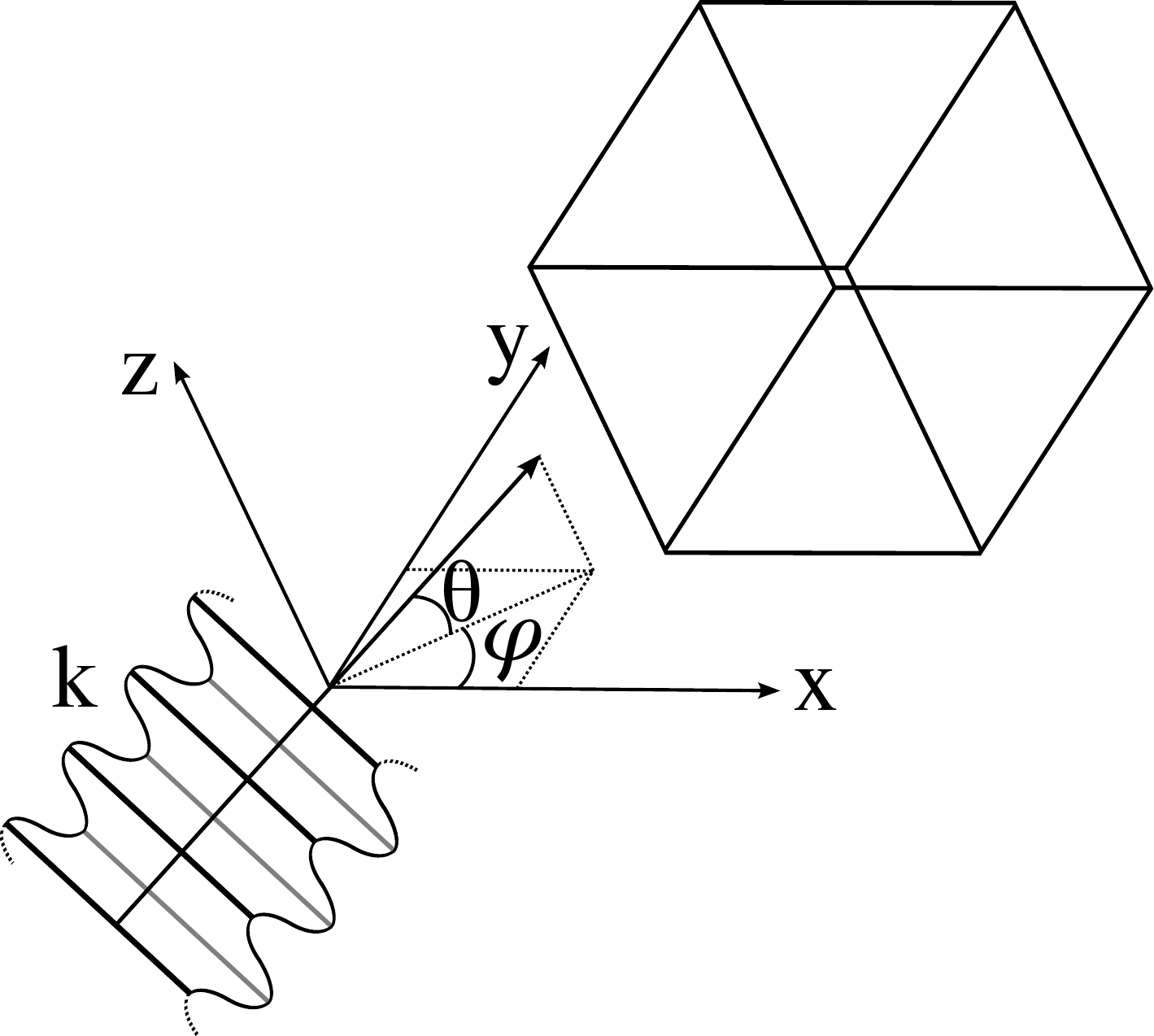}
			\caption{Schematic showing inclined plane wave passing through a cell with a geometrically transformed rectilinear stencil.}
			\label{fig:3D_schematic}
		\end{subfigure}
		\caption{Linear advection schematic for two and three dimensions.}	
	\end{figure}

The procedure for investigating the dispersion and dissipation properties of finite element methods has been laid out in some detail by Huynh~\cite{Huynh2007}, Hesthaven and Warburton~\cite{Hesthaven2007}, and Vincent~\etal~\cite{Vincent2011}. It is broadly classified as a von Neumann analysis. The procedure was, however, only performed in 1D, with critical insight into the analytical performance of FR when applied to more realistic problems overlooked. Extension of the analysis to higher dimension domains was performed by Lele~\cite{Lele1992} for various finite difference schemes. This did, however, avoid the increased complexity of finite element von Neumann analysis. To begin our extension we introduce the 2D linear advection equation:
	\begin{align}
		\px{u}{t} + \nabla\cdot\mathbf{F} &= 0 \\ \label{eq:2d_lin_ad} 
		\mathbf{F} = \begin{bmatrix} f\\ g \end{bmatrix} &=  u\mathbf{a} = \begin{bmatrix} au \\ bu \end{bmatrix}
	\end{align} 	
	Flux reconstruction then uses the superposition of the discontinuous and corrected flux divergence, meaning Eq.~(\ref{eq:2d_lin_ad}) can be rewritten as:
		\begin{equation}
			\px{\mathbf{u}_{i,j}}{t} = - \nabla\cdot\mathbf{F}^{\delta D}_{i,j} - \nabla\cdot\mathbf{F}^{\delta C}_{i,j}  
		\end{equation}
		Taking the following definition of the Jacobian, the computational-physical domain transformation can be defined:
		\begin{align}
			\mathbf{G} &= \begin{bmatrix} \px{x}{\xi} \: \px{y}{\xi} \\ \px{x}{\eta} \: \px{y}{\eta}  \end{bmatrix} = \begin{bmatrix} G_1 \: G_2 \\ G_3 \: G_4 \end{bmatrix}
			~			
			\quad \mathrm{and} \quad J = |\mathbf{G}| \\
			u &= J^{-1} \hat{u}, \quad \mathbf{F} = J^{-1} \mathbf{G} \hb{F}, \quad \nabla\cdot\mathbf{F} = J^{-1}\hat{\nabla}\cdot\hb{F} 
		\end{align}
		where we use $\hnabla$ to mean $[\px{}{\xi},\px{}{\eta}]^T$ in 2D. From the work of Huynh~\cite{Huynh2007}, Castonguay~\cite{Castonguay2012}, and Sheshadri~\etal~\cite{Sheshadri2016a}, Eq.(\ref{eq:2d_lin_ad}) is written in two dimensions as:
		\begin{equation}
			\hnabla\cdot\hb{F}^{\delta D} = \sum^p_{i=0}\sum^p_{j=0} \hfD_{i,j}\px{l_i(\xi)}{\xi}l_j(\eta) + \sum^p_{i=0}\sum^p_{j=0} \hgD_{i,j}\px{l_j(\eta)}{\eta}l_i(\xi) \\ \label{eq:dis_div_2d}
			\end{equation}
		\begin{equation}
			\begin{split}	
				\hnabla\cdot\hb{F}^{\delta C} = \sum^p_{i=0}\Bigg(& 
			  		(\hfI_{L,i} - \hfD_{L,i})\dx{h_{L,i}}{\xi} +
			  		(\hfI_{R,i} - \hfD_{R,i})\dx{h_{R,i}}{\xi} + \\ \nonumber
			  		&(\hgI_{B,i} - \hgD_{B,i})\dx{h_{B,i}}{\eta} +
			  		(\hgI_{T,i} - \hgD_{T,i})\dx{h_{T,i}}{\eta}\Bigg)
			\end{split}		
		\end{equation}		
		where we use $L$, $R$, $B$, and $T$ subscripts to mean left, right, bottom, and top respectively.  We will then impose that grid transformations are purely rectilinear, \emph{i.e.} $G_2 = G_3 = 0$. This is to reduce the number of dependent variables will still allowing an important form of grid deformation to be investigated.  We may now use Eq.~(\ref{eq:dis_div_2d}) and convert it into a matrix form:
		\begin{align}
			\hnabla\cdot\hb{F}^{\delta D} & = \mathbf{D}_{\xi}\hbd{f}_{i,j} + \mathbf{D}_{\eta}\hbd{g}_{i,j} \\
			\nabla\cdot\mathbf{F}^{\delta D} & = G^{-1}_{1,i,j}\mathbf{D}_{\xi}\bd{f}_{i,j} + G^{-1}_{4,i,j}\mathbf{D}_{\eta}\bd{g}_{i,j}
		\end{align}
		To apply the correction function, we need to calculate the interface values around the element. For the case of generalised central/upwinding with upwinding ratio $\alpha$, the common interface fluxes may be written as:		
		\begin{alignat}{3}
			&G^{-1}_{4,i,j}\hfI_L &&= a\big(\alpha G^{-1}_{4,i-1,j}\hud_{i-1,j,R} &&+ (1-\alpha)G^{-1}_{4,i  ,j}\hud_{i  ,j,L}\big)\\
			~
			&G^{-1}_{4,i,j}\hfI_R &&= a\big(\alpha G^{-1}_{4,i,j}\hud_{i,j,R} &&+ (1-\alpha)G^{-1}_{4,i+1,j}\hud_{i+1,j,L}\big)\\
			~			
			&G^{-1}_{1,i,j}\hgI_B &&= b\big(\alpha G^{-1}_{1,i,j-1}\hud_{i,j-1,T} &&+ (1-\alpha)G^{-1}_{1,i  ,j}\hud_{i  ,j,B} \big)\\
			~
			&G^{-1}_{1,i,j}\hgI_T &&= b\big(\alpha G^{-1}_{1,i,j}\hud_{i,j,T} &&+ (1-\alpha)G^{-1}_{1,i,j+1}\hud_{i,j+1,B}\big)
		\end{alignat} 
		where $\alpha = 1$ gives rise to upwinding and $\alpha =0.5$ produces central difference. Hence the divergence correction can be written as:
		\begin{equation}
		\begin{split}
			\hnabla\cdot\hb{F}^{\delta C}_{i,j} =			
			a\alpha\bigg(G^{-1}_{4,i-1,j}\mathbf{h_L}\mathbf{l_R}^T\bhud_{i-1,j} - G^{-1}_{4,i,j}\mathbf{h_L}\mathbf{l_L}^T\bhud_{i,j}\bigg)& + \\			
		    a(1-\alpha)\bigg(G^{-1}_{4,i+1,j}\mathbf{h_R}\mathbf{l_L}^T\bhud_{i+1,j} - G^{-1}_{4,i,j}\mathbf{h_R}\mathbf{l_R}^T\bhud_{i,j}\bigg)& + \\
			b\alpha\bigg(G^{-1}_{1,i,j-1}\mathbf{h_B}\mathbf{l_T}^T\bhud_{i,j-1} - G^{-1}_{1,i,j}\mathbf{h_B}\mathbf{l_B}^T\bhud_{i,j}\bigg)& + \\															b(1-\alpha)\bigg(G^{-1}_{1,i,j+1}\mathbf{h_T}\mathbf{l_B}^T\bhud_{i,j+1} - G^{-1}_{1,i,j}\mathbf{h_T}\mathbf{l_T}^T\bhud_{i,j}\bigg)&
		\end{split} 
		\end{equation}
		where $\mathbf{h_L}$ is the \emph{gradient} of the left correction function at the solution points and again $\mathbf{l_L}$ are the values of the polynomial basis at the left interface and so on for $R$, $T$, and $B$. Therefore, by grouping terms by their cell indexing and transforming each term into the physical domain:
		\begin{equation}
		\begin{split}
			\px{\mathbf{u}_{i,j}}{t} = &
			- a\Big ( G_{1,i-1,j}^{-1} \mathbf{C}_L \bd{u}_{i-1,j} 
			+ G_{1,i,j}^{-1} \mathbf{C}_{0\xi} \bd{u}_{i,j}
			+ G_{1,i+1,j}^{-1} \mathbf{C}_R \bd{u}_{i+1,j} \Big ) \\
			&- b\Big( G_{4,i,j-1}^{-1} \mathbf{C}_B \bd{u}_{i,j-1} 
			+ G_{4,i,j}^{-1} \mathbf{C}_{0\eta} \bd{u}_{i,j}
			+ G_{4,i,j+1}^{-1} \mathbf{C}_T \bd{u}_{i,j+1} \Big )	
		\end{split}\label{eq:2d_linear_matrix}
		\end{equation}
		where
		\begin{alignat}{4}\label{eq:Cxmat}
			&\mathbf{C}_L &&= \alpha\mathbf{h_L}\mathbf{l_R}^T \quad
			\mathbf{C}_R &&= (1-\alpha)\mathbf{h_R}\mathbf{l_L}^T \quad
			\mathbf{C}_{0\xi} &&= \mathbf{D}_{\xi} - \alpha\mathbf{h_L}\mathbf{l_L}^T - (1-\alpha)\mathbf{h_R}\mathbf{l_R}^T\\
			&\mathbf{C}_B &&= \alpha\mathbf{h_B}\mathbf{l_T}^T\quad
			\mathbf{C}_T &&= (1-\alpha)\mathbf{h_T}\mathbf{l_B}^T\quad
			\mathbf{C}_{0\eta} &&= \mathbf{D}_{\eta} - \alpha\mathbf{h_B}\mathbf{l_B}^T - (1-\alpha)\mathbf{h_T}\mathbf{l_T}^T\label{eq:Cymat}
		\end{alignat}
	Finally, we are here interested in the frequency response of the system, and, importantly to engineers and technicians, how the cell's orientation relative to an oncoming wave affects performance. Therefore, we impose a trial solution of the form: 
	\begin{equation}\label{eq:2d_bloch}
		u(x,y;t) = \exp{(ik(x\cos\theta + y\sin\theta - ct))}
	\end{equation}
	and by substitution into Eq.(\ref{eq:2d_lin_ad}), the advection velocity, $\mathbf{a}$, can be found, which is shown schematically in Fig.~\ref{fig:2D_schematic}.
	\begin{equation}
		\mathbf{a} = \begin{bmatrix} a\\ b \end{bmatrix} = \begin{bmatrix} \cos\theta \\ \sin\theta \end{bmatrix}
	\end{equation}

	The plane wave can then be projected into the computational domain and discretised as:
	\begin{equation}\label{eq:2d_project_bloch}
		\mathbf{u}_{i,j} = \mathbf{v}\exp{\bigg(ik\Big(\big(0.5(\xi +1)\delta_{i} + x_i\big)\cos\theta + \big(0.5(\eta +1)\delta_{j} + y_i\big)\sin\theta -ct\Big)\bigg)} 
	\end{equation}
	where, for brevity, $\delta_i = x_i-x_{i-1}$ and $\delta_j = y_j - y_{j-1}$ are defined. Inserting Eq.(\ref{eq:2d_project_bloch}) into Eq.(\ref{eq:2d_linear_matrix}), an Eigenvalue problem can be obtained as: 
	\begin{equation}\label{eq:2deigen}
	\begin{split}
		-ikc(k)\mathbf{v} = &
			- \cos\theta\Big(G_{1,i-1,j}^{-1} \mathbf{C}_L\exp{\big(-ik\delta_{i-1}\cos\theta\big)}
			+                G_{1,i,j}^{-1} \mathbf{C}_{0\xi}
			+                G_{1,i+1,j}^{-1} \mathbf{C}_R\exp{\big(ik\delta_{i}\cos\theta\big)} \Big)\mathbf{v} \\
			&-\sin\theta\Big(G_{4,i,j-1}^{-1} \mathbf{C}_B\exp{\big(-ik\delta_{j-1}\sin\theta\big)} 
			+                G_{4,i,j}^{-1} \mathbf{C}_{0\eta} 
			+                G_{4,i,j+1}^{-1} \mathbf{C}_T\exp{\big(ik\delta_j\sin\theta\big)} \Big)\mathbf{v}
	\end{split}
	\end{equation}
	where $\Re{k(c(k))}=\Re{(\omega)}$ and $\Im{(kc(k))}=\Im{(\omega)}$ are the dispersion and dissipation, respectively, and $\omega$ is the modified angular frequency response of the system. BY studying the trial solution of Eq.(\ref{eq:2d_bloch}) it can be understood that if $\Im{(\omega)} > 0$ then the amplitude of the wave will increase and vice versa. Furthermore, if $\Re{(\omega)}\ne k$ then a wave will move at a different speed compared to the other waves inside the packet, causing the quality of the interpolation to be affected as the solution is advanced in time as the components that make it up move at different speeds. An important point is the difference between phase velocity $\omega/k$ and group velocity $\mathrm{d}\omega/\mathrm{d}k$. Phase velocity is the speed of a wave in a packet of waves. Group velocity is the speed of the packet. Therefore, changes to $\mathrm{d}\omega/\mathrm{d}k$, i.e can be thought of more as a change to the physics due to the numerical method.
	
    Equation~(\ref{eq:2deigen}) can alternatively be cast in the form of an update equation. If initially Eq.(\ref{eq:2d_linear_matrix}) is combined with Eq.(\ref{eq:2d_project_bloch}), then a new matrix, $\mathbf{Q}_{i,j}$, can be defined:
	\begin{equation}
		\px{\mathbf{u}_{i,j}}{t} = \mathbf{Q}_{i,j}\mathbf{u}_{i,j}		
	\end{equation}
	\begin{equation}\label{eq:Q_matrix}
		\begin{split}
		\mathbf{Q}_{i,j}  =&  
			- J_{i,j}\cos\theta\Big(G_{4,i-1,j}^{-1} \mathbf{C}_L \exp{\big(-ik\delta_{i-1}\cos\theta\big)}
			+                G_{4,i,j}^{-1} \mathbf{C}_{0\xi}
			+                G_{4,i+1,j}^{-1} \mathbf{C}_R\exp{\big(ik\delta_{i}\cos\theta\big)} \Big) \\
			& - J_{i,j}\sin\theta\Big(G_{1,i,j-1}^{-1} \mathbf{C}_B \exp{\big(-ik\delta_{j-1}\sin\theta\big)} 
			+                G_{1,i,j}^{-1} \mathbf{C}_{0\eta} 
			+                G_{1,i,j+1}^{-1} \mathbf{C}_T \exp{\big(ik\delta_j\sin\theta\big)} \Big)
		\end{split}	
	\end{equation}
	This definition of the semi-discrete FR operator, $\mathbf{Q}$, can then be used to form what is called the update equation by imposing some temporal discretisation. As such we may write:
	\begin{align}
		\mathbf{u}^{n+1}_{i,j} &= \mathbf{R}(\mathbf{Q}_{i,j})\mathbf{u}^n_{i,j} \label{eq:update_eq} \\
		\mathbf{R}_{33} &= \mathbf{I} + \frac{\tau\mathbf{Q}_{i,j}}{1!} + \frac{(\tau\mathbf{Q}_{i,j})^2}{2!} + \frac{(\tau\mathbf{Q}_{i,j})^3}{3!} \label{eq:RK33}  	
	\end{align}
	where the superscript denotes the time level, and our update matrix is $\mathbf{R}$. Shown here is also an example definition for $\mathbf{R}$ for a 3-step $3^{\mathrm{rd}}$-order Runge-Kutta time integration scheme. Finally, in keeping with von Neumann's theorems \cite{Isaacson1994,Lax1956} and Banach's fixed point theorem~\cite{Kress1998}, the spectral radius of $\mathbf{R}$ has to be less than or equal to 1 for stability. $\rho(\mathbf{R}) \leqslant 1 \: \forall \: k \in \mathbb{R}$.		
		
	In recent works by Vermerie~\etal~\cite{Vermeire2017a} and Trojak~\etal~\cite{Trojak2018c}, the Fourier analysis was extended by fully discretising the equation. This is performed by taking Eq.(\ref{eq:update_eq}) and again applying Eq.(\ref{eq:2d_bloch}). This results in:
	\begin{equation}
		\exp{(-ik(c-1)\tau)}\mathbf{v} = \lambda\mathbf{v} = \exp{(ik\tau)}\mathbf{R}(k,\tau)\mathbf{v}
	\end{equation}
	where the time step from $n$ to $n+1$ is $\tau$. Hence, rearranging for the modified wave speed:
	\begin{equation}
		c = \frac{i\log{(\lambda)}}{k\tau} + 1
	\end{equation}
	where $\lambda$ are the eigenvalues of $\exp{(ik\tau)}\mathbf{R}$. The advantage of this further analysis is that it gives the dispersion and dissipation relations of the full scheme as would be experienced when applied as implicit LES.
	
	Lastly, the linear FR  operator matrix, $\mathbf{Q}_{i,j}$, can be diagonalised as:
	\begin{equation}
		\mathbf{Q}_{i,j} = ik\mathbf{W}_{i,j}\mathbf{\Lambda}_{i,j}\mathbf{W}^{-1}_{i,j}
	\end{equation}
	where $\mathbf{W}$ is a matrix of eigenvectors and $\mathbf{\Lambda}$ is a diagonal matrix of normalised eigenvalues. From this, the weights of modes used to reconstruct the solution can be found as:
	\begin{equation}
		\mathbf{u}_{i,j} = \mathbf{W}_{i,j}\pmb{\beta}_{i,j}
	\end{equation}
	where $\pmb{\beta}$ is an array of mode weights used to project $\mathbf{u}_{i,j}$ into the functional space of FR. The posedness of the projection can then be measured for wavenumbers using the number of the matrix of modes defined as:
	\begin{equation}
		\kappa(\mathbf{W}_{i,j}) = \frac{\sigma_\mathrm{max}(\mathbf{W}_{i,j})}{\sigma_\mathrm{min}(\mathbf{W}_{i,j})}
	\end{equation}
	where $\sigma(\mathbf{W}_{i,j})$ is a singular value of $\mathbf{W}_{i,j}$, with the matrix becoming singular as $\kappa \rightarrow \infty$.
	
	The results of this section can then be extended to n-dimensions and, in particular, the three--dimensional case will be investigated to show the continuation of trends with higher dimensionality. The analysis can broadly be repeated and is excluded for brevity, but importantly the prescribed solution is taken as:
	\begin{equation}
		u = \exp{\big(ik(x\cos\phi\cos\theta + y\cos\phi\sin\theta + z\sin\phi - ct)\big)}
	\end{equation}
	where the angles are as shown in Fig.~\ref{fig:3D_schematic}, and hence the 3D convective velocities for linear advection are:
	\begin{equation}
		\mathbf{a} = \begin{bmatrix} \cos\phi\cos\theta \\ \cos\phi\sin\theta \\ \sin\phi\end{bmatrix}
	\end{equation}


\section{Analytical Findings}\label{sec:results_analytical}
	The analytical methods presented in section~\ref{sec:vn_method} allow us to investigate many properties of FR, however from Eq.(\ref{eq:Q_matrix}-\ref{eq:RK33}) it can be seen that the functional space of $\mathbf{Q}$ is 8 dimensional, leading to the functional space of $\rho(\mathbf{R})$ being 9 dimensional ($\tau,\gamma_x,\gamma_y,\Delta_x,\Delta_y,k,\theta,\iota,p$). Therefore we need to restrict our investigation to some key results relating to grid deformation. Firstly, understanding the dispersion and dissipation ($\Re{(\omega)}$ \& $\Im{(\omega)}$) in 2D for both uniform and stretched grids will be important. Secondly, we wish to understand how higher dimensionality and grid deformations affect the temporal stability of FR through evaluation of the CFL limits~\cite{Courant1967}. Here, the dispersion and dissipation relations will be useful in explaining the trends seen and will aid in linking this work to that of Trojak~\etal~\cite{Trojak2017a}. The definition of the CFL number in higher dimensions will be taken as:
	\begin{equation}
		\mathrm{CFL}_d = \tau\sum^d_{i=1}\frac{a_i}{\Delta_i}
	\end{equation}
	where $d$ is the dimensionality, $\tau$ is the time step and, $a_i$ and $\Delta_i$ are the characteristic velocity and grid spacing in the $i^{\mathrm{th}}$ dimension, respectively. The CFL limit is then the maximum value of CFL at which the scheme is stable in a von Neumann sense. Finally, we wish to understand if correction functions can be used to alleviate any effects of deformation by understanding how the scheme properties vary with correction function. Within this study, the link to the posed nature of the linear system will be explored with regard to how this relates to the other properties.
	
	\subsection{Review of 1D Grid Expansion}\label{sec:1dvn}
\begin{figure}
		\centering
		\begin{subfigure}[b]{0.45\linewidth}
			\centering
			\includegraphics[width=\linewidth]{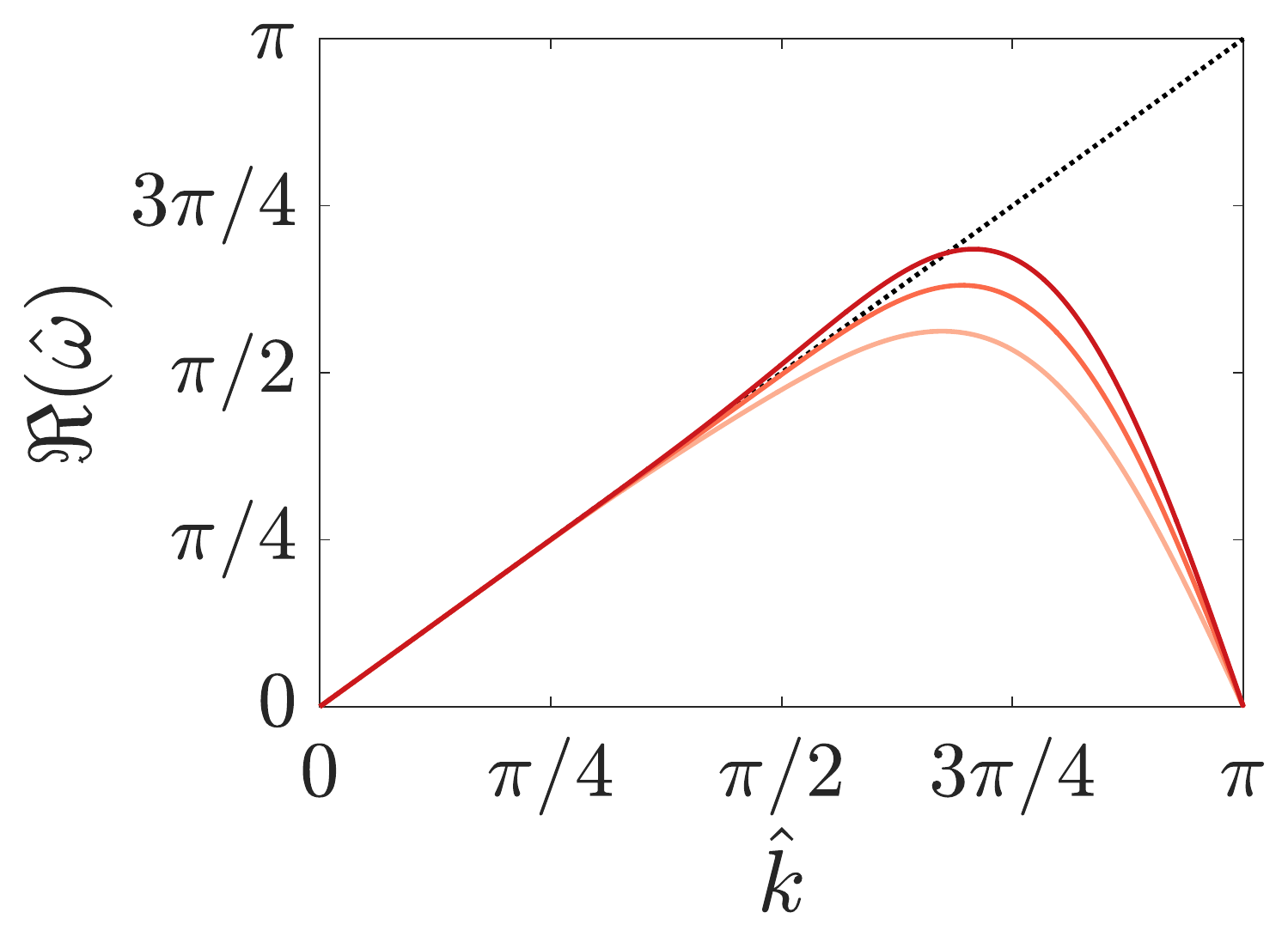}
			\caption{Dispersion}
			\label{fig:1D_FR4_disp}
		\end{subfigure}
		~
		\begin{subfigure}[b]{0.45\linewidth}
			\centering
			\includegraphics[width=\linewidth]{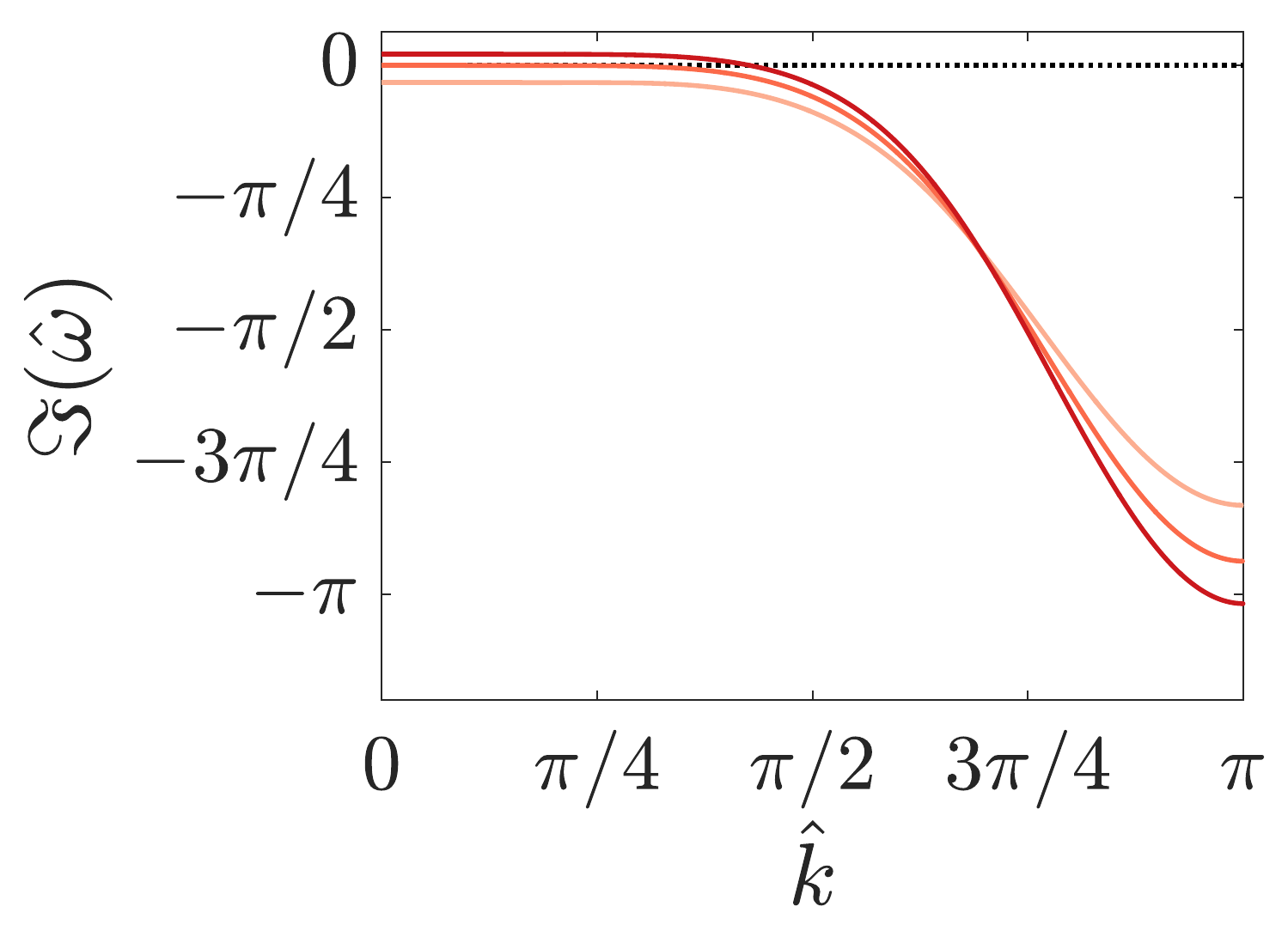}
			\caption{Dissipation}
			\label{fig:1D_FR4_diss}
		\end{subfigure}
		~
		\begin{subfigure}[b]{0.48\linewidth}
			\centering
			\includegraphics[width=\linewidth]{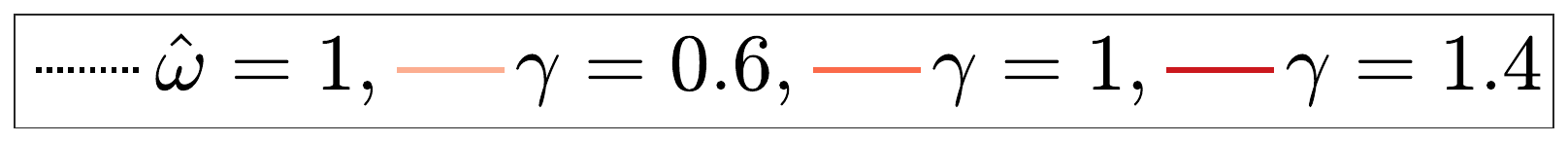}
		\end{subfigure}
		\caption{Upwinded 1D FR, $p=4$, with Huynh $g_2$ correction functions for several geometric expansion ratios.}
		\label{fig:FR1D_p4}
	\end{figure}

		Before commencing with the Fourier/von Neumann analysis in 2D, we will give a brief review of the behaviour exhibited in one dimension, for which more detail can be found in Trojak~\etal~\cite{Trojak2017a}. Figure~\ref{fig:FR1D_p4} shows the results for Huynh $g_2$ correction functions for $p=4$ on a few geometrically expanding grids. 
	
	Starting by considering the dispersion, Fig.~\ref{fig:1D_FR4_disp} shows that expanding grids cause dispersion overshoot, as is seen by the dispersion relation going above the line $\hat{\omega}=1$. While contractions cause dispersion to undershoot. More interesting is the impact of the grid on the dissipation as is shown in  Fig.~\ref{fig:1D_FR4_diss}. Here expanding grids are shown to have a region where $\Im{(\hat{\omega})} > 0$ at low wavenumbers, while contracting grids in the same wavenumber range have $\Im{(\hat{\omega})} < 0$. This implies that a wave propagating across a contracting grid is more quickly dissipated, whereas on an expanding grid the amplitude of the wave is increased. Hence, expanding grids are strictly unstable. This result was first presented in \cite{Trojak2017a} and was confirmed with numerical experiments.
	
	This characteristic is obviously of importance as it could lead to the scheme being non-conservative. Hence, the importance of studying this phenomena in higher dimensions.
	
	\subsection{Effect of Grid Expansion on Dispersion and Dissipation}\label{sec:dissdisp}
	For the higher dimensional case, we begin by considering the dispersion and dissipation on a uniform grid in two dimensions. We are concerned here with the primary mode --- as FR has multiple modes, this is the one that physically represents the wave. Although, as was found by Asthana~\cite{Asthana2017}, this may not be how the energy distributes itself. We identify the physical mode as the mode whose dispersion relation that goes through zero and dissipation relation similar to those seen in Fig.~\ref{fig:1D_FR4_diss}. 
	
\begin{figure}
		\centering
		\begin{subfigure}[b]{0.4\linewidth}
			\centering
			\includegraphics[width=\linewidth]{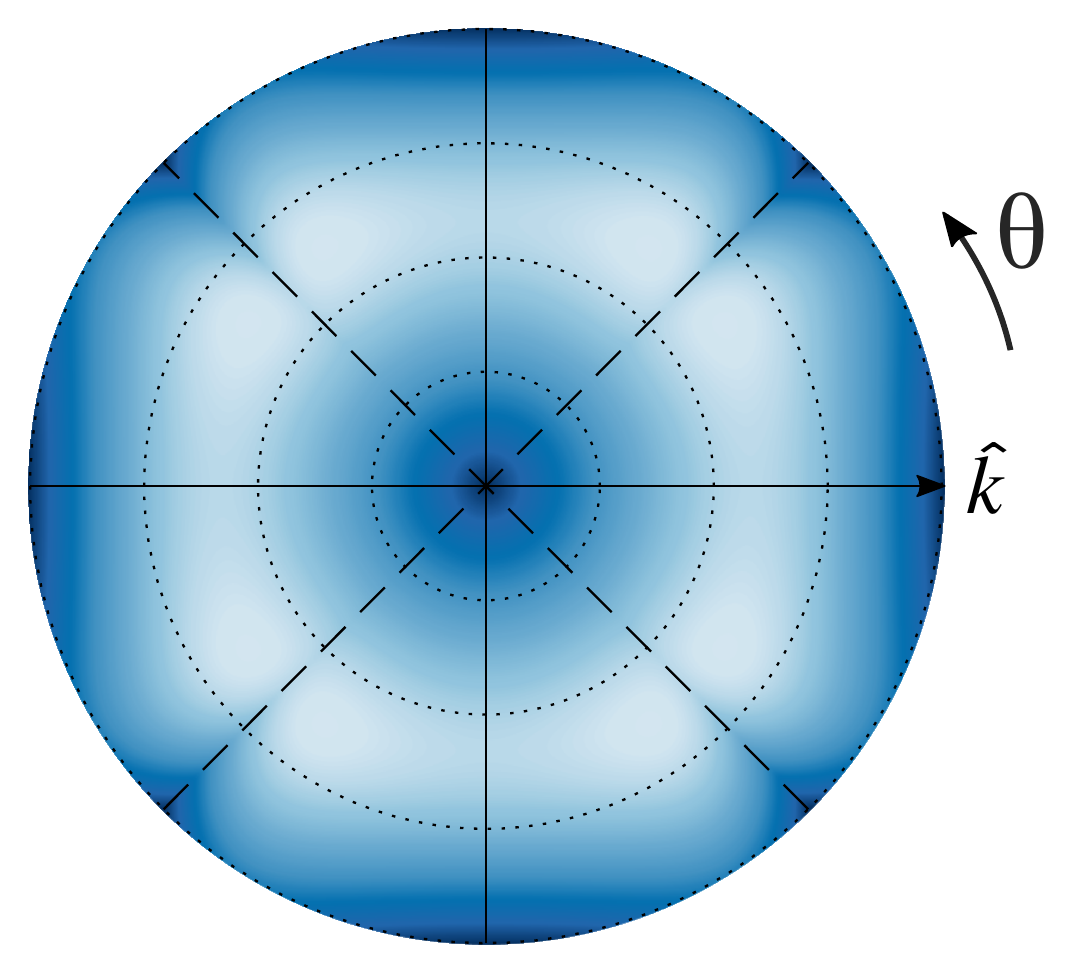}
			\caption{$p=2$}
			\label{fig:FR2_disp}
		\end{subfigure}
		~
		\begin{subfigure}[b]{0.4\linewidth}
			\centering
			\includegraphics[width=\linewidth]{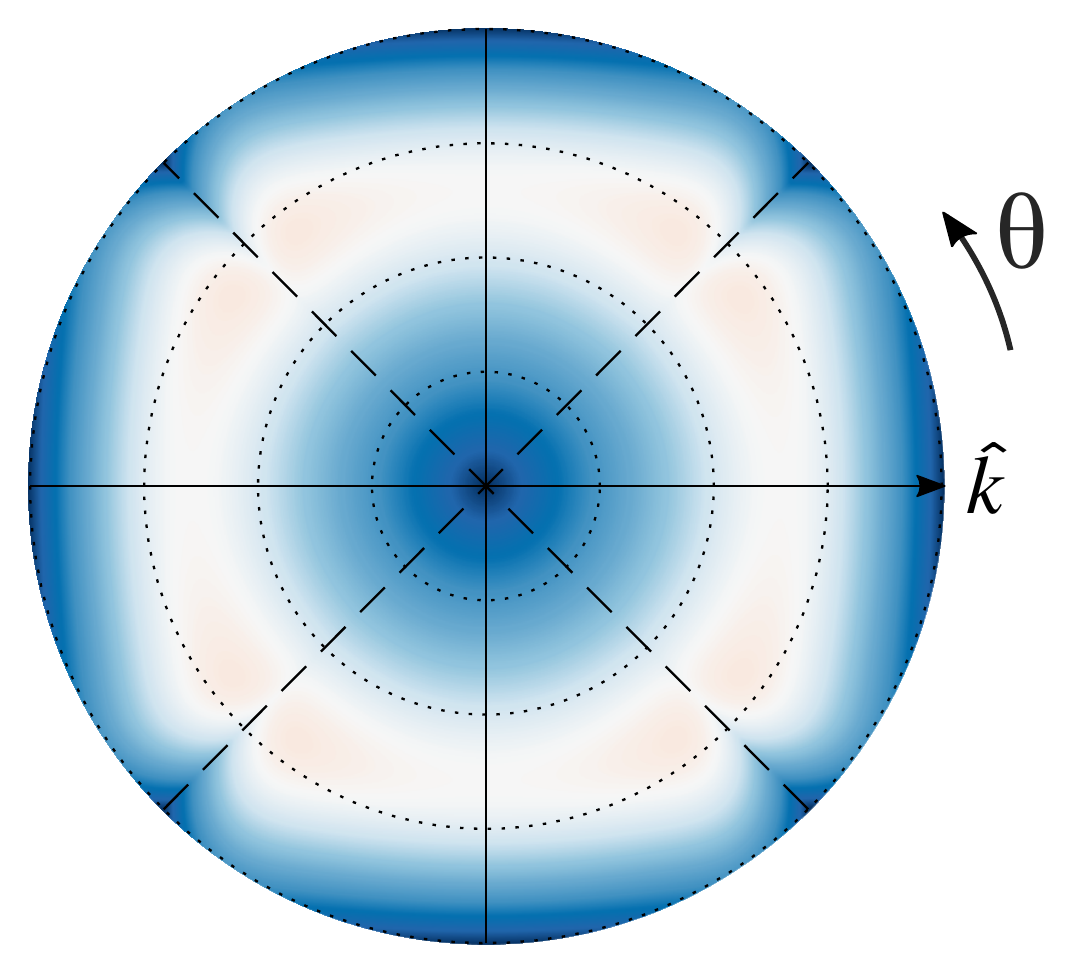}
			\caption{$p=3$}
			\label{fig:FR3_disp}
		\end{subfigure}
		~
		\begin{subfigure}[b]{0.4\linewidth}
			\centering
			\includegraphics[width=\linewidth]{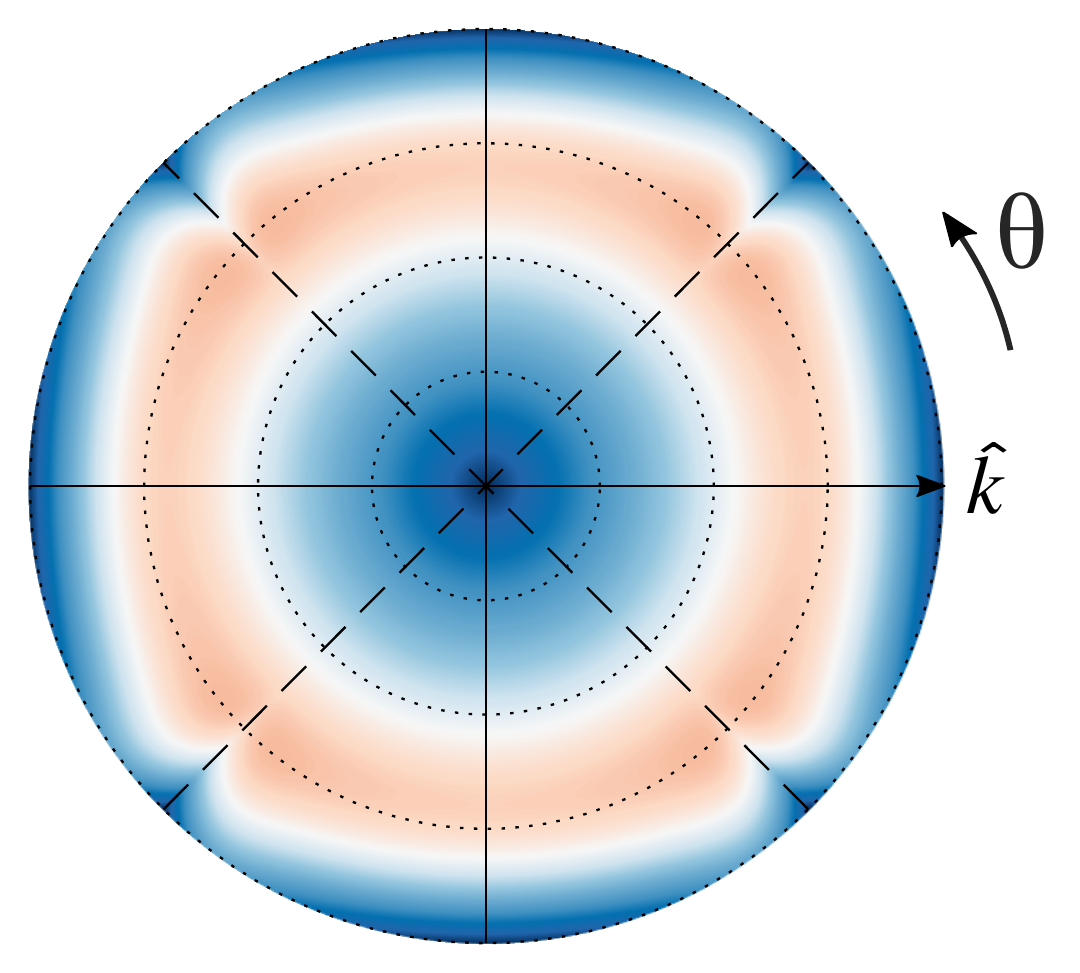}
			\caption{$p=4$}
			\label{fig:FR4_disp}
		\end{subfigure}
		~
		\begin{subfigure}[b]{0.4\linewidth}
			\centering
			\includegraphics[width=\linewidth]{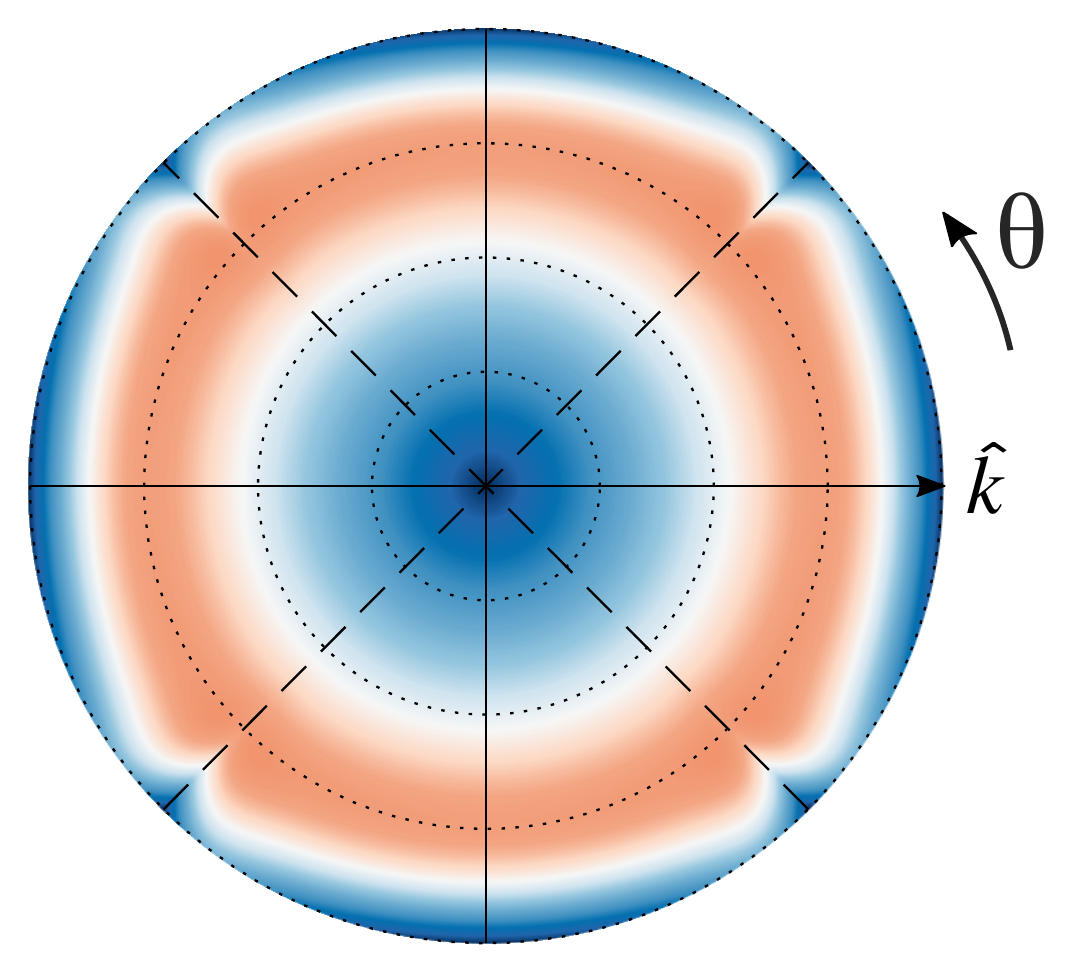}
			\caption{$p=5$}
			\label{fig:FR5_disp}
		\end{subfigure}
		~
		\begin{subfigure}[b]{0.3\linewidth}
			\centering
			\includegraphics[width=\linewidth]{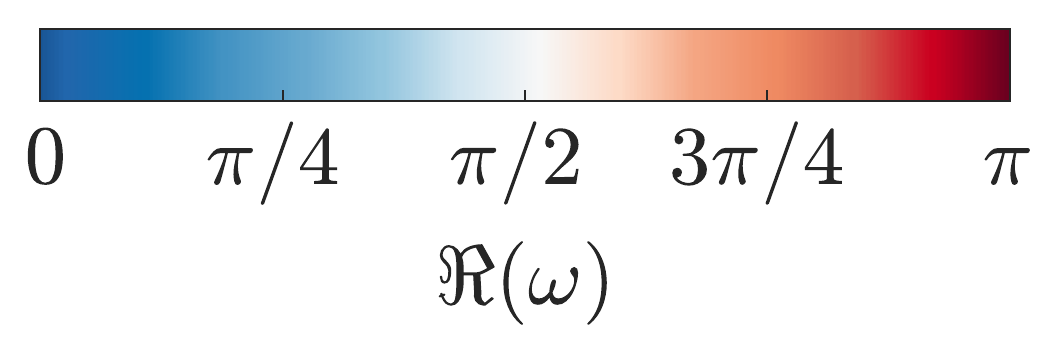}
		\end{subfigure}
		\caption{Primary mode dispersion for 2D upwinded FR, with Huynh $g_2$ corrections, at various orders.  Normalised wavenumber as radial distance (markers at $\pi/4$ intervals), and element angle of incidence as azimuthal distance.}
		\label{fig:FR_polar_disp}
	\end{figure}

\begin{figure}
		\centering
		\begin{subfigure}[b]{0.4\linewidth}
			\centering
			\includegraphics[width=\linewidth]{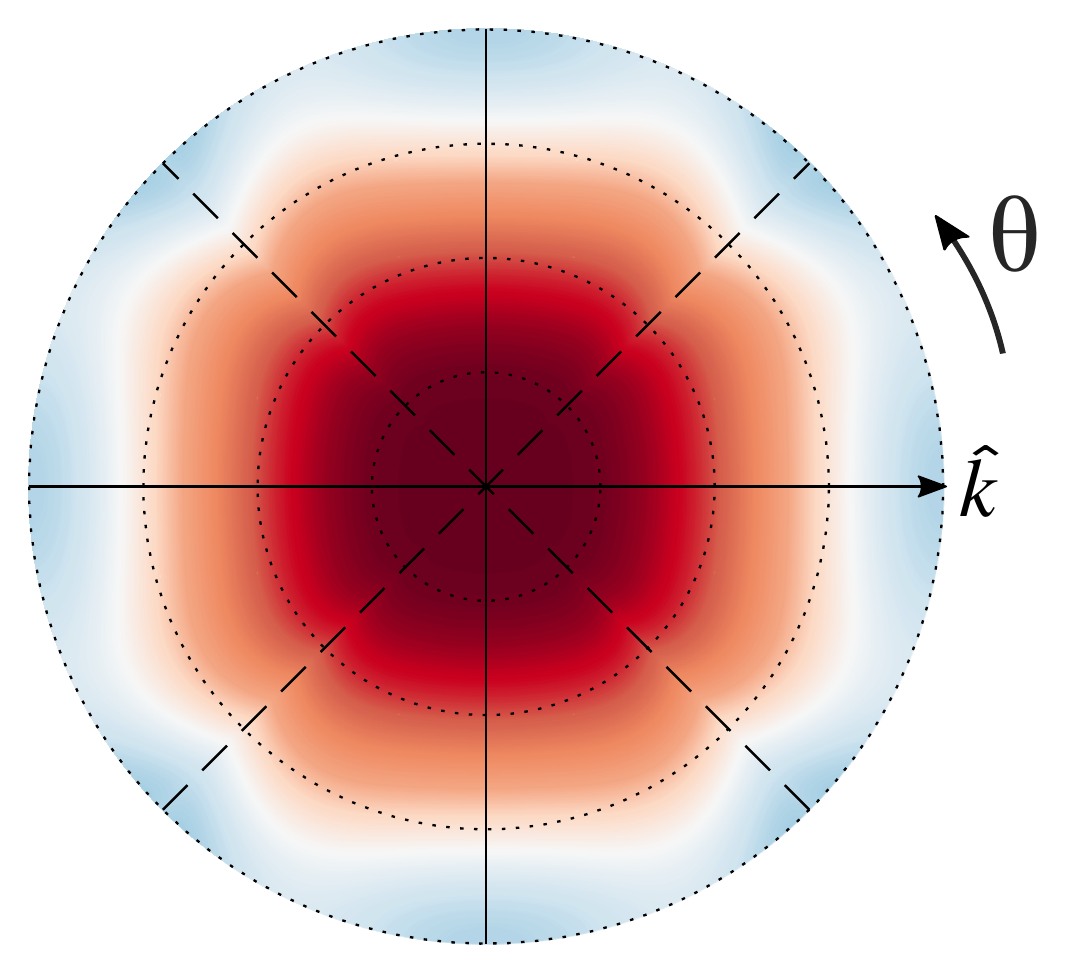}
			\caption{$p=2$}
			\label{fig:FR2_diss}
		\end{subfigure}
		~
		\begin{subfigure}[b]{0.4\linewidth}
			\centering
			\includegraphics[width=\linewidth]{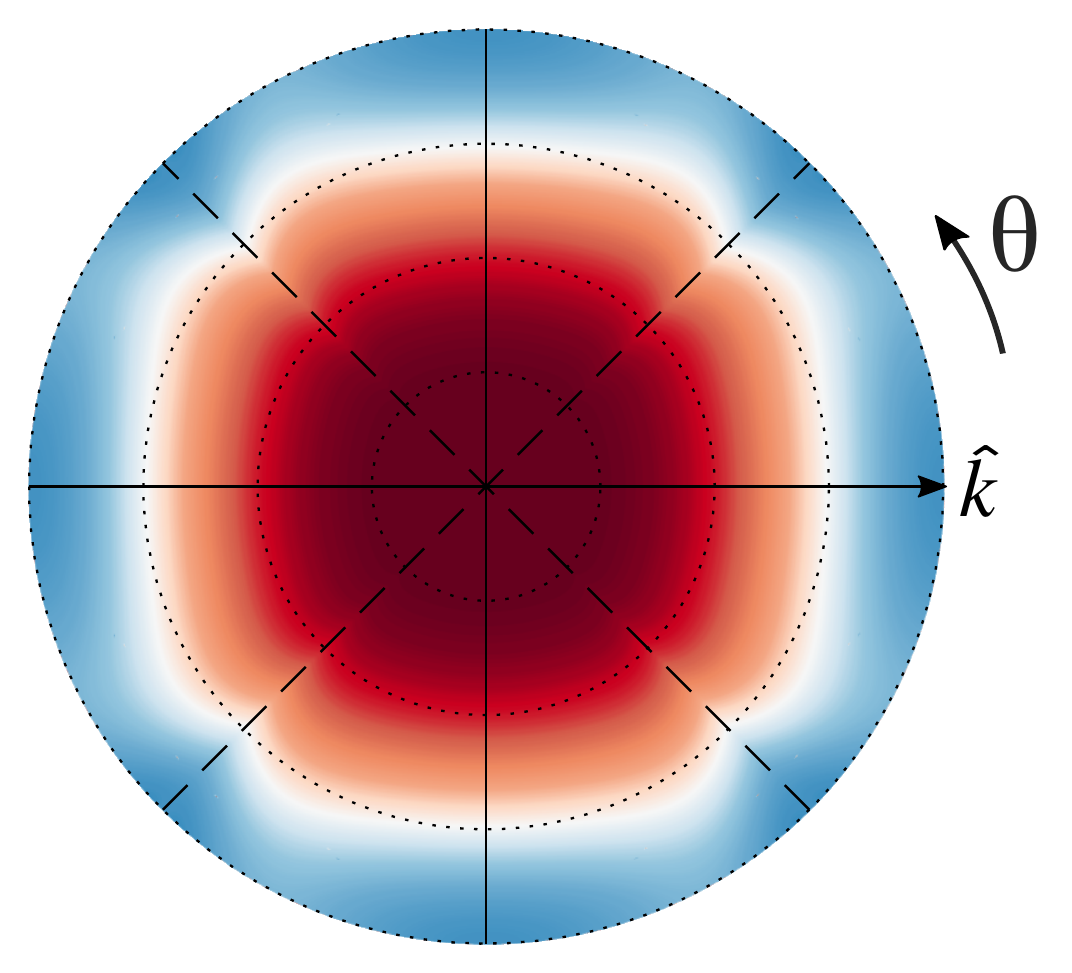}
			\caption{$p=3$}
			\label{fig:FR3_diss}
		\end{subfigure}
		~
		\begin{subfigure}[b]{0.4\linewidth}
			\centering
			\includegraphics[width=\linewidth]{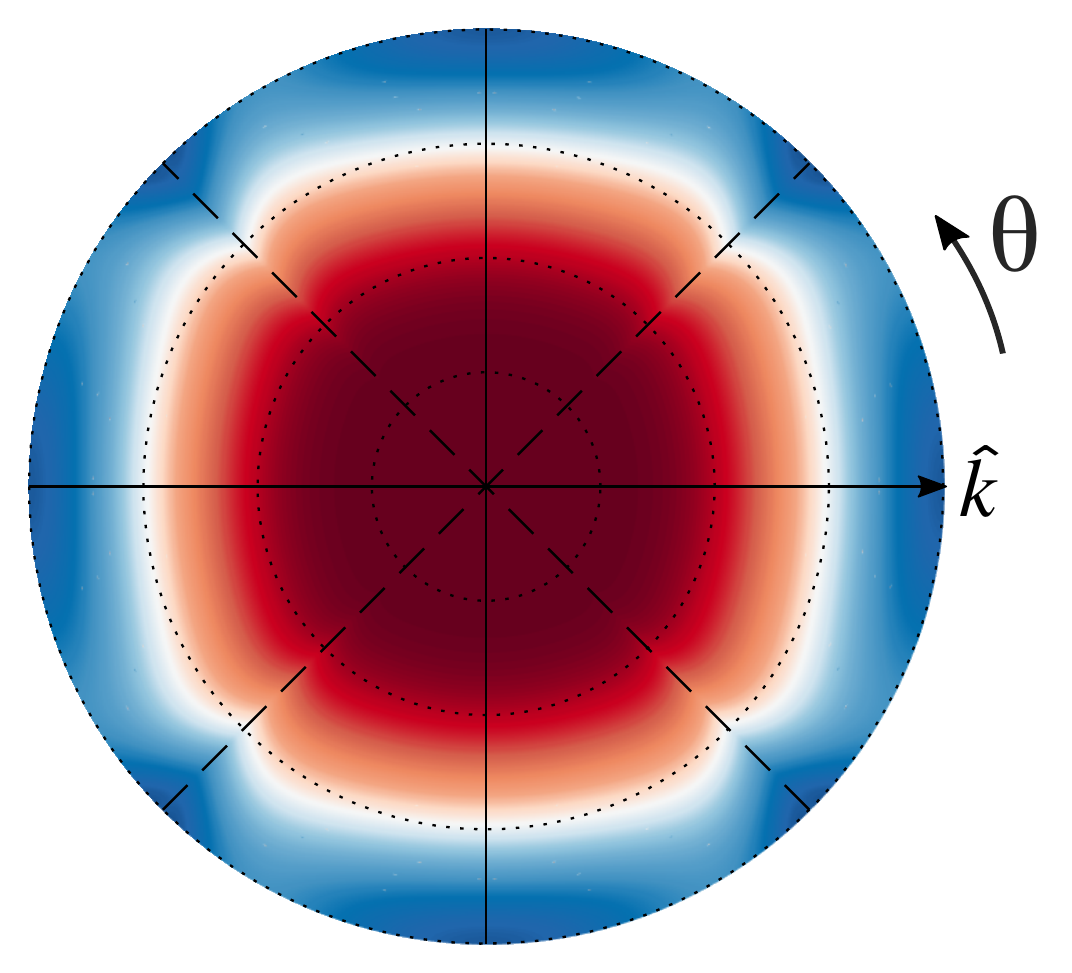}
			\caption{$p=4$}
			\label{fig:FR4_diss}
		\end{subfigure}
		~
		\begin{subfigure}[b]{0.4\linewidth}
			\centering
			\includegraphics[width=\linewidth]{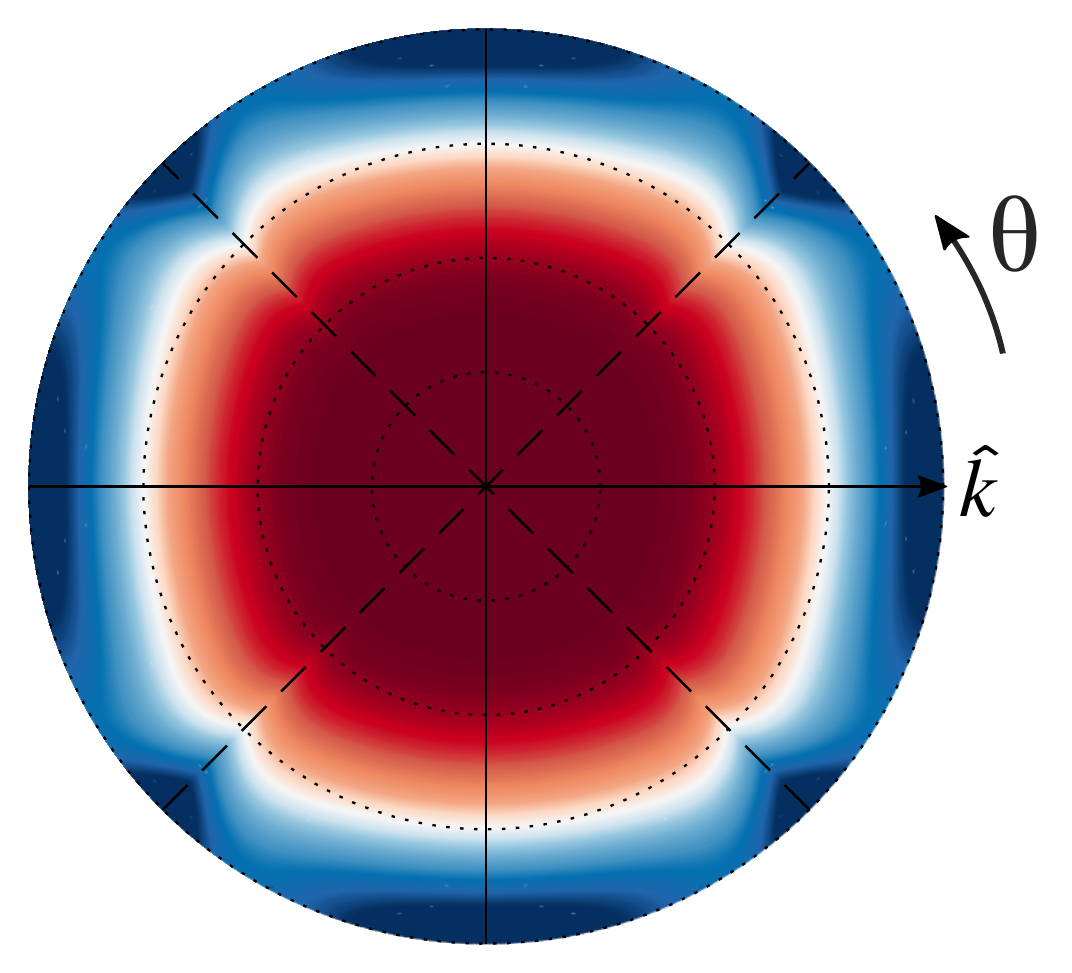}
			\caption{$p=5$}
			\label{fig:FR5_diss}
		\end{subfigure}
		~
		\begin{subfigure}[b]{0.3\linewidth}
			\centering
			\includegraphics[width=\linewidth]{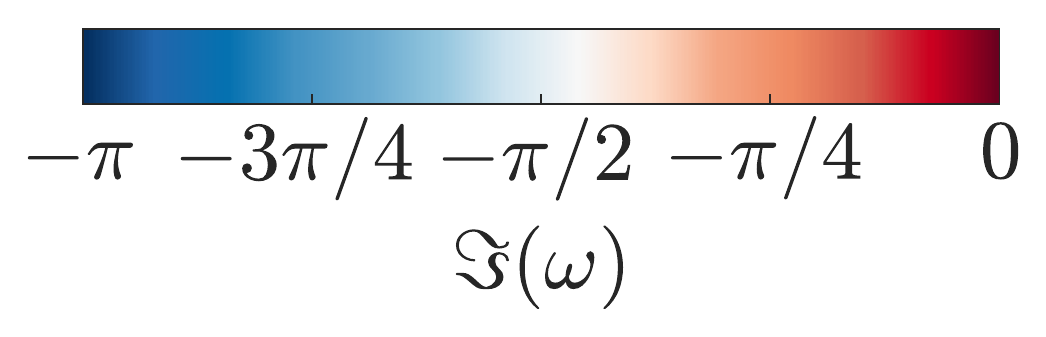}
		\end{subfigure}
		\caption{Primary mode dissipation for 2D upwinded FR, with Huynh $g_2$ corrections, at various orders.  Normalised wavenumber as radial distance (markers at $\pi/4$ intervals), and element angle of incidence as azimuthal distance.}
		\label{fig:FR_polar_diss}
	\end{figure}

	The dispersion and dissipation relations are then shown in Figs.~\ref{fig:FR_polar_disp}~and~\ref{fig:FR_polar_diss}. It is clear that for all orders FR becomes more dispersive and dissipative at $\theta=45^{\circ}$. This is more easily seen for the dissipation relation, but in the case of dispersion is displayed by the indent at $\hat{k}\approx 3\pi/4$. Furthermore, there doesn't seem to be any widening in the range of angles over which FR becomes more dispersive and dissipative as the order is reduced. By comparison with the results of Lele~\cite{Lele1992}, where a similar test is performed for standard and compact difference schemes, FR shows a comparatively smaller change in performance as the angle is varied. It is thought that this due to the method of polynomial fitting used by FR, namely that this implementation of FR used a tensor grid of monomials \emph{i.e.}, the number of solution points is $(p+1)^d$ and hence the monomials in the interpolation go from $(\xi^0\eta^0,\xi^1\eta^0\dots\xi^p\eta^p)$. By contrast, finite differences do not include the cross product terms, which will become increasingly dominant as the angle is increased. 
	
	Moving on, we then consider the impact of non-uniform grids on the character of the dispersion and dissipation. There are two cases that have been identified from previous work as being of interest. Firstly, when the grid is expanding, does this give the same positive dissipation in higher dimensions as seen in the 1D case, Fig.~\ref{fig:1D_FR4_diss}? Secondly, if positive dissipation is seen in the 2D case, will an orthogonal contraction help to stabilise the grid? For example, if $\gamma_x=1.1$ then will setting $\gamma_y=0.9$ help to reduce the positive dissipation?
	
	\begin{figure}
		\centering
		\begin{subfigure}[b]{0.4\linewidth}
			\centering
			\includegraphics[width=\linewidth]{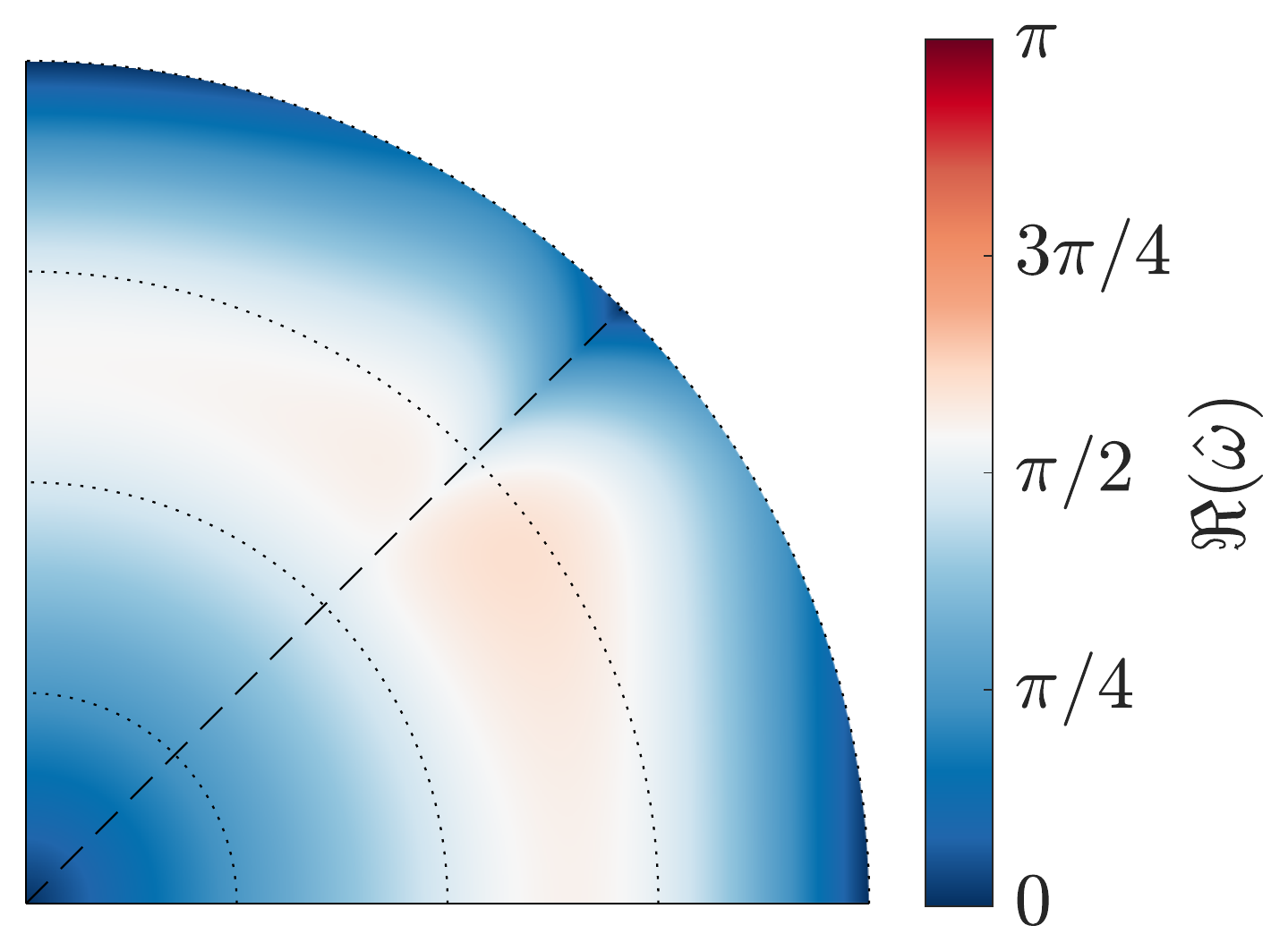}
			\caption{Dispersion, $\gamma_x=1.1$, $\gamma_y=1$}
			\label{fig:FR3_disp_11}
		\end{subfigure}
		~
		\begin{subfigure}[b]{0.4\linewidth}
			\centering
			\includegraphics[width=\linewidth]{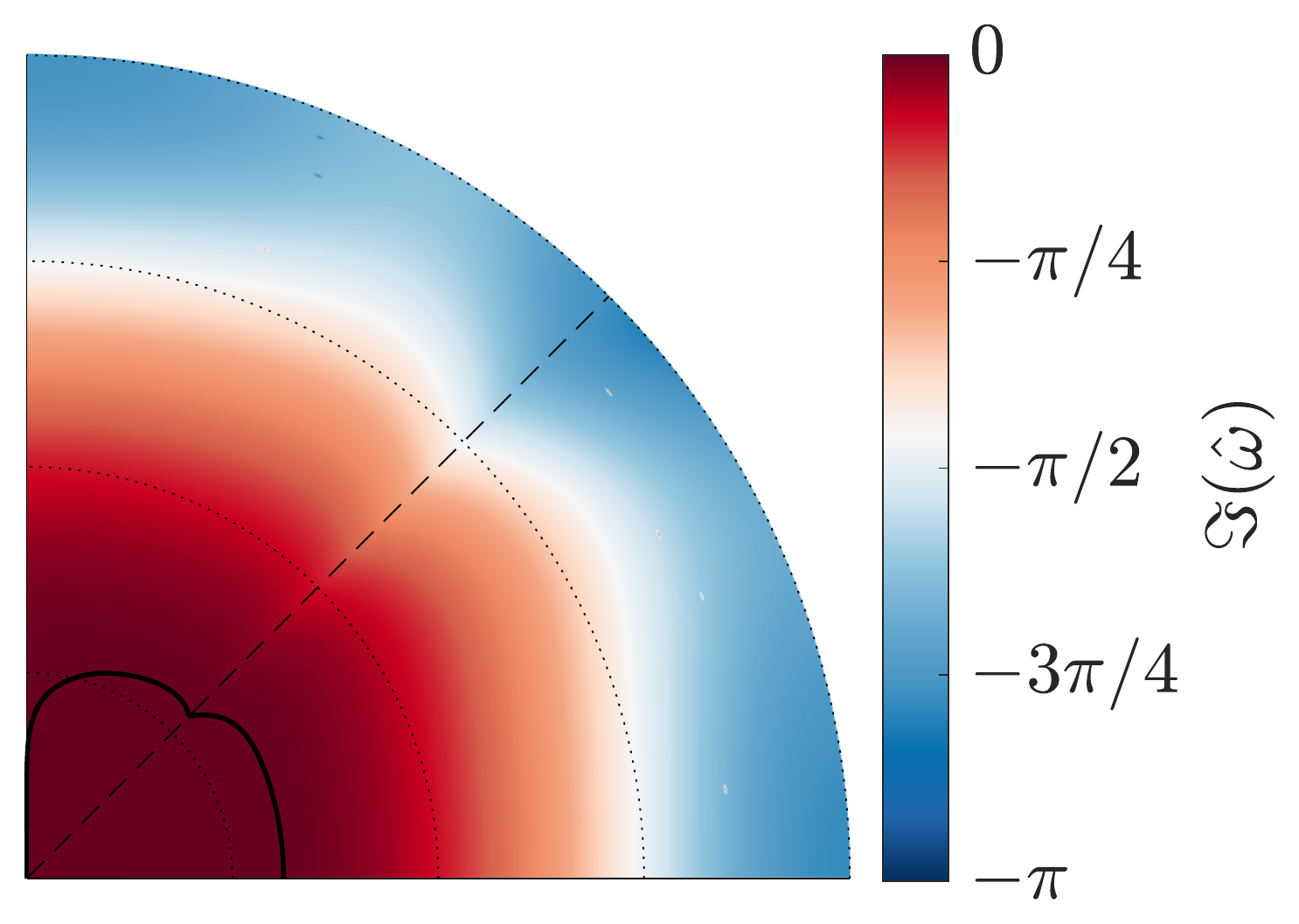}
			\caption{Dissipation, $\gamma_x=1.1$, $\gamma_y=1$}
			\label{fig:FR3_diss_11}
		\end{subfigure}
		~
		\begin{subfigure}[b]{0.4\linewidth}
			\centering
			\includegraphics[width=\linewidth]{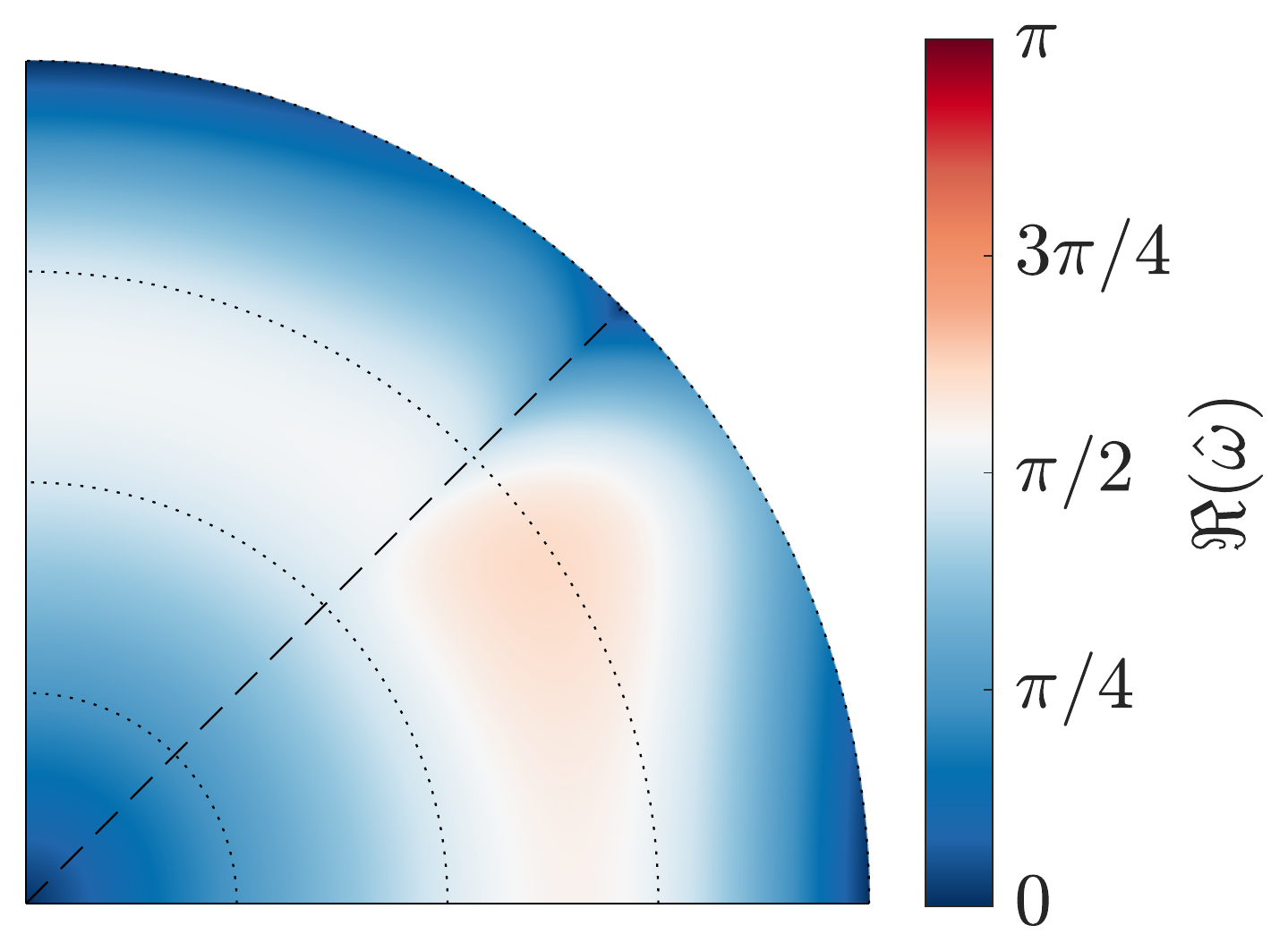}
			\caption{Dispersion, $\gamma_x=1.1$, $\gamma_y=0.9$}
			\label{fig:FR3_disp_1109}
		\end{subfigure}
		~
		\begin{subfigure}[b]{0.4\linewidth}
			\centering
			\includegraphics[width=\linewidth]{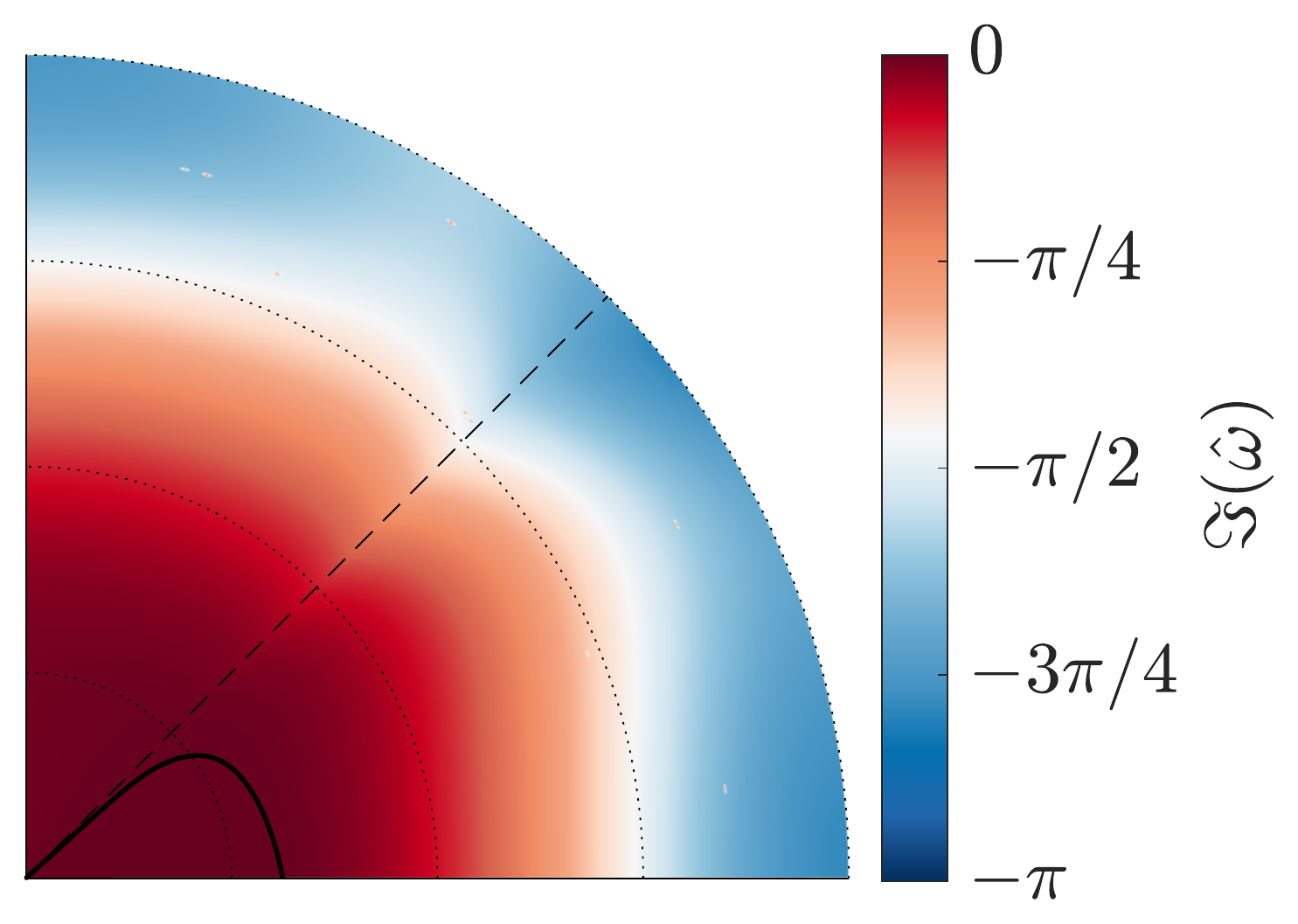}
			\caption{Dissipation, $\gamma_x=1.1$, $\gamma_y=0.9$}
			\label{fig:FR3_diss_1109}
		\end{subfigure}
		\caption{Two dimensional upwinded FR, $p=3$ with Huynh $g_2$ corrections, for different grid expansion factors. Normalised wavenumber as radial distance (markers at $\pi/4$ intervals), and element angle of incidence as azimuthal distance. The solid black line on the dissipation plots is the contour of zero dissipation.}
		\label{fig:FR_polar_disp_stretched}
	\end{figure}
	
\begin{figure}
		\centering
		\begin{subfigure}[b]{0.45\linewidth}
			\centering
			\includegraphics[width=\linewidth]{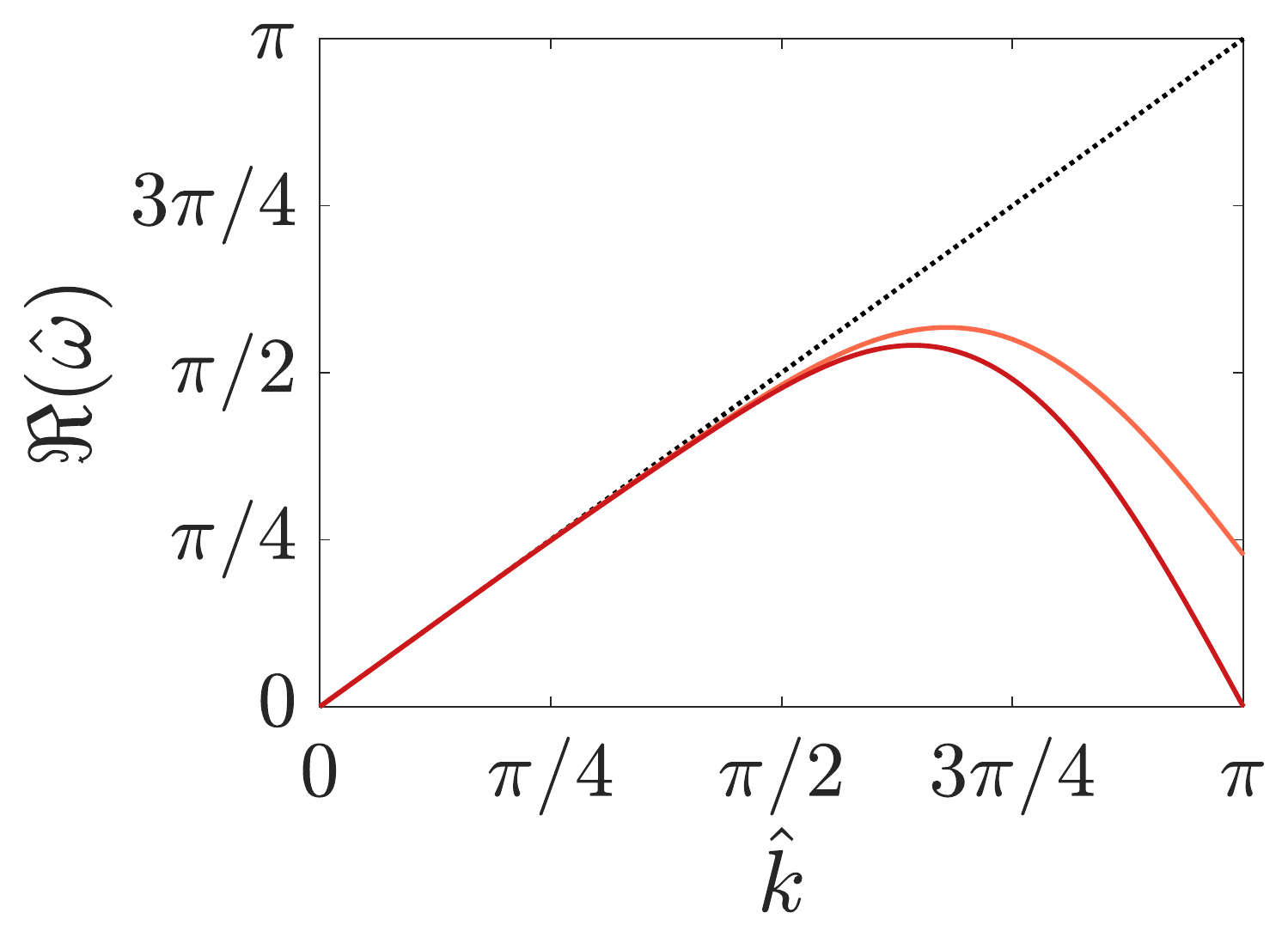}
			\caption{Dispersion: $\gamma_x=\gamma_y = 1$}
			\label{fig:FRp4_re_1010angle}
		\end{subfigure}
		~
		\begin{subfigure}[b]{0.45\linewidth}
			\centering
			\includegraphics[width=\linewidth]{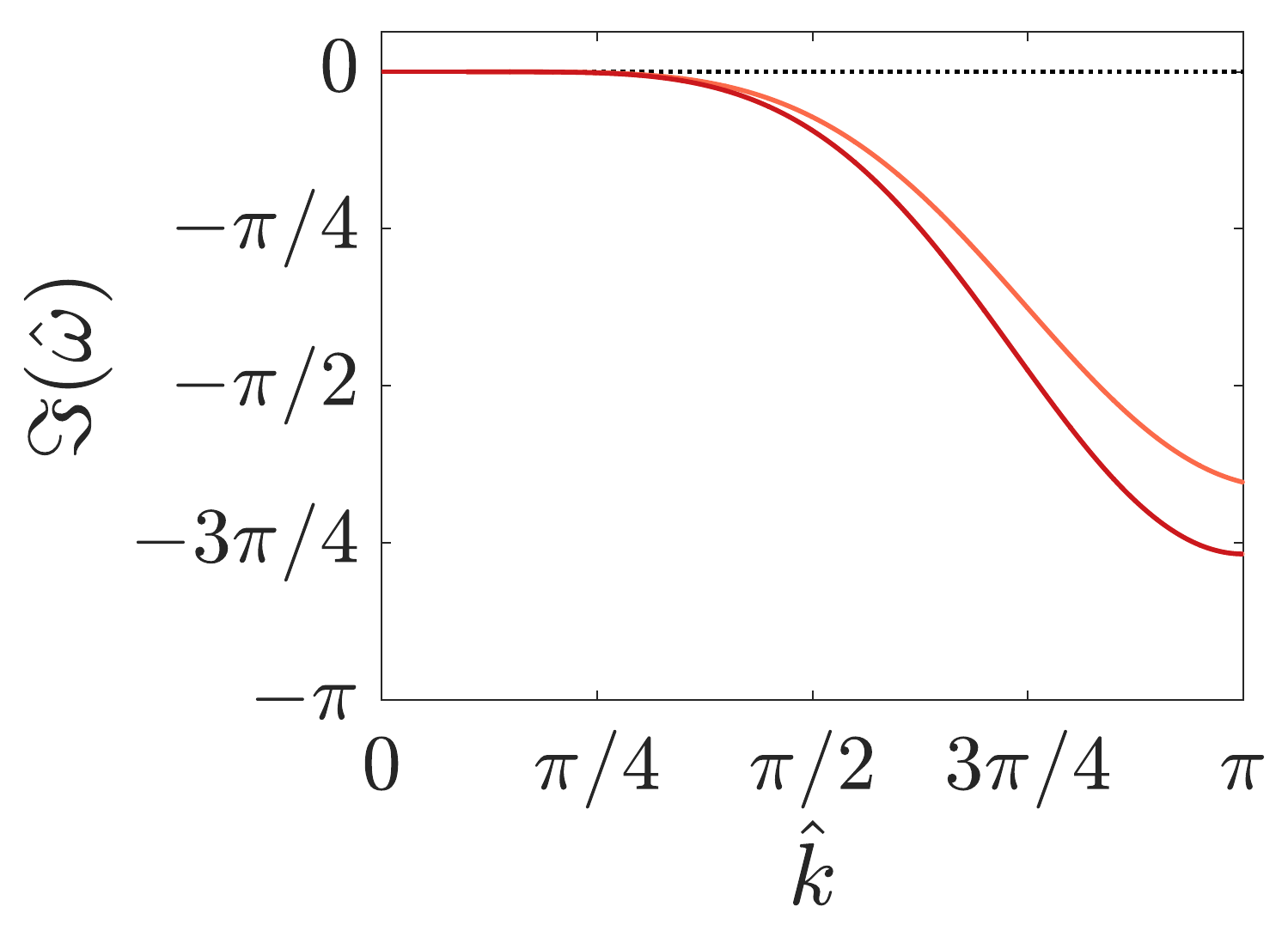}
			\caption{Dissipation: $\gamma_x=\gamma_y = 1$}
			\label{fig:FRp4_im_1010angle}
		\end{subfigure}
		~
		\begin{subfigure}[b]{0.45\linewidth}
			\centering
			\includegraphics[width=\linewidth]{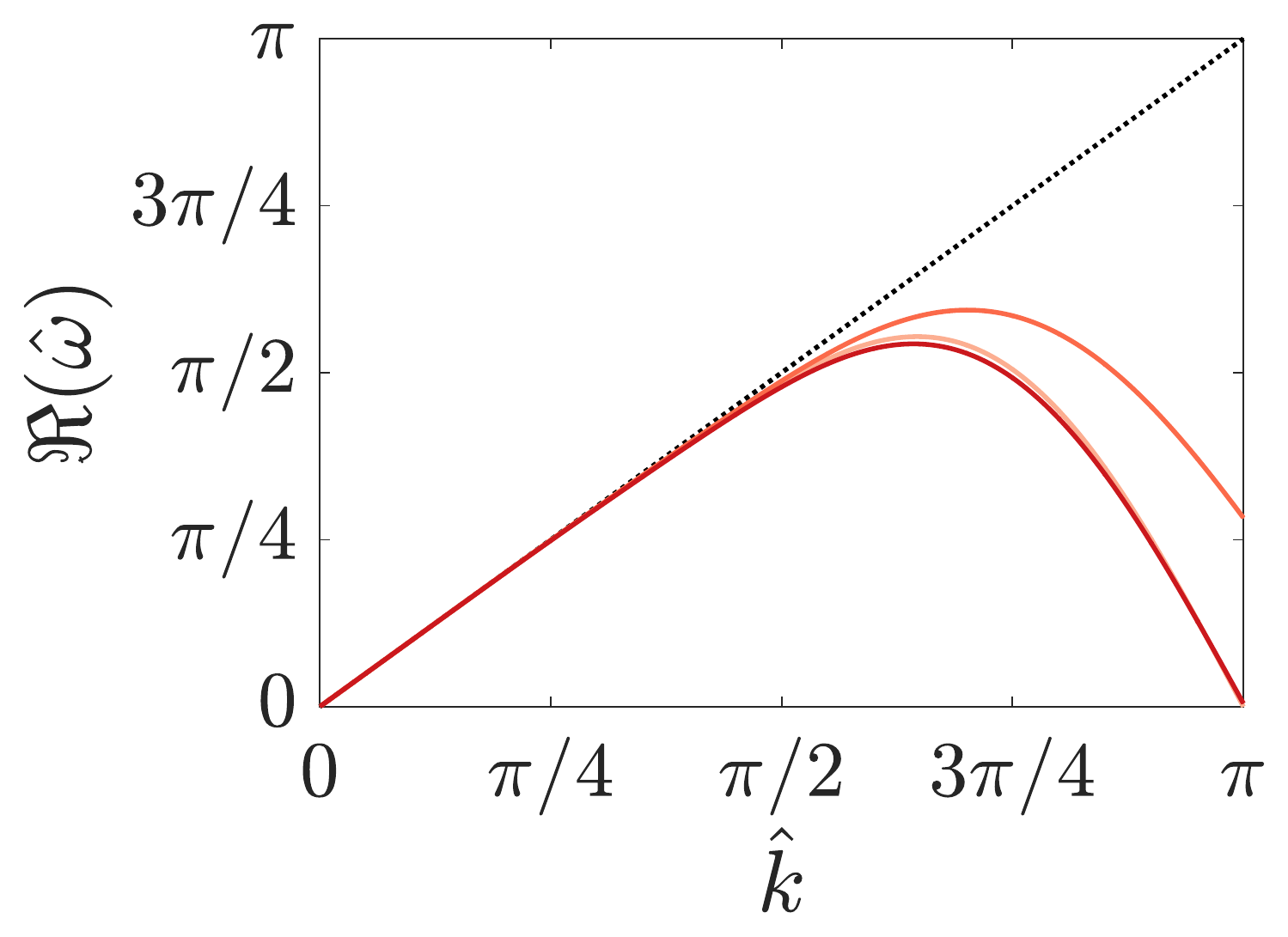}
			\caption{Dispersion: $\gamma_x=1.1$, $\gamma_y = 1$}
			\label{fig:FRp4_re_1110angle}
		\end{subfigure}
		~
		\begin{subfigure}[b]{0.45\linewidth}
			\centering
			\includegraphics[width=\linewidth]{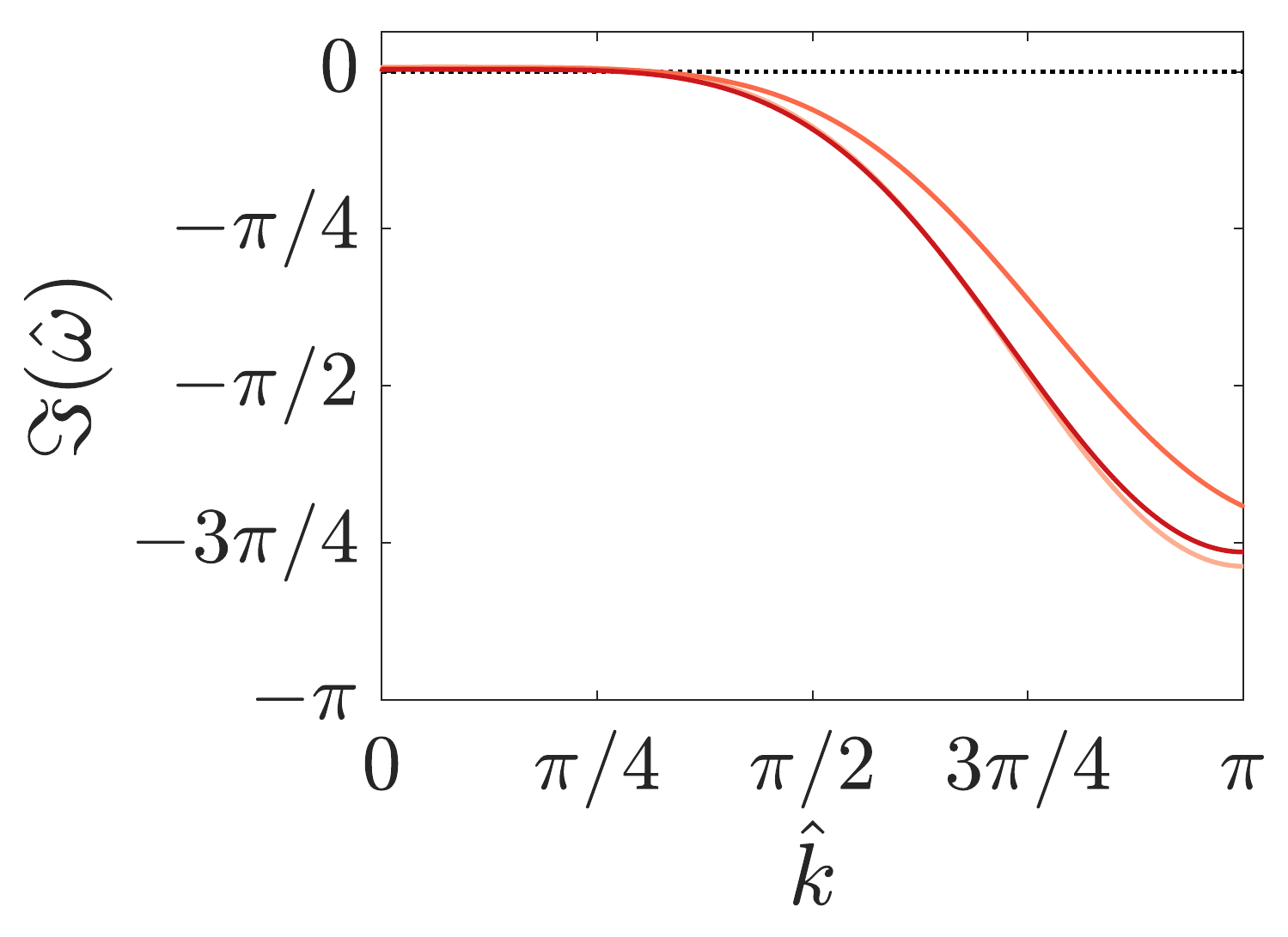}
			\caption{Dissipation: $\gamma_x=1.1$, $\gamma_y = 1$}
			\label{fig:FRp4_im_1110angle}
		\end{subfigure}
		~
		\begin{subfigure}[b]{0.48\linewidth}
			\centering
			\includegraphics[width=\linewidth]{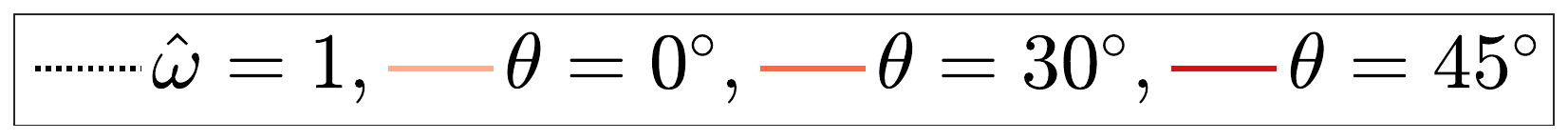}
		\end{subfigure}
		\caption{Two dimensional upwinded FR, $p=3$, with Huynh $g_2$ corrections at selected incident angles.}
		\label{fig:FRp4_angle}
	\end{figure}

	The first of these questions is explored in Figs.~\ref{fig:FR3_disp_11}~\&~\ref{fig:FR3_diss_11}. It is observed in Fig.~\ref{fig:FR3_diss_11} that expanding grids do cause positive dissipation in higher dimensions. This is more clearly displayed for some specific angles in Fig.~\ref{fig:FRp4_im_1110angle}, where the dissipation is seen to be slightly positive at low wavenumbers. In one dimension, this behaviour was previously explained that as a wave moves through elements of different size the group velocity ($\mathrm{d}\omega/\mathrm{d}k$) will change. For an expanding grid, this leads to low wavenumber energy collecting in elements. The same mechanism looks to be responsible in higher dimensions.
	
	The next question of whether a contraction orthogonal to the expansion will help to stabilise the situation is considered in Figs.~\ref{fig:FR3_disp_1109} and \ref{fig:FR3_diss_1109}. From Fig.~\ref{fig:FR3_diss_1109}, it can be seen that the answer is yes. However, there is still a region of positive dissipation for $\theta<45^\circ$ and as the incidence angle of the wave approaches $\theta=0$ the stabilisation brought about by the contraction decays. The $\theta=0$ case is identical to the case of $\gamma_y=1$.

\begin{figure}
		\centering
		\begin{subfigure}[b]{0.4\linewidth}
			\centering
			\includegraphics[width=\linewidth]{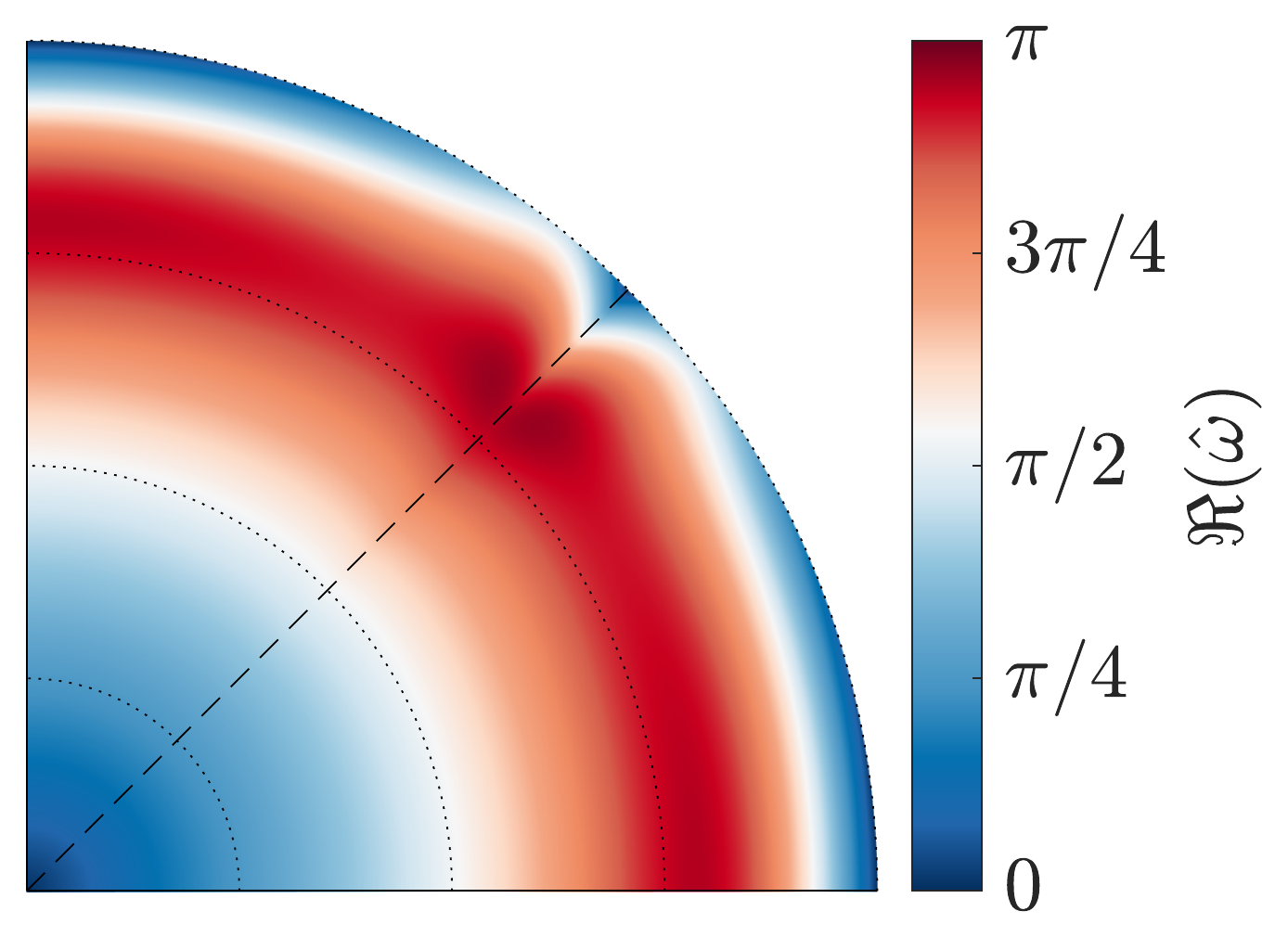}
			\caption{Dispersion, $\gamma_x=\gamma_y=1.0$}
			\label{fig:FR3_disp_1010_dg}
		\end{subfigure}
		~
		\begin{subfigure}[b]{0.4\linewidth}
			\centering
			\includegraphics[width=\linewidth]{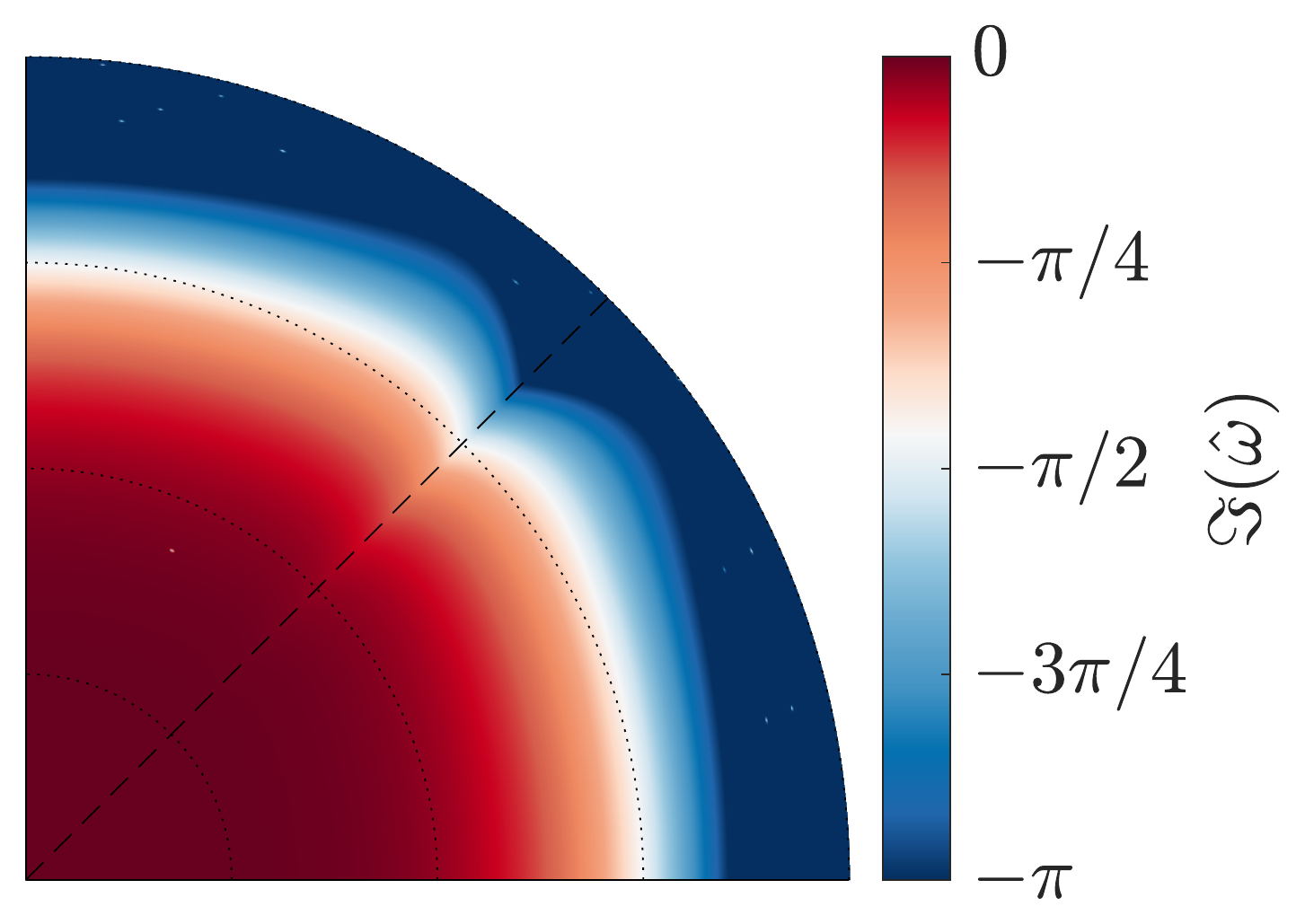}
			\caption{Dissipation, $\gamma_x=\gamma_y=1.0$}
			\label{fig:FR3_diss_1010_dg}
		\end{subfigure}
		~
		\begin{subfigure}[b]{0.4\linewidth}
			\centering
			\includegraphics[width=\linewidth]{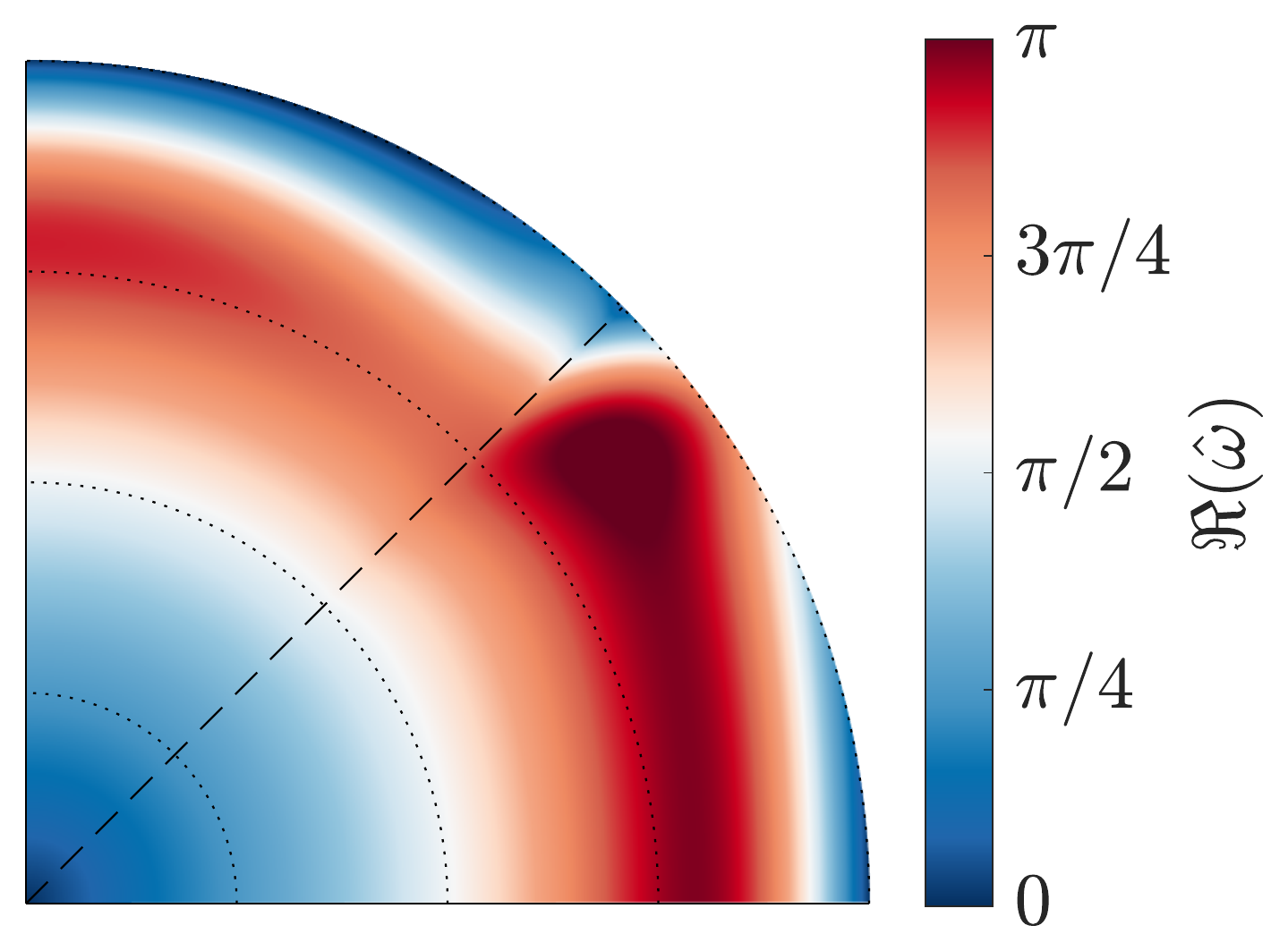}
			\caption{Dispersion, $\gamma_x=1.1$, $\gamma_y=0.9$}
			\label{fig:FR3_disp_1109_dg}
		\end{subfigure}
		~
		\begin{subfigure}[b]{0.4\linewidth}
			\centering
			\includegraphics[width=\linewidth]{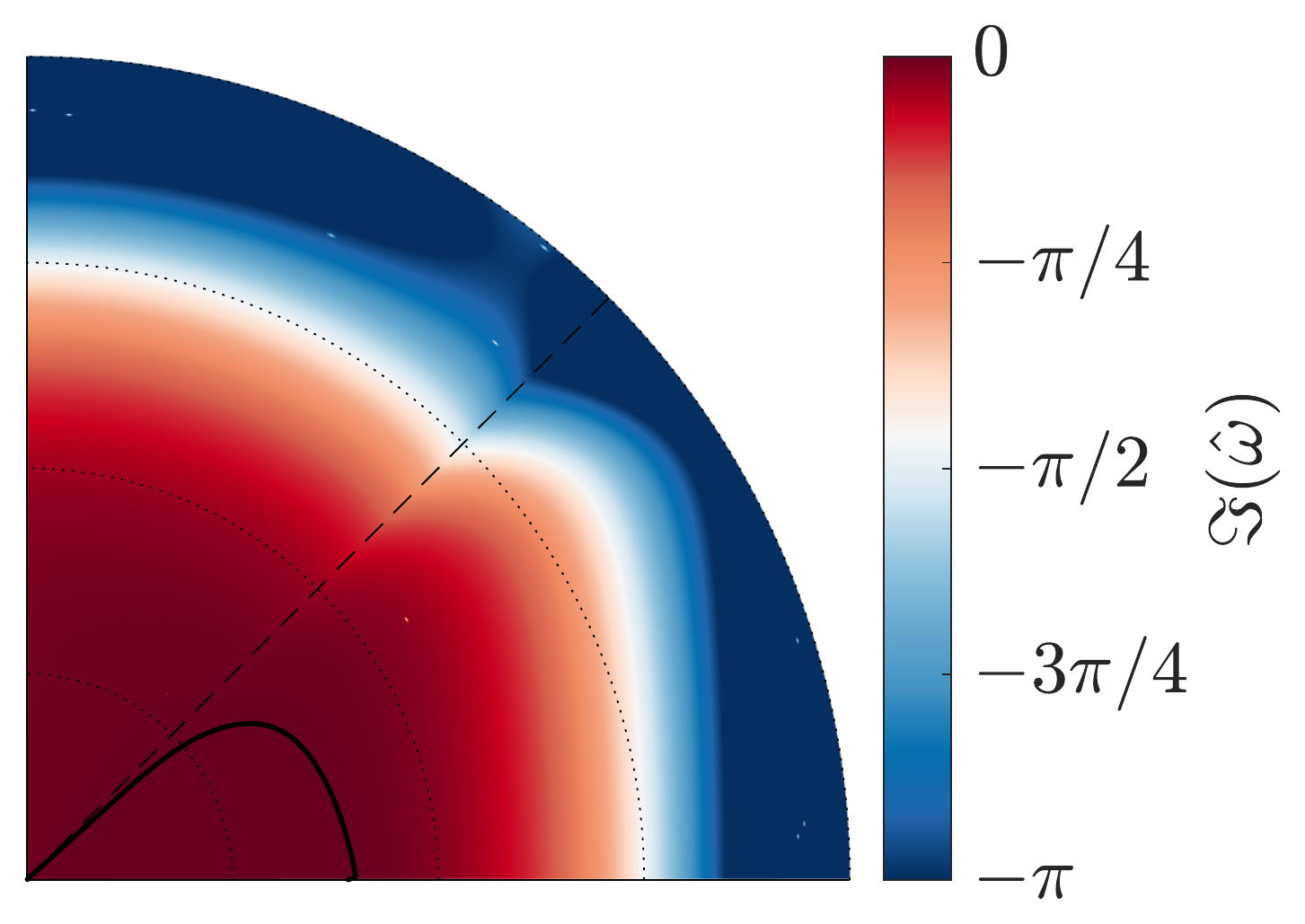}
			\caption{Dissipation, $\gamma_x=1.1$, $\gamma_y=0.9$}
			\label{fig:FR3_diss_1109_dg}
		\end{subfigure}
		\caption{Two dimensional upwinded FR, $p=3$ with NDG corrections, for different grid expansion factors. Normalised wavenumber as radial distance (markers at $\pi/4$ intervals), and element angle of incidence as azimuthal distance. The solid black line on the dissipation plots is the contour of zero dissipation.}
		\label{fig:FR_polar_disp_stretched_dg}
	\end{figure}

     For completeness, we have included plots of the dispersion and dissipation for NDG in Fig.~\ref{fig:FR_polar_disp_stretched_dg}. By comparing the un-stretched and stretched results of Fig.~\ref{fig:FR3_disp_1010_dg}-\ref{fig:FR3_diss_1010_dg} to those of Fig.~\ref{fig:FR3_disp_1109_dg}-\ref{fig:FR3_diss_1109_dg} it is clear that the effect of stretching is similar to that seen with the Huynh $g_2$ corrections in Fig.~\ref{fig:FR_polar_disp_stretched}. However, the notable between NDG and Huynh's $g_2$ corrections is that NDG exhibits more dissipation at higher wavenumbers, as well as dispersion overshoot. This is difficult to see precisely in Fig.~\ref{fig:FR_polar_disp_stretched_dg}, but is indicated by the darker band at $\hat{k}\approx 3\pi/4$. Furthermore, the location of the zero dissipation contour is different between the two correction functions, which for NDG is at a higher wavenumber. This is caused by the unstretched NDG having lower dissipation at lower wavenumbers. As is seen by comparing Fig.~\ref{fig:FR3_diss_1010_dg} and Fig.~\ref{fig:FR3_diss_11}. 
	
	Lastly, for this investigation into the dispersion and dissipation characteristics of FR, we wish to make a note on the Nyquist frequency of the elements. As was also the case for 1D stretched elements, the Nyquist frequency has a dependency on the expansion ratio. This is found from the harmonic mean of the 1D Nyquist frequencies, then normalised by the element size at that angle. Hence, the normalised wavenumber is then:
	\begin{equation}
		\hat{k} = k/k_\mathrm{nq} = k\max{\big\{\cos{(\theta)},\sin{(\theta)}\big\}}\bigg(\frac{1}{p+1}\bigg)\sqrt{\bigg(\frac{\cos{(\theta)}}{\gamma_x}\bigg)^2 + \bigg(\frac{\sin{(\theta)}}{\gamma_y}\bigg)^2}
	\end{equation}
	
	\subsection{Effect of Grid Aspect Ratio on CFL Limit}
	Beyond the resolution of the scheme, is the question of setting up a case and running it on some machine. For this, knowledge of the temporal stability is key and we will begin by looking at the effect of the relative size of an element in $x$ and $y$ on the CFL limits of FR. In this case, the grid is not expanding or contracting, merely the ratio of $\Delta_x$ to $\Delta_y$ is varied.    

\begin{figure}
		\centering
		\begin{subfigure}[b]{0.45\linewidth}
			\centering
			\includegraphics[width=\linewidth]{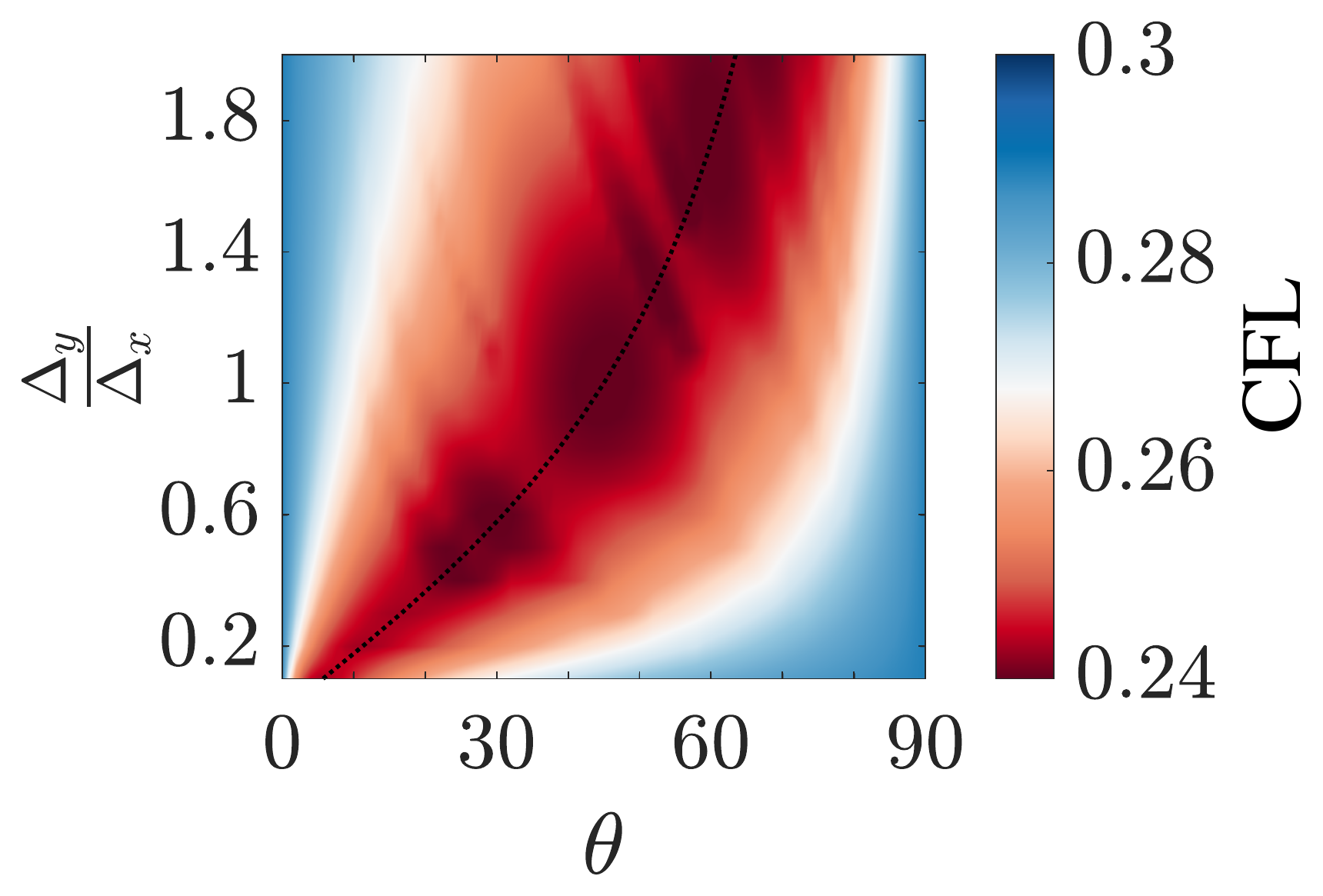}
			\caption{$p=3$}
			\label{fig:FR3_RK44_dy_theta}
		\end{subfigure}
		~
		\begin{subfigure}[b]{0.46\linewidth}
			\centering
			\includegraphics[width=\linewidth]{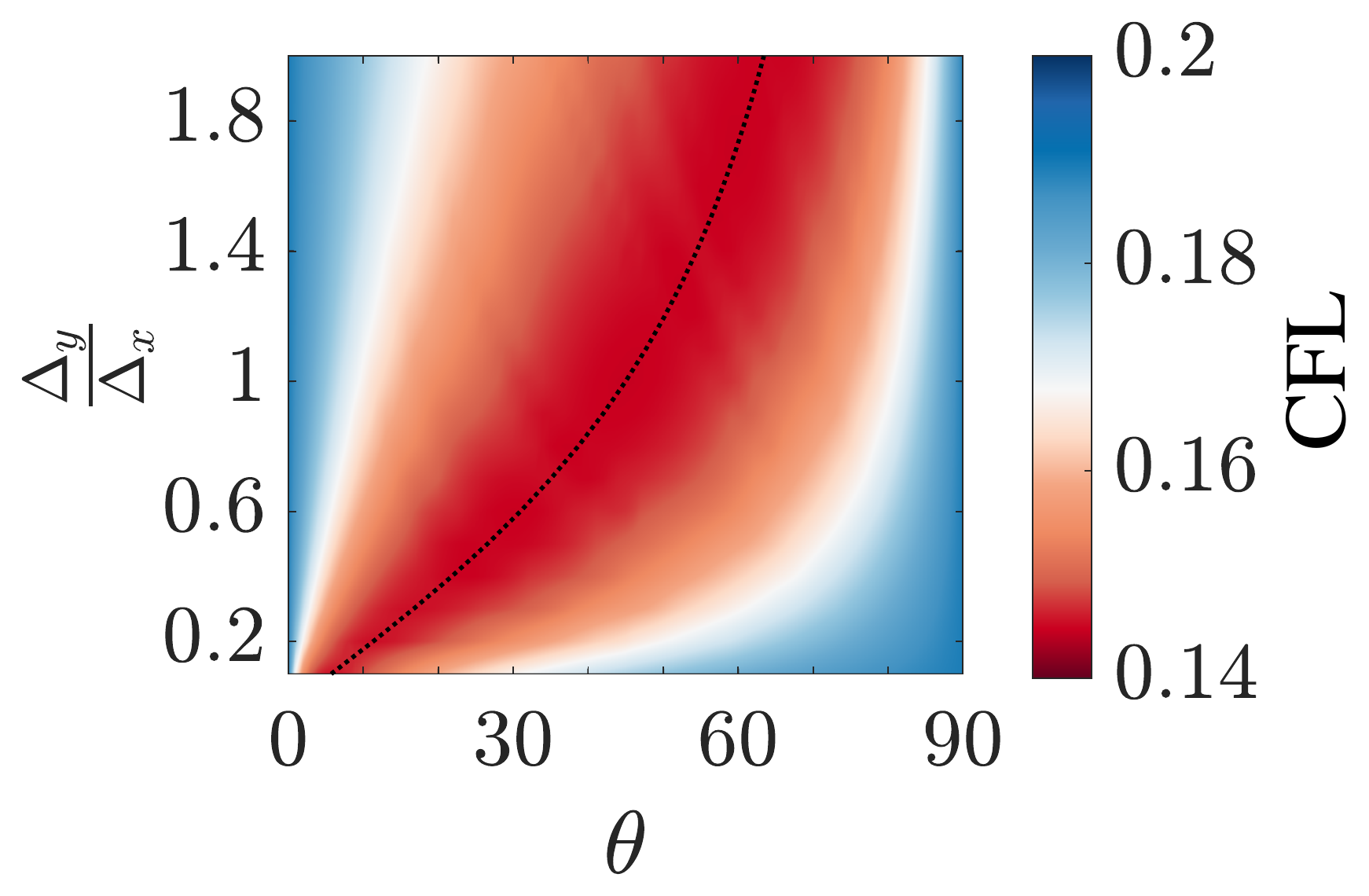}
			\caption{$p=4$}
			\label{fig:FR4_RK44_dy_theta}
		\end{subfigure}
		\caption{Effect of varying: the relative size an element in $x$ and $y$; the angle of incidence; and CFL in 2D for upwinded FR. The method of temporal integration used here was RK44. (\dotted) $\Delta_y/\Delta_x = \tan{\theta}$}
		\label{fig:FR_dy_theta}
	\end{figure}

    As can be seen from Fig.~\ref{fig:FR_dy_theta} there is a clear impact on the CFL limit of FR when elements are rectangular, with no change in the CFL limit when the waves are aligned with the grid. What is evident is that the angle of incidence where the CFL limit is smallest, for a given size ratio, is when a wave is incident at an angle of $\tan^{-1}{(\Delta_y/\Delta_x)}$. This angle corresponds to the maximum length across the element and hence, when the wave is decomposed into $x$ and $y$ components, the wavenumbers are lowest --- \emph{i.e.}, $\inf_{\theta\in\mathbb{R}}({\sup{\{\Delta_y\cos{(\theta)},\Delta_x\sin{(\theta)}\}})}$  when $\theta=\tan^{-1}{(\Delta_y/\Delta_x)}$. As the wavenumbers will be at their lowest necessary to form the wave, from the 1D dissipation of FR, the dissipation will also be at its lowest. Therefore, there is less dissipation in the spatial scheme available to counteract the negative dissipation of the temporal scheme, hence reducing temporal stability at $\theta=\tan^{-1}{(\Delta_y/\Delta_x)}$.  
	
	\subsection{Effect of Grid Expansion on CFL Limit}\label{sec:geo_cfl}
\begin{figure}
		\centering
		\begin{subfigure}[b]{0.4\linewidth}
			\centering
			\includegraphics[width=\linewidth]{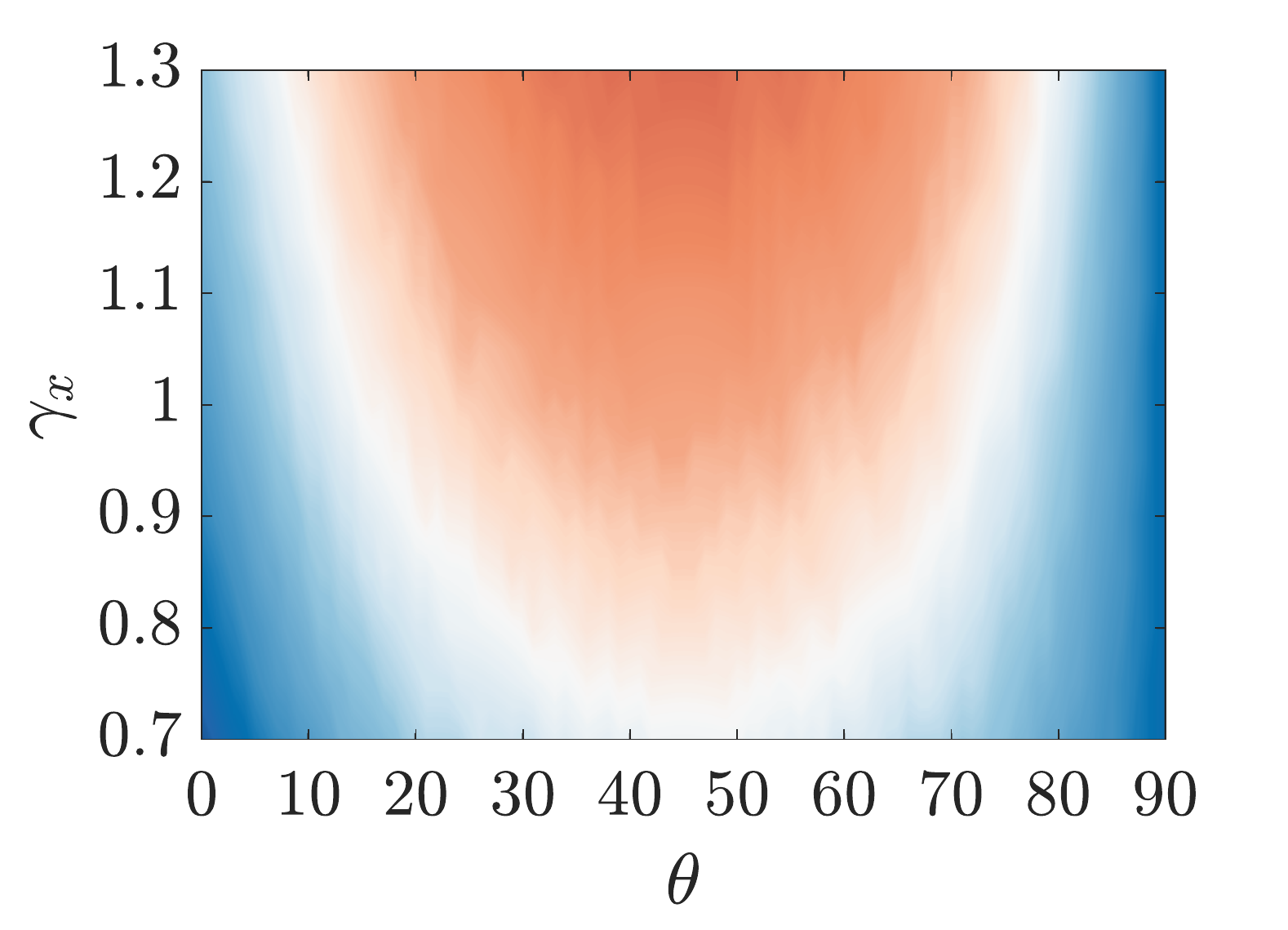}
			\caption{$\gamma_y = 0.8$}
			\label{fig:FR4_2d_hu_08_x}
		\end{subfigure}
		~
		\begin{subfigure}[b]{0.4\linewidth}
			\centering
			\includegraphics[width=\linewidth]{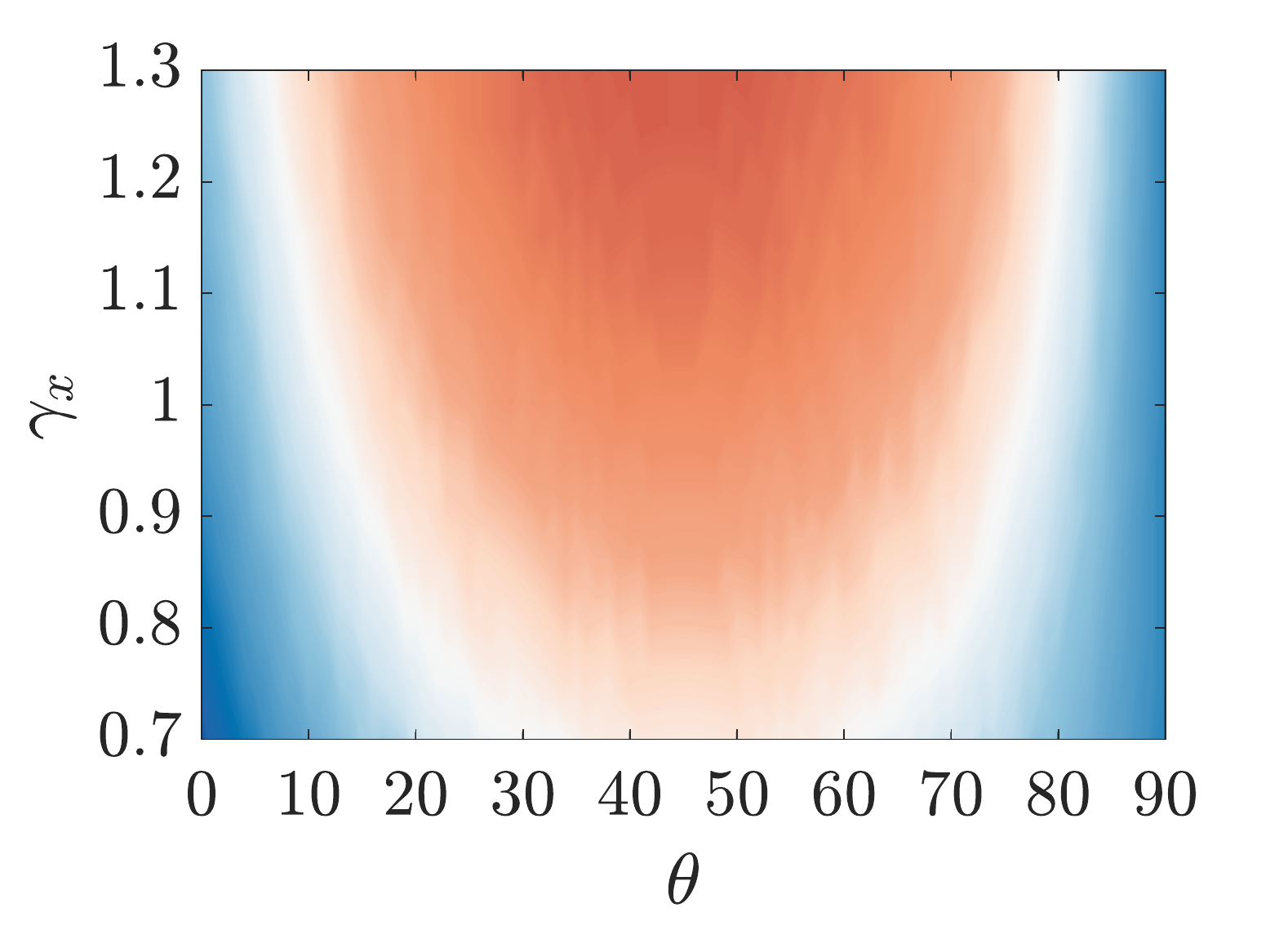}
			\caption{$\gamma_y = 0.9$}
			\label{fig:FR4_2d_hu_09_x}
		\end{subfigure}
		~
		\begin{subfigure}[b]{0.4\linewidth}
			\centering
			\includegraphics[width=\linewidth]{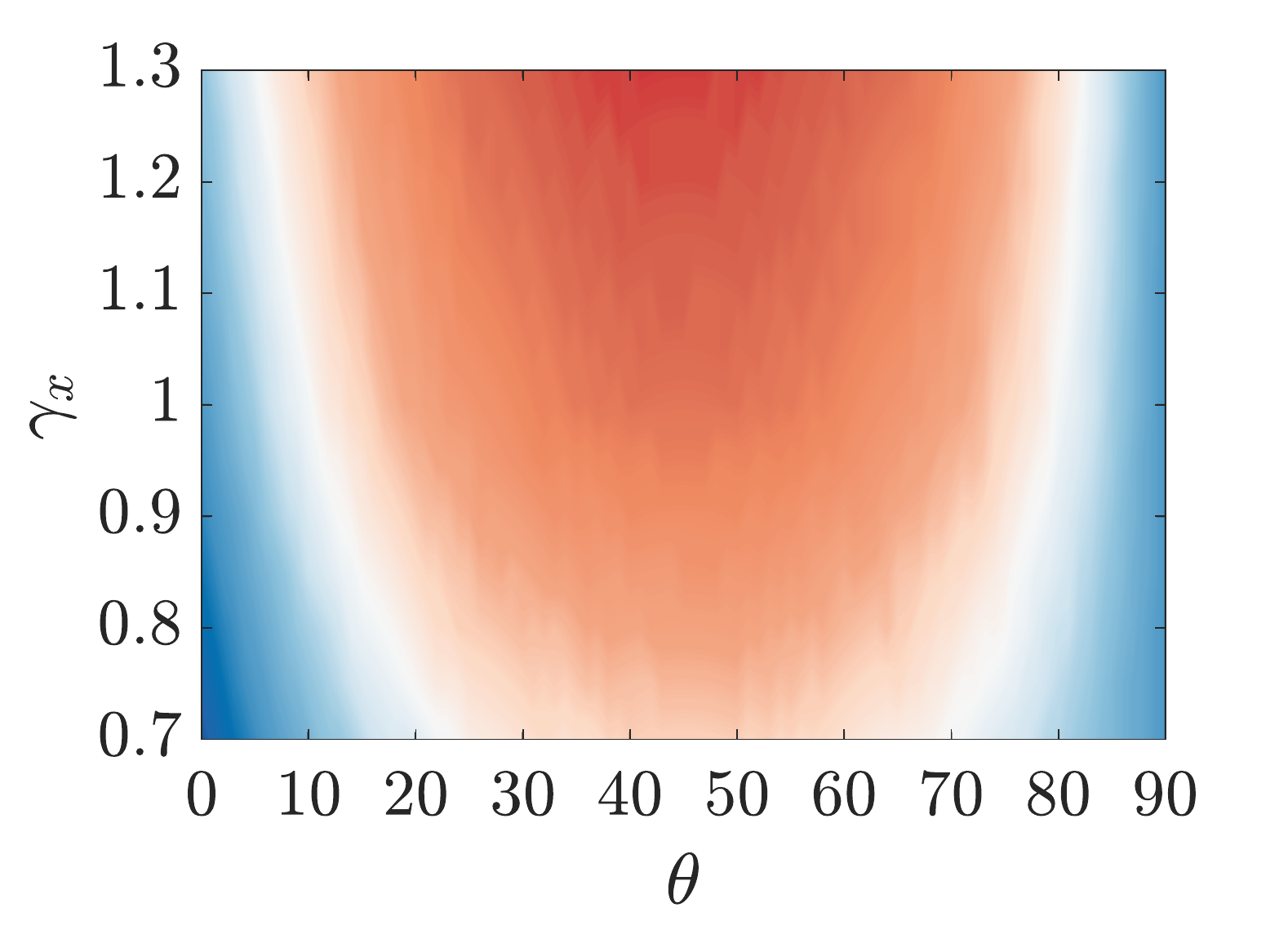}
			\caption{$\gamma_y = 1.0$}
			\label{fig:FR4_2d_hu_10_x}
		\end{subfigure}
		~
		\begin{subfigure}[b]{0.4\linewidth}
			\centering
			\includegraphics[width=\linewidth]{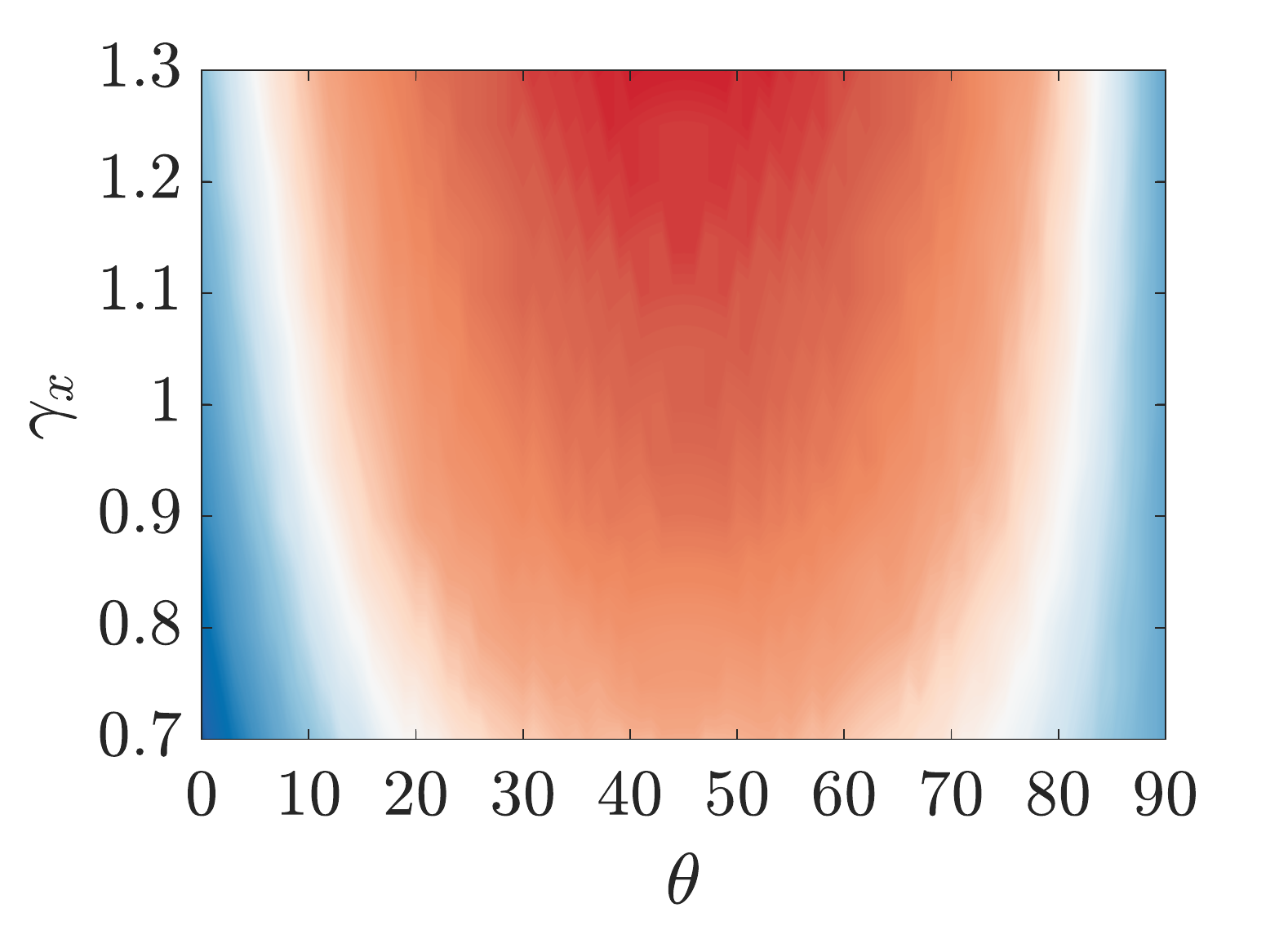}
			\caption{$\gamma_y = 1.1$}
			\label{fig:FR4_2d_hu_11_x}
		\end{subfigure}
		~
		\begin{subfigure}[b]{0.3\linewidth}
			\centering
			\includegraphics[width=\linewidth]{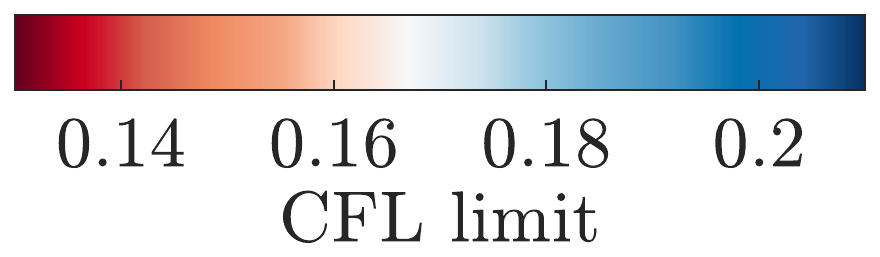}
		\end{subfigure}
		\caption{CFL limit for 2D linear advection with FR ($p=4$) using Huynh correction function, showing variation with $\theta$ and $\gamma_x$ for some set values of $\gamma_y$. Time integration is RK44.}
		\label{fig:FR4_2D_HU_CFL}
	\end{figure}

	Now we introduce to the grid an expansion or contraction in $x$ and $y$, with varying incident angles. The results of these distortions are shown in Fig.~\ref{fig:FR4_2D_HU_CFL}.  The minimum CFL limit is seen to be at $\theta=45^{\circ}$, with temporal performance peaking as expected at $\theta=0^{\circ},90^{\circ}$. This result agrees with that of Fig.~\ref{fig:FR_dy_theta}. However, this also shows that the angle of minimum CFL is only dependent on the local element shape, and in the case investigated here the central element is always square. Furthermore, when the CFL limit in the quasi-1D case ($\theta=0^{\circ},90^{\circ}$) is compared to the results of Trojak~\etal~\cite{Trojak2017a}, the CFL limit is found to be lower than the 1D case. This may be due in part to the increased modes of the system and their coupling leading to a less stable system. This decrease is corroborated by numerical tests, and the analytical reasoning will follow shortly.

    A second point to note is that for non-grid-aligned waves the expansion or contraction in both components affects stability. With a contraction orthogonal to expansion again helping to stabilise the scheme. This point is subtle, because if the decomposition of the wave into $x$ and $y$ were linearly independent then it would be expected that the lowest CFL limit would dominate. However, this result demonstrates that there is coupling between the $x$ and $y$ components --- which could be used advantageously, as was discussed earlier.
	
	Throughout the analytical tests in which waves were injected at incidence on a square central element, the Nyquist wavenumber was found to be $k_{nq,\theta} = k_{nq,0}/\cos(\theta)$ for $0^{\circ} \leqslant \theta \leqslant 45^{\circ}$. This result can be understood in one of two ways. Firstly that a wave at an angle can draw on more points in the normal direction to form a fit of higher wavenumbers. Or that the wave can be thought of as being a coupled decomposition of the wave into the $x$ and $y$ directions, and although it has been said that these are not independent, this does mean that higher wavenumbers can be supported.
	
	To understand further why expanding meshes are less stable and contracting meshes more stable it can be illuminating to consider the 1D linear advection equation, which for upwinded FR can be written as:
	\begin{equation}
		\frac{\partial \mathbf{u}_j}{\partial t} = -J^{-1}_j\mathbf{C}_0\mathbf{u}_j - J^{-1}_{j-1}\mathbf{C}_{-1}\mathbf{u}_{j-1}
	\end{equation}
	where $\mathbf{C}_{-1}$ and $\mathbf{C}_{0}$ are defined similarly to $\mathbf{C}_L$ and $\mathbf{C}_{0,\xi}$ in Eq.(\ref{eq:Cxmat}). Then using Euler's method to temporally integrate this equation, we find:
	\begin{equation}
		\mathbf{u}^{n+1}_j = \mathbf{u}^n_j - \tau \big(J^{-1}_j\mathbf{C}_0\mathbf{u}^n_j + J^{-1}_{j-1}\mathbf{C}_{-1}\mathbf{u}^n_{j-1}\big) = \big (\mathbf{I} - \tau J^{-1}_j\mathbf{C}_0\big )\mathbf{u}^n_j - \tau J^{-1}_{j-1}\mathbf{C}_{-1}\mathbf{u}^n_{j-1}
	\end{equation}	 
	where the superscript $n$ denotes the time step. If $u^{n-m}_{j-m-1} (\forall \: m \in \mathbb{N}$) is then recursively substituted, the final form is then:
	\begin{equation}\label{eq:upwind_infinite}
		\mathbf{u}^{n+1}_{j} = \frac{1}{2}\sum^{\infty}_{m=0}(\mathbf{I}-\tau(2\gamma)^{m}\mathbf{C_0})(\underbrace{-2\gamma\tau\mathbf{C}_{-1}}_{\mathbf{T}})^m\mathbf{u}^{n-m}_{j-m}
	\end{equation}
	 where we assume the solution is on a geometrically expanding grid in order to substitute for the Jacobian --- hence being only valid for linearly transformed elements. If we consider that $\mathbf{C}_0$ and $\mathbf{C}_{-1}$ are linear operators, then rather than prescribing a solution, the dynamics of linear operators can be used~\cite{Bayart2009}. So, if the mesh extends infinitely downwind, then it is sufficient to say that Eq.(\ref{eq:upwind_infinite}) is stable when it is a hypercylic orbit. Hence, the stability criterion is that $-2\gamma\tau\mathbf{C}_{-1} = \mathbf{T}$ is a matrixable linear hypercylic operator. The definition of which is that $\sup{\|\mathbf{T}^n\|} \leqslant 1, \forall\: n\in\mathbb{N}$, which in turn implies that $\rho(\mathbf{T}) = 1$. What this aims to show is that the stability criterion is dependent on the product of $\tau$ and $\gamma$, as well as on $\mathbf{C}_0$. Therefore, as $\gamma$ increases the maximum stable $\tau$ decreases, for constant $\mathbf{C}_{-1}$. This also explains that although the scheme may be formally unstable, correct setting of $\tau$ for a given $\gamma$ can lead to a $\mathbf{T}$ that is still hypercyclic and give a bounded solution. However because of the $(-1)^m$ this does mean that if $\mathbf{U}_{n-1} \subset \mathbb{R}$ is the set of solutions at some time step $n-1$, then the solution $\mathbf{u}^n_j \nsubseteq \mathbf{U}_{n-1}$ \emph{i.e}, for a sinusoidal solution the computed value may exceed the prescribed magnitude.
	
	\subsection{Effect of Grid Expansion of Fully Discrete Dispersion and Dissipation}	
		Following on from the exploration of grid expansion on temporal stability limits, we will present the fully discretised Fourier analysis in two dimensions. In this investigation, we look at the dispersion and dissipation of Huynh $g_2$ correction functions for several expansion ratios. Throughout this investigation as the angle is swept through from $0-90^\circ$ the time step will be held constant as this is reflective of practical applications.
		
	\begin{figure}
		\centering
		\begin{subfigure}[b]{0.45\linewidth}
			\centering
			\includegraphics[width=\linewidth]{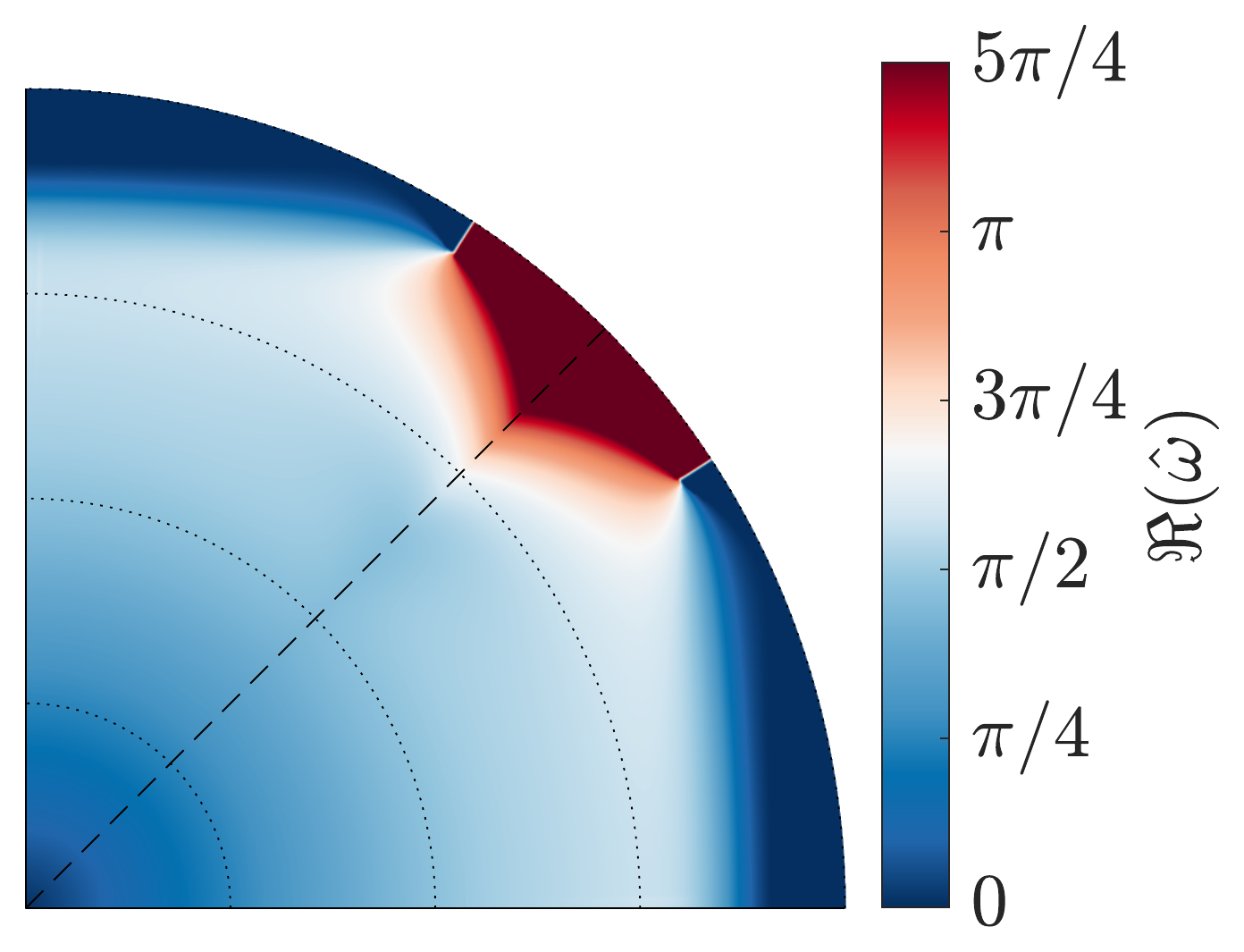}
			\caption{Dispersion: $\gamma_x=\gamma_y = 1$}
			\label{fig:FRp3_re_1010_t18}
		\end{subfigure}
		~
		\begin{subfigure}[b]{0.45\linewidth}
			\centering
			\includegraphics[width=\linewidth]{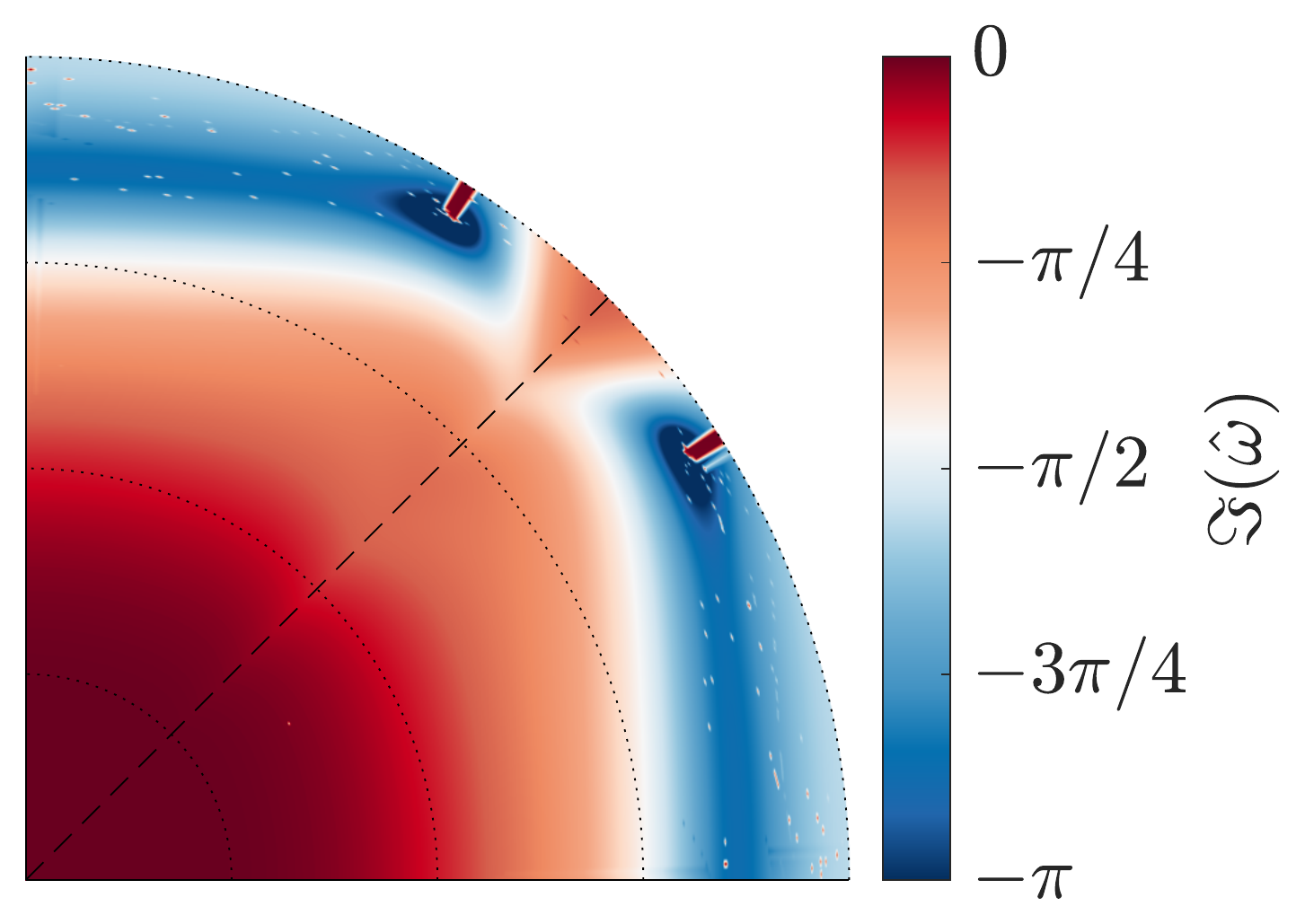}
			\caption{Dissipation: $\gamma_x=\gamma_y = 1$}
			\label{fig:FRp3_im_1010_t18}
		\end{subfigure}
		~
		\begin{subfigure}[b]{0.45\linewidth}
			\centering
			\includegraphics[width=\linewidth]{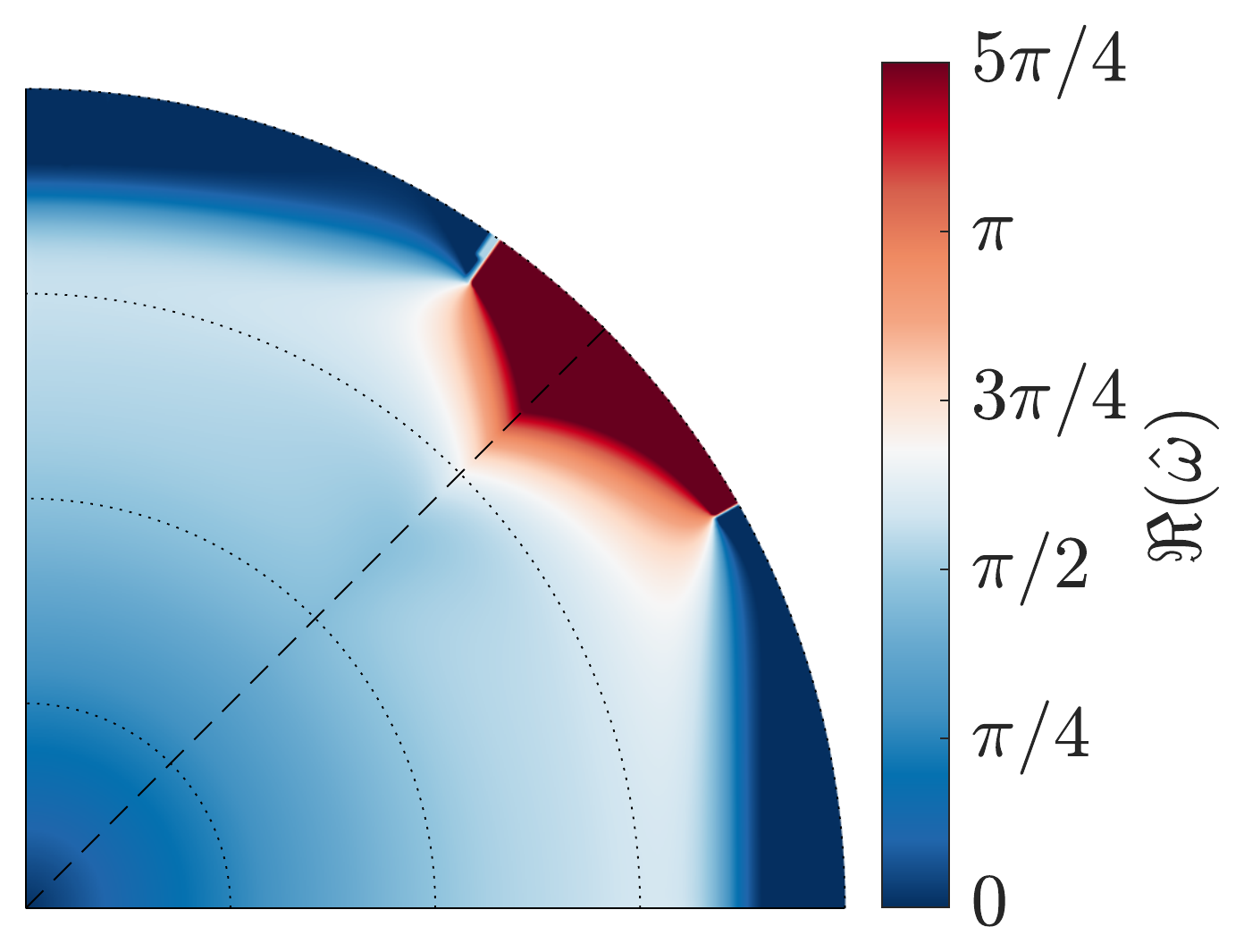}
			\caption{Dispersion: $\gamma_x=1.1$, $\gamma_y = 1$}
			\label{fig:FRp3_re_1110_t18}
		\end{subfigure}
		~
		\begin{subfigure}[b]{0.45\linewidth}
			\centering
			\includegraphics[width=\linewidth]{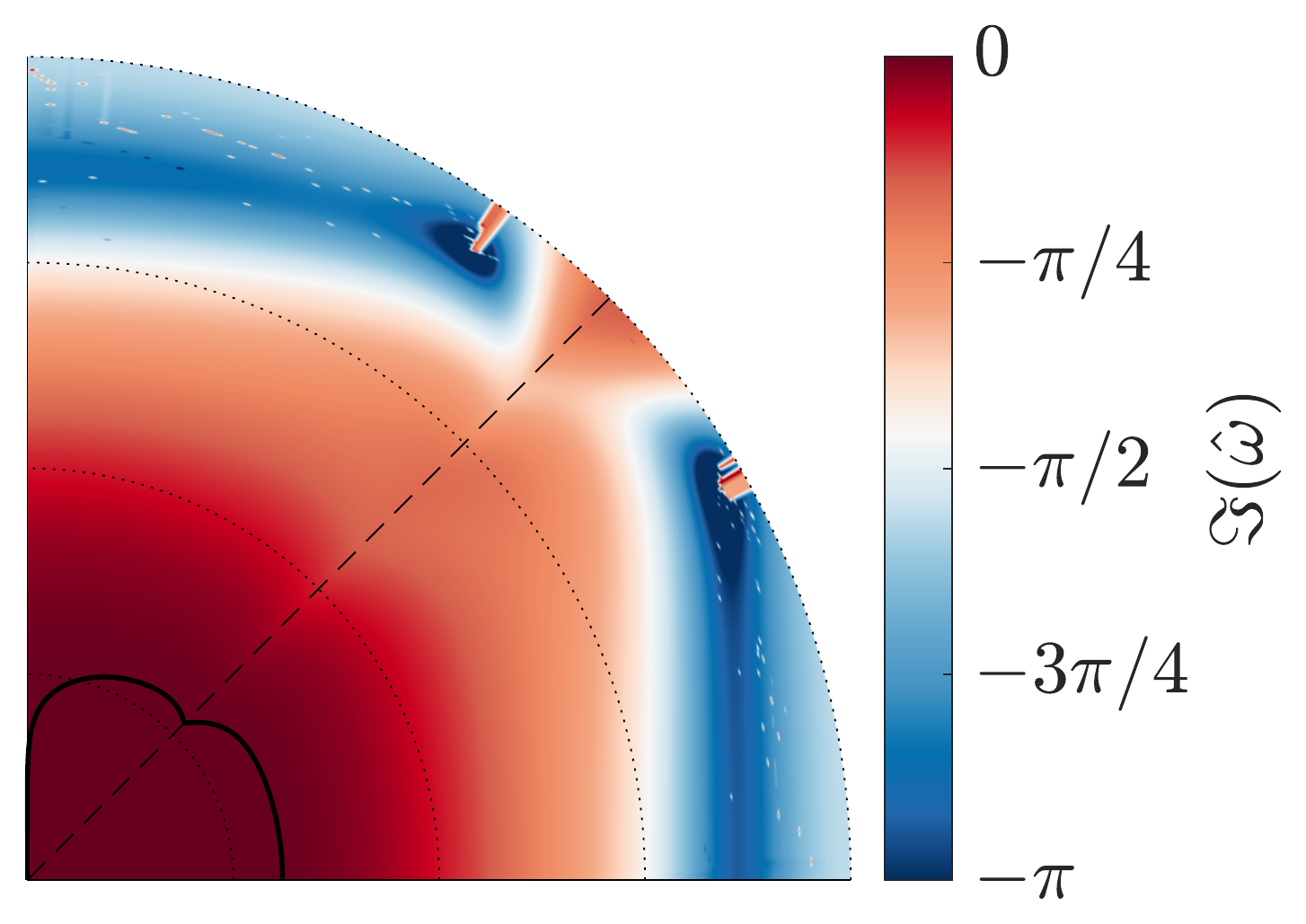}
			\caption{Dissipation: $\gamma_x=1.1$, $\gamma_y = 1$}
			\label{fig:FRp3_im_1110_t18}
		\end{subfigure}
		~
		\begin{subfigure}[b]{0.45\linewidth}
			\centering
			\includegraphics[width=\linewidth]{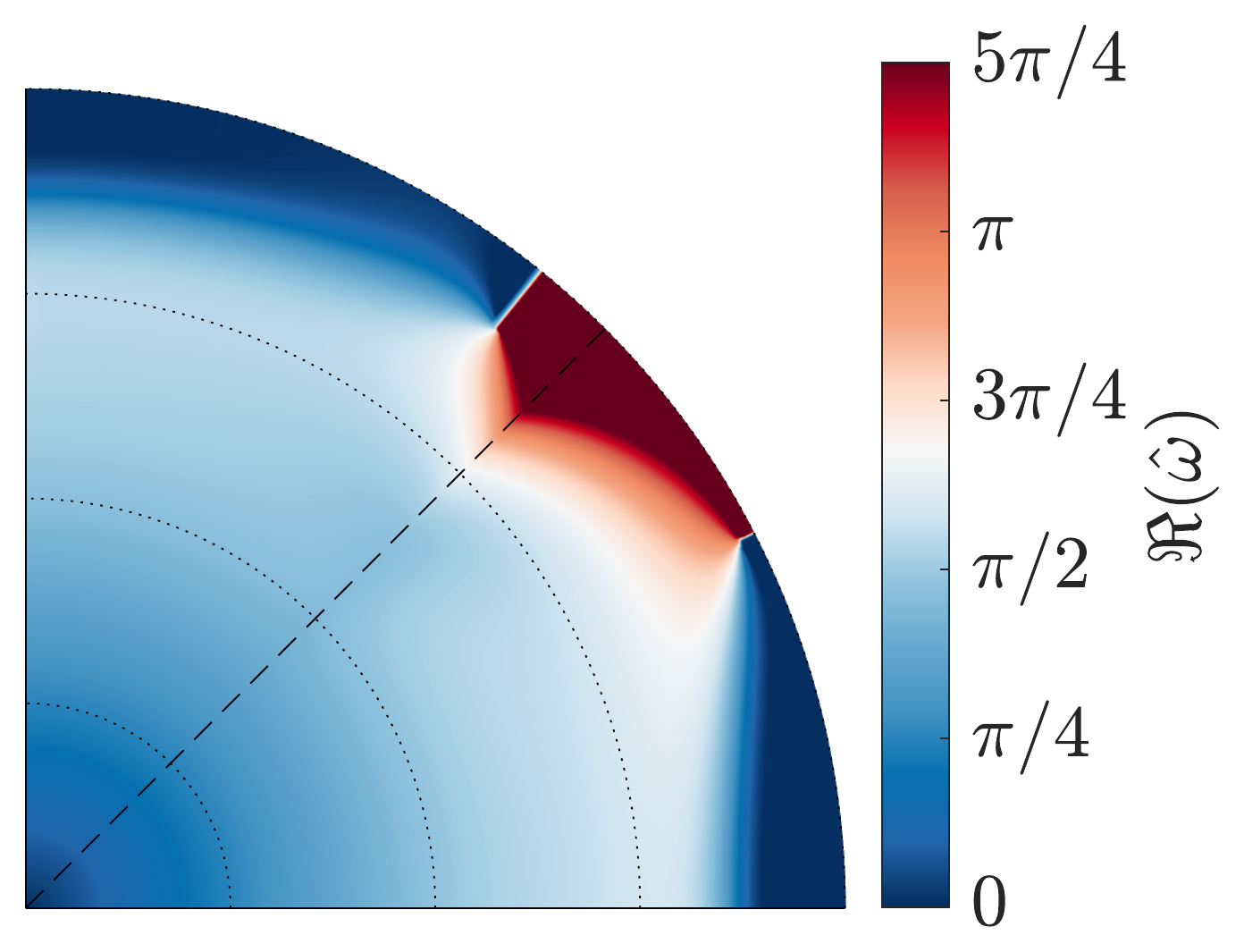}
			\caption{Dispersion: $\gamma_x=1.1$, $\gamma_y = 0.9$}
			\label{fig:FRp3_re_1109_t18}
		\end{subfigure}
		~
		\begin{subfigure}[b]{0.45\linewidth}
			\centering
			\includegraphics[width=\linewidth]{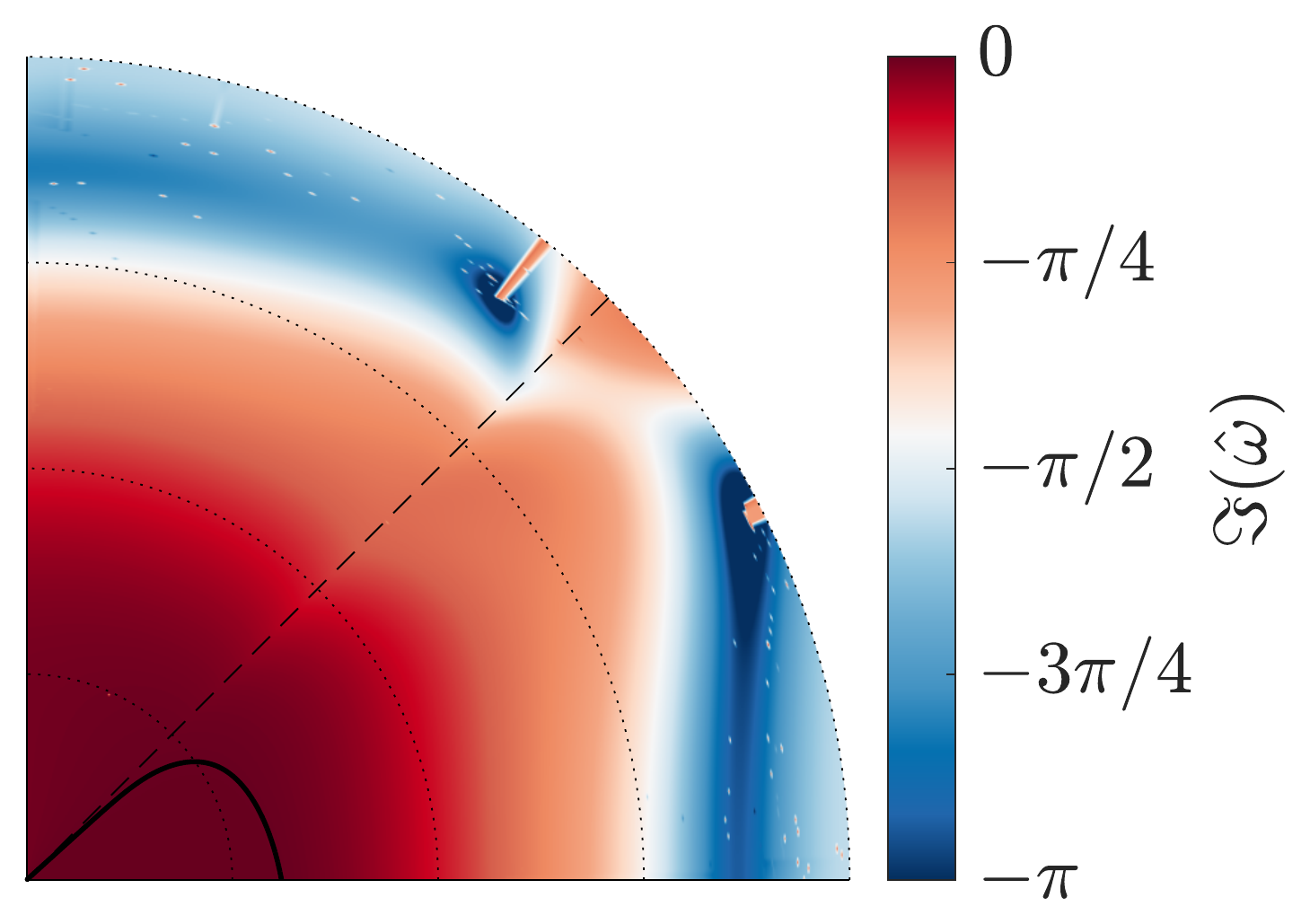}
			\caption{Dissipation: $\gamma_x=1.1$, $\gamma_y = 0.9$}
			\label{fig:FRp3_im_1109_t18}
		\end{subfigure}
		\caption{Dispersion and dissipation of upwinded FR, $p=3$, with Huynh $g_2$ corrections and explicit RK44 temporal integration, $\tau=0.18$. The radial distance is the normalised wavenumber (including the effect of angle), and the azimuthal distance is the angle of incidence to the element. The solid black line on the dissipation plots is the contour of zero dissipation.}
		\label{fig:FRp3_hu_fd}
	\end{figure}	
			
	Figure~\ref{fig:FRp3_hu_fd} displays the dispersion and dissipation relation for Huynh $g_2$ correction functions when fully discretised with RK44 explicit temporal integration. It should be noted that there are some anomalous artefacts in the data due to the complexity of sorting and selecting the eigenvalues. Primarily this shows that when fully discretised, the instability on expanding grids is still present. This is then coupled to the same behaviour that was observed by Vermerie~\etal~\cite{Vermeire2017a} and Trojak~\etal~\cite{Trojak2018c}. Namely that the gradient of the dispersion can be much larger and hence the magnitude group velocity can be very large. Also observed, which is seen here, is that as the explicit time step is increased the dissipation is reduced, with the largest difference seen at higher wavenumbers. As a result the wavenumber at which $\Im{(\hat{\omega})}=0$ (the solid black lines) in Fig~\ref{fig:FRp3_im_1110_t18}~\&~\ref{fig:FRp3_im_1109_t18} is not appreciably changed compared to the semi-discrete results. However, in Fig~\ref{fig:FRp3_hu_fd} there are regions of significantly increased and decreased dissipation at high wavenumbers. The general impact that these results show is that, when fully discretised, FR becomes more heterogeneous.
	
	\subsection{Effect of Correction Functions with Grid Expansion}
\begin{figure}
		\centering
		\begin{subfigure}[b]{0.4\linewidth}
			\centering
			\includegraphics[width=\linewidth]{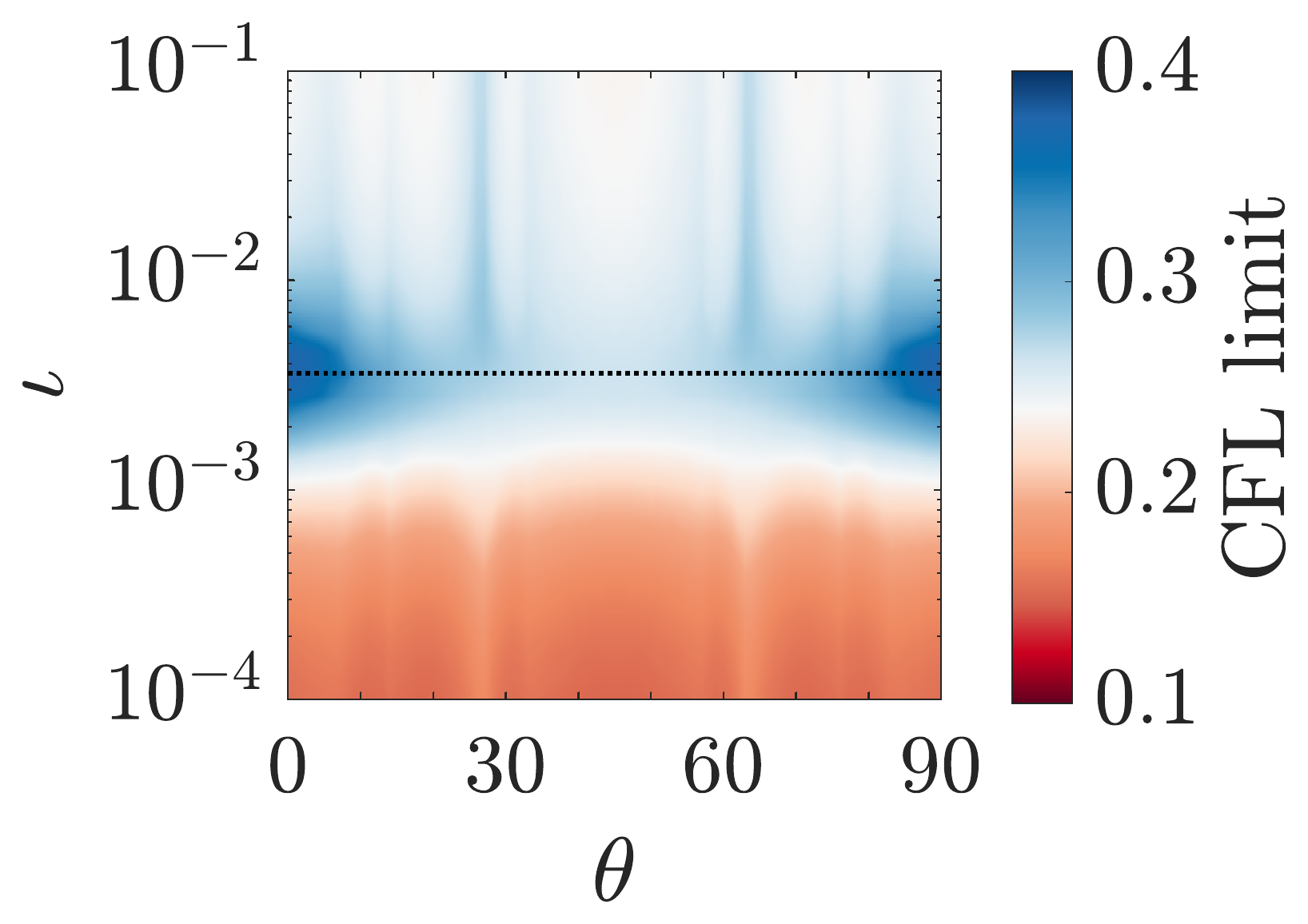}
			\caption{$p=3$, $\gamma_x = 1.0$ \& $\gamma_y=1.0$}
			\label{fig:FR3_RK44_s_y10_x10}
		\end{subfigure}
		~
		\begin{subfigure}[b]{0.4\linewidth}
			\centering
			\includegraphics[width=\linewidth]{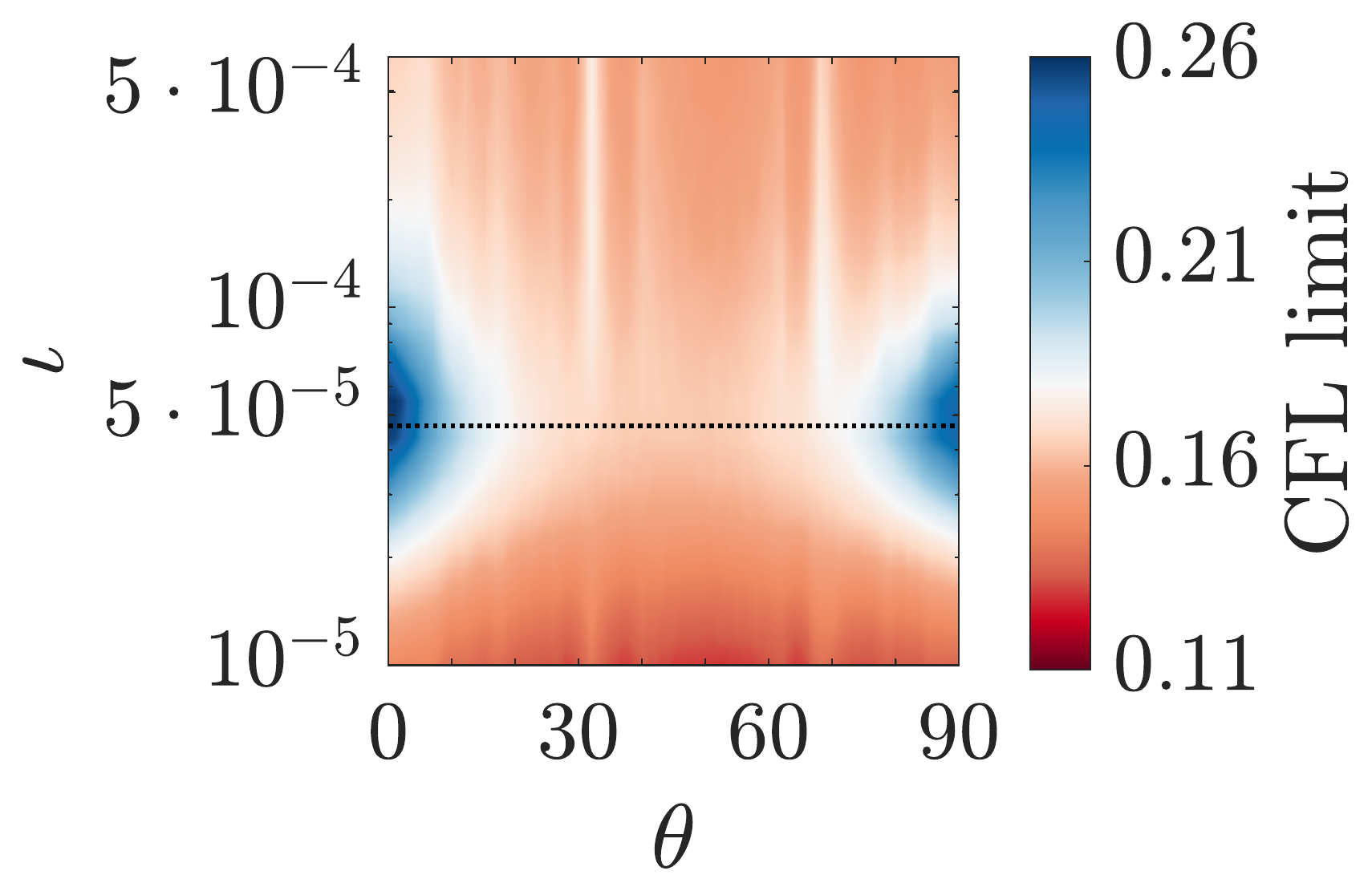}
			\caption{$p=4$, $\gamma_x = 0.8$ \& $\gamma_y=1.0$}
			\label{fig:FR4_RK44_s_y10_x08}
		\end{subfigure}
		~
		\begin{subfigure}[b]{0.4\linewidth}
			\centering
			\includegraphics[width=\linewidth]{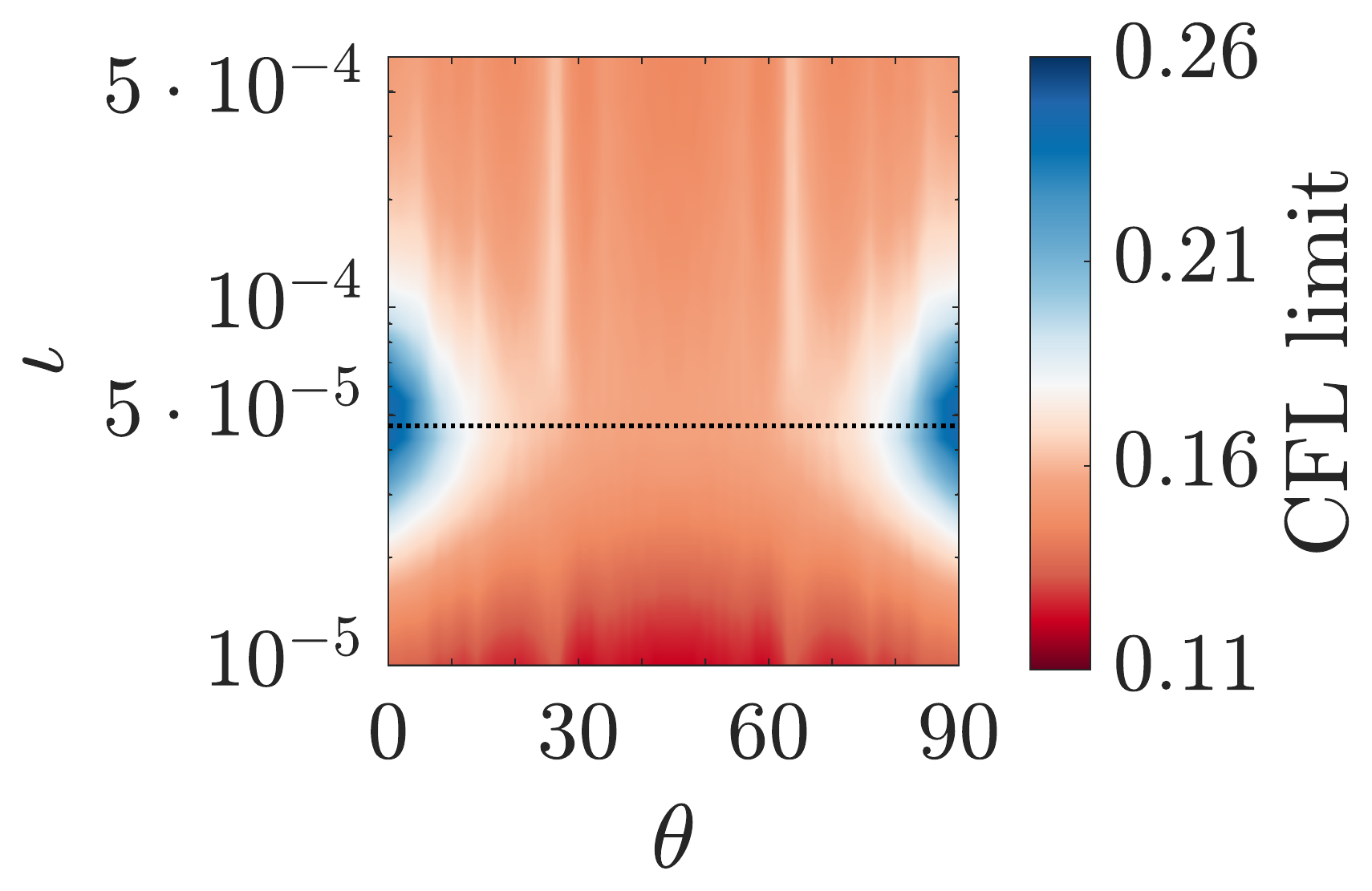}
			\caption{$p=4$, $\gamma_x = 1.0$ \& $\gamma_y=1.0$}
			\label{fig:FR4_RK44_s_y10_x10}
		\end{subfigure}
		~
		\begin{subfigure}[b]{0.4\linewidth}
			\centering
			\includegraphics[width=\linewidth]{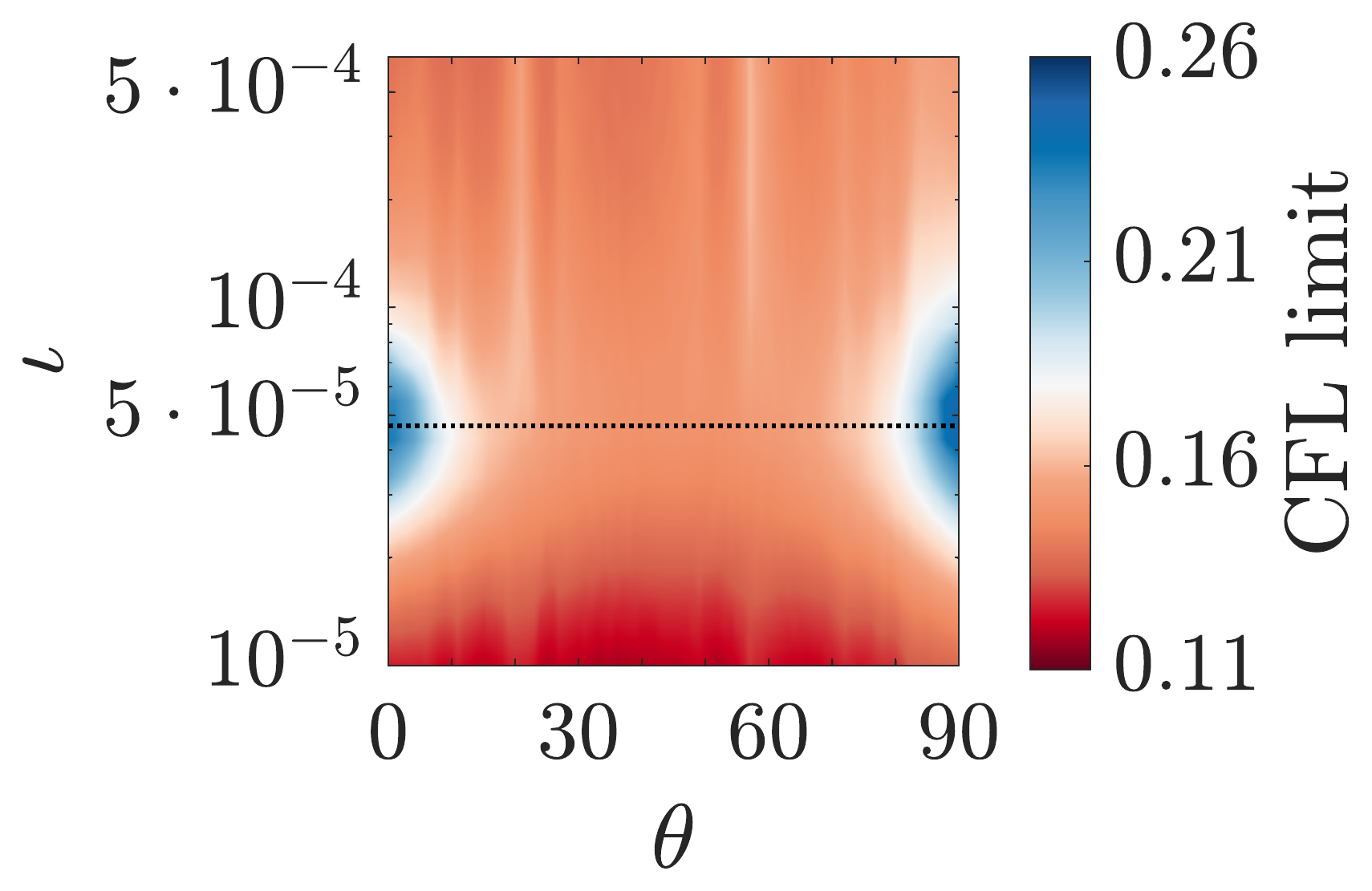}
			\caption{$p=4$, $\gamma_x = 1.3$ \& $\gamma_y=1.0$}
			\label{fig:FR4_RK44_s_y10_x13}
		\end{subfigure}
		\caption{CFL limit for 2D linear advection, at several orders. Varying correction function parameter, $\iota$, and angle $\theta$. Time integration is RK44. The value of $\iota_+$ is shown as a dashed black line.}
		\label{fig:FR_S_theta_RK44}
	\end{figure}

	As was mentioned in section~\ref{sec:intro}, a series of correction functions with peak temporal stability and spatial accuracy was proposed by Vincent~\etal~\cite{Vincent2011}, defined by a correction parameter, $\iota_+$. These correction functions exhibited the superconvergence expected of Nodal DG (Cockburn~\etal~\cite{Cockburn1999}), however $\iota_+ \neq \iota_{\mathrm{DG}}$ as they account for the variation in the dispersion caused by the discrete temporal integration. We wish to investigate if the advantage of this family of correction functions is maintained in 2D or if an analogue can be found. Therefore, the correction function parameter is varied for different angles and grid expansion rate, the results of which are presented in Fig.~\ref{fig:FR_S_theta_RK44}.
	Using Fig.~\ref{fig:FR4_RK44_s_y10_x10} as an example, the peak CFL at $\iota_+$ can be seen clearly in the case of $\theta = 0^{\circ}$, with its peak value reduced in comparison to the 1D case. (This was discussed in section~\ref{sec:geo_cfl}. However, the clear peak at $\theta=0^{\circ},90^{\circ}$ does not significantly persist as the wave angle increases to $\theta = 45^{\circ}$, with the peak becoming substantially flattened. Therefore, the balancing effect that the modification of correction function has on the dispersion of the scheme seems to have a limited scope. Also, it seems there is no other correction function able to achieve the same effect at the intermediate range of angles. Furthermore, the persistence of the peak CFL limit is not seen as the expansion rate and order is varied, with $\iota_+$, in fact, suffering the most appreciable decay in performance as the angle is varied, when compared to the other correction functions.

    To investigate further the impact of varying correction functions on higher dimensional problems, we will consider the projection of the solution into the functional space of FR for linear advection. The method for understanding this is by studying the posedness of the linear operators. The process of projection was outlined in section~\ref{sec:vn_method}. Presented here are the several results showing how the posedness varies with angle, correction function, and order. 
	\begin{figure}
		\centering
		\begin{subfigure}[b]{0.4\linewidth}
			\centering
			\includegraphics[width=\linewidth]{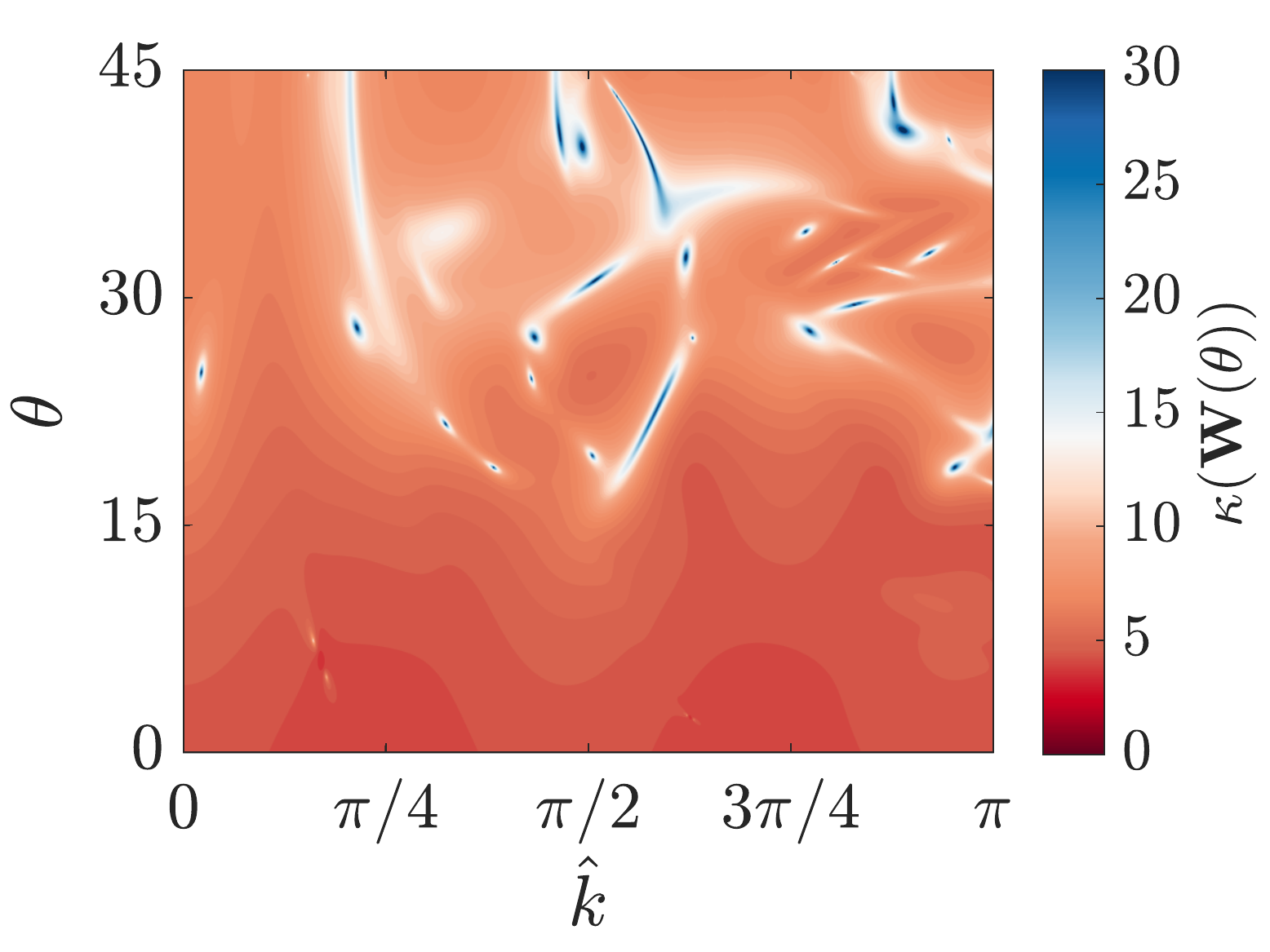}
			\caption{Variation of mode condition number with angle and wavenumber for upwinded 2D FR, $p=2$, with Huynh $g_2$ correction functions.}
			\label{fig:FR2HU_kappa}
		\end{subfigure}
		~
		\begin{subfigure}[b]{0.4\linewidth}
			\centering
			\includegraphics[width=\linewidth]{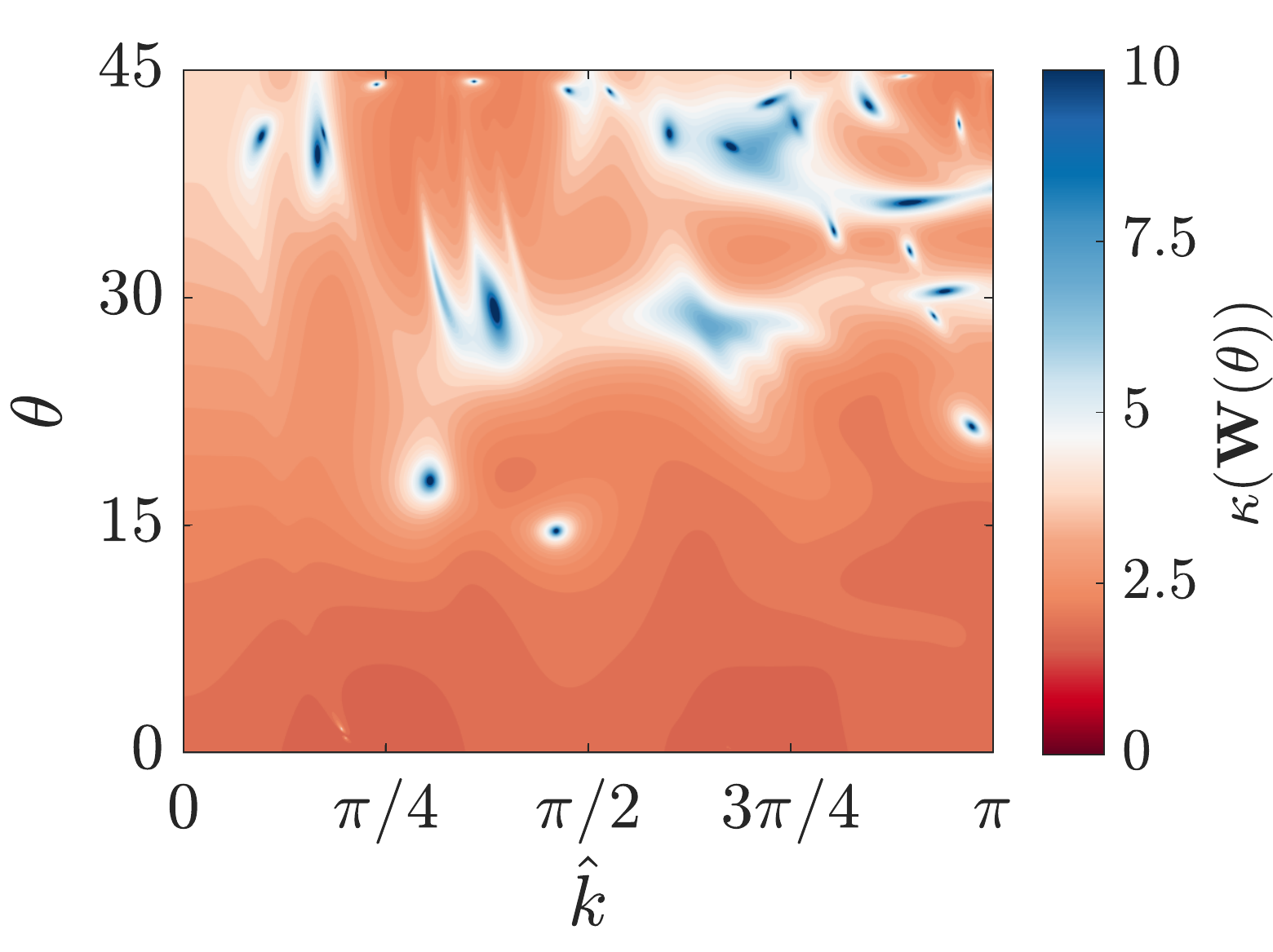}
			\caption{Variation of mode condition number with angle and wavenumber for upwinded 2D FR, $p=2$, with Nodal DG correction functions.}
			\label{fig:FR2DG_kappa}
		\end{subfigure}
		~
		\begin{subfigure}[b]{0.4\linewidth}
			\centering
			\includegraphics[width=\linewidth]{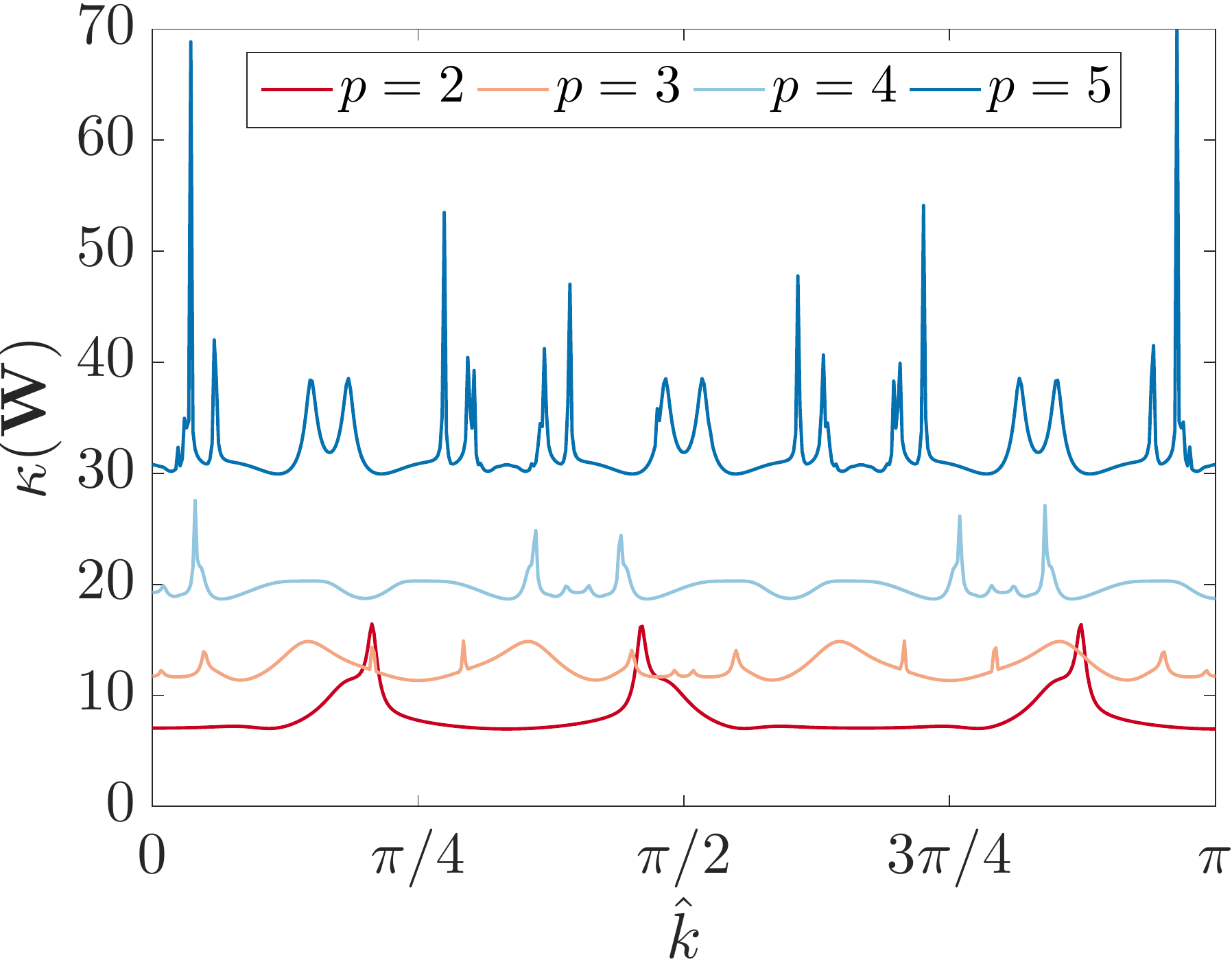}
			\caption{Variation of mode condition number against order, $p$, for upwinded 2D FR, $\theta=45^{\circ}$, with Huynh $g_2$ correction functions.}
			\label{fig:FR_45_kappa}
		\end{subfigure}
		~
		\begin{subfigure}[b]{0.4\linewidth}
			\centering
			\includegraphics[width=\linewidth]{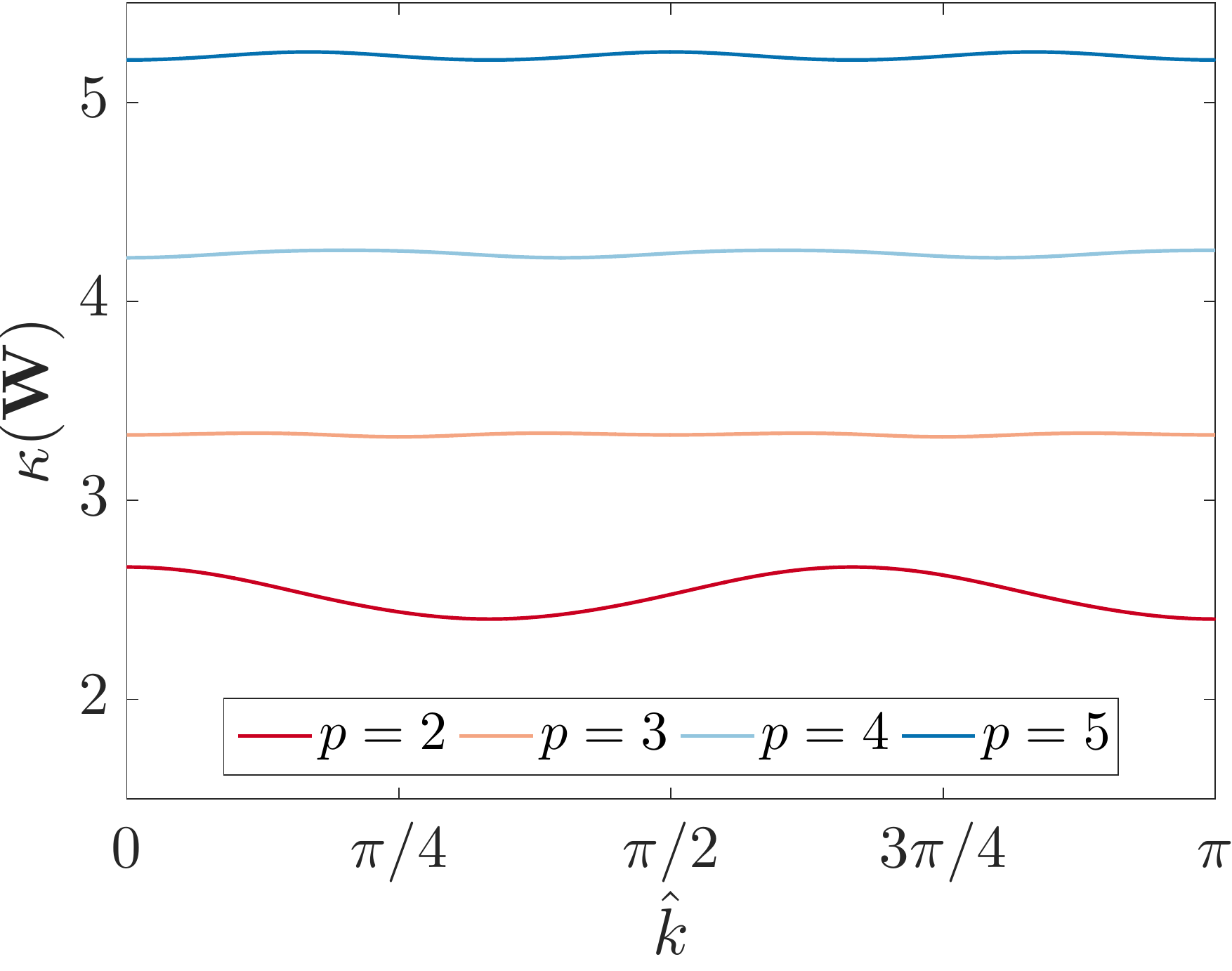}
			\caption{Variation of mode condition number against order, $p$, for upwinded 1D FR, with Huynh $g_2$ correction functions.}
			\label{fig:FR1D_kappa}
		\end{subfigure}
		\caption{Condition number of various linear FR configurations}
		\label{fig:FR_kappa}
	\end{figure}	
	
	Several insights into the different multidimensional behaviour of FR can be gained by studying Fig.~\ref{fig:FR_kappa}. By comparison of Fig.~\ref{fig:FR2HU_kappa}~\&~\ref{fig:FR2DG_kappa} it can be seen that Huynh's $g_2$ correction function causes the projection to be more ill-posed on average compared to Nodal DG methods. Hence, the superconvergent DG recovering scheme of Vincent~\etal~\cite{Vincent2011} has decreased temporal stability. This is because the ill-posedness indicates the sensitivity of the reconstruction to change, a more sensitive reconstruction means that error can result in energy being transferred to more dissipative modes. This point is important as it means to show that ill-conditioning is not directly the mechanism for loss. But the movement of energy to other more dissipative modes is, and the condition number is symptomatic of that.
	
    Furthermore, in both cases as, across all wavenumbers, the condition number can be seen to increase with incident angle. When compared with Fig.~\ref{fig:FR4_2D_HU_CFL} and Fig.~\ref{fig:FR4_RK44_s_y10_x10} it can be seen that for a given correction function the CFL limit reaches its minimum at $\theta = 45^{\circ}$, therefore a stark increase in condition number can also cause a decrease in temporal stability, potentially due to too much transport between modes.

\begin{figure}
		\centering
		\includegraphics[width=0.7\linewidth]{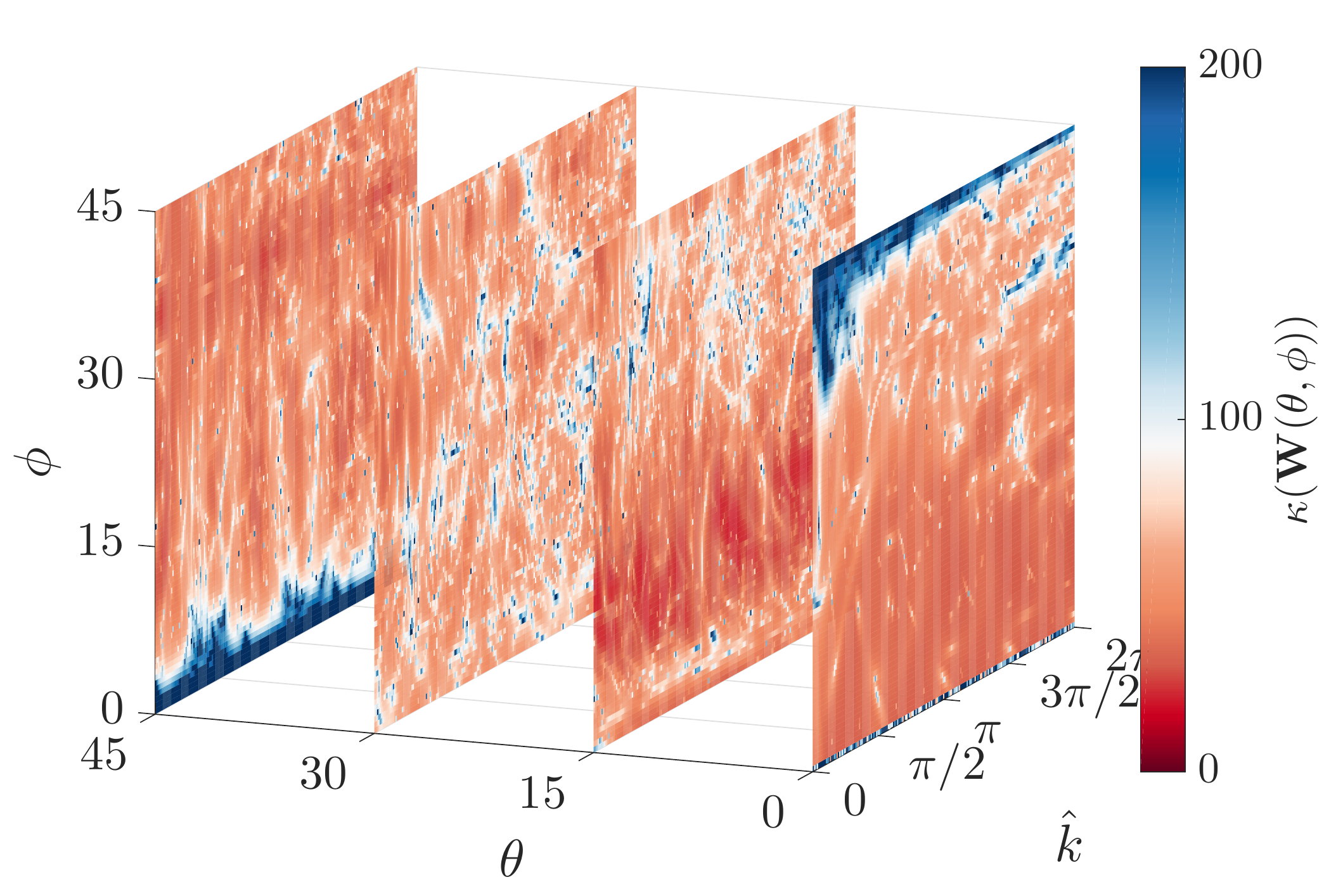}
		\caption{Variation of condition number with wave angles and wavenumber for upwinded 3D FR, $p=3$, $\gamma_x = \gamma_y = \gamma_z = 1$, using Huynh $g_2$ correction functions.}
		\label{fig:FR3_3D_possednes}
	\end{figure}

    A result presented by Trojak~\etal~\cite{Trojak2017a} was that for FR the best points per wavelength (PPW) performance was seen for $p=4$ with the PPW increasing for orders higher than this. A result that is exhibited by Fig.~\ref{fig:FR_45_kappa} is that the condition number for $p=5$ schemes is higher than for lower orders and may have passed a point where increased order is out weighted by inaccuracy in ill-conditioning. This may explain the optimal result seen and is touching on the fundamental problem characterised by Runge's phenomena, that high-order may introduce accuracy through high-order but may also introduce inaccuracy reflected by a high condition. A result that cannot be clearly seen in either Fig.~\ref{fig:FR2HU_kappa}~or~\ref{fig:FR2DG_kappa}, however, was exhibited by FR was for $\theta=0$ the condition number was significantly higher than that found for 1D FR owing to the naturally poor conditioning of a 2D system acting as a quasi-1D one. This is linked to the lower CFL number experienced in this case, as was shown in Fig.~\ref{fig:FR4_2D_HU_CFL}. To show that the results found in this section extend to higher dimensions, the von Neumann analysis was repeated for 3D 'hexi-linear' grids. The primary result of interest is the increased condition number of the functional projection and this can be seen in Fig.\ref{fig:FR3_3D_possednes}. The message is that, as expected, the ill-posedness of the reconstruction increases with the order. While small increases in the condition number can give increased temporal stability, larger increases in condition number tend to act to destabilise the coupled spatial-temporal scheme. 
	
\section{Non-linear Navier-Stokes Equations with Randomised Grids}
    It is common within the CFD community to use the canonical Taylor-Green Vortex~(TGV)~\cite{Taylor1937} test case to assess the numerics of a solver applied to the Navier-Stokes equations with turbulence --- and to that end, there is a plethora of DNS data available for comparison~\cite{Brachet1983,DeBonis2013}. However, this case is quite contrived and ultimately will favour spectral or structured methods due to the Cartesian and periodic domain, whilst also being unrepresentative of engineering flows that are often wall bounded and/or have complex geometries. Hence, we propose linearly deforming the elements of the mesh by jittering the corner nodes to be more representative of real mesh conditions. Importantly, these deformations will introduce cross multiplication into the Jacobian, as well as local regions of expansion and contraction. 

    The initial conditions of the TGV being used here are those of DeBonis~\cite{DeBonis2013}, where the character of the flow is controlled by the non-dimensional parameters defined as:
	\begin{equation}		
		R_e = \frac{\rho_0U_0L}{\mu}, \quad\quad
		P_r = 0.71 = \frac{\mu\gamma R}{\kappa(\gamma-1)}, \quad\quad
		M_a = 0.08 = \frac{U_0}{\sqrt{\gamma R T_0}} 
	\end{equation}
	where we will use the standard set of free-variables for the velocity, density, pressure, and gas characteristics:
	\begin{equation}
		U_0 = 1, \quad\quad \rho_0 = 1, \quad\quad p_0 = 100, \quad\quad R = 1, \quad\quad \gamma = 1.4, \quad\quad L = 1
	\end{equation}
	Here, due to the solver implementation, we use a specific gas constant of unity and hence, to achieve the required Reynolds and Prandtl numbers, the dynamic viscosity and thermal conductivity can be set appropriately. 
	
\begin{figure}
		\centering
		\begin{subfigure}[b]{0.3\linewidth}
			\centering
			\includegraphics[width=\linewidth]{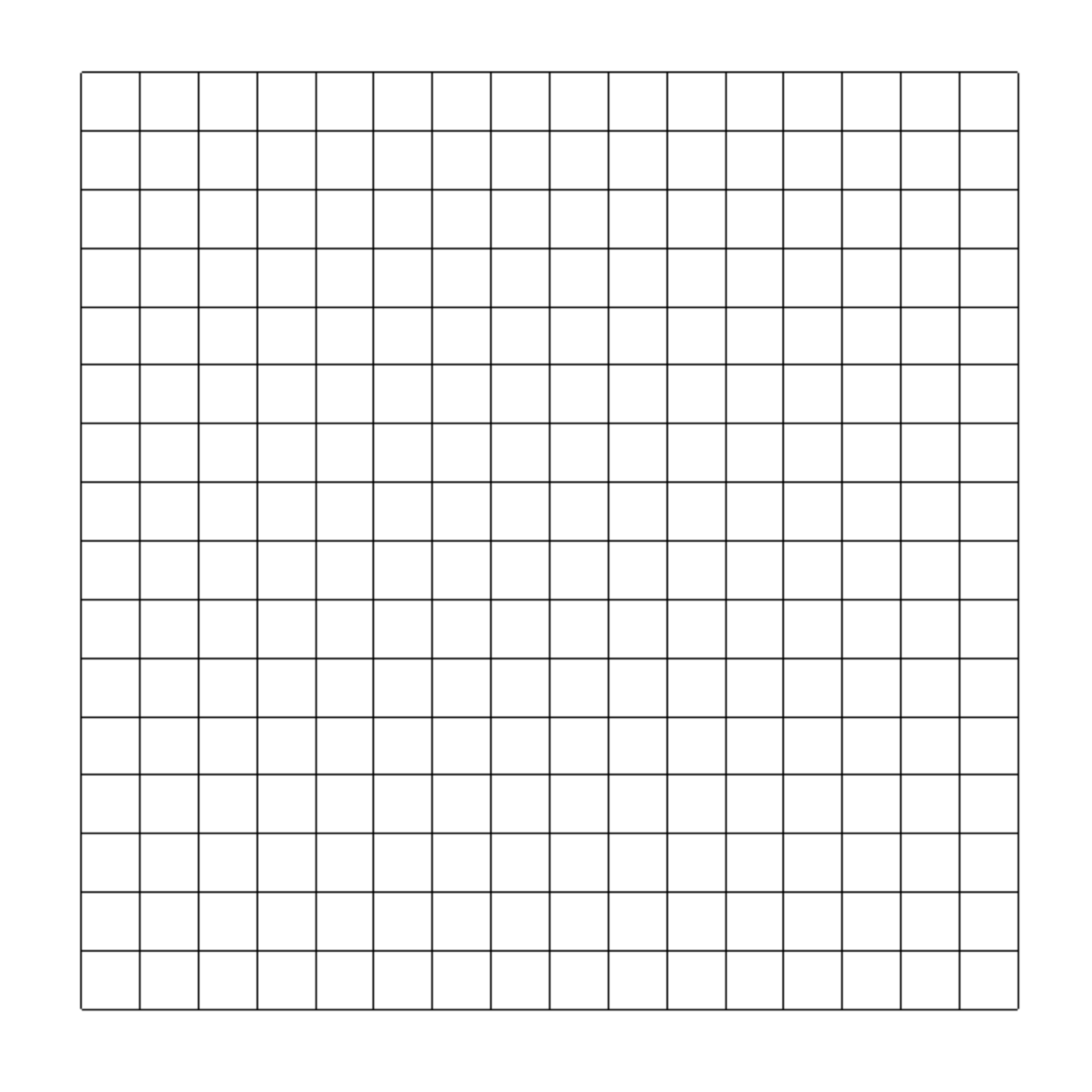}
			\caption{$q_h=\sqrt{\pi/6}\approx0.7236$}
			\label{fig:mesh1}
		\end{subfigure}
		~
		\begin{subfigure}[b]{0.3\linewidth}
			\centering
			\includegraphics[width=\linewidth]{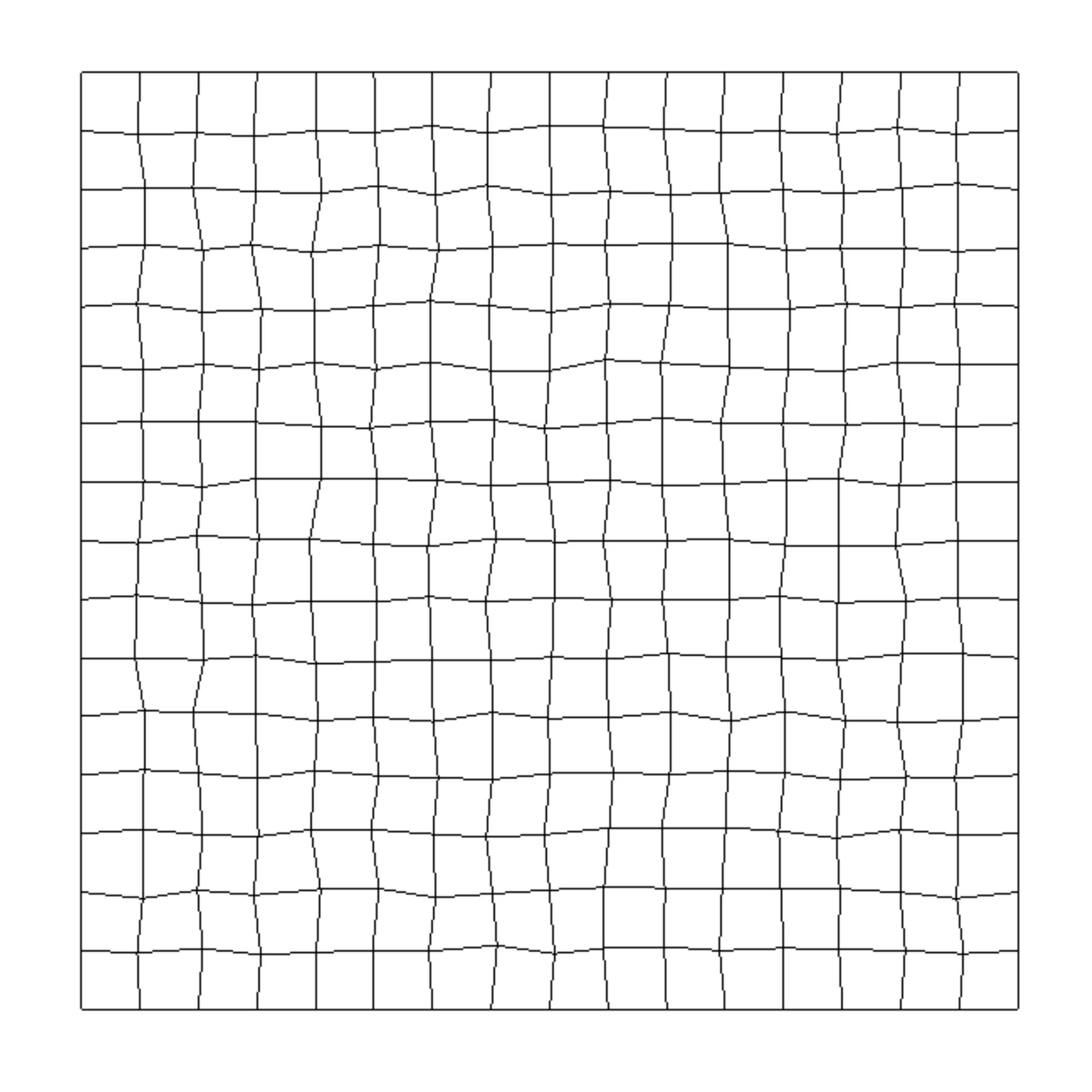}
			\caption{$q_h=0.7201$}
			\label{fig:mesh2}
		\end{subfigure}
		~
		\begin{subfigure}[b]{0.3\linewidth}
			\centering
			\includegraphics[width=\linewidth]{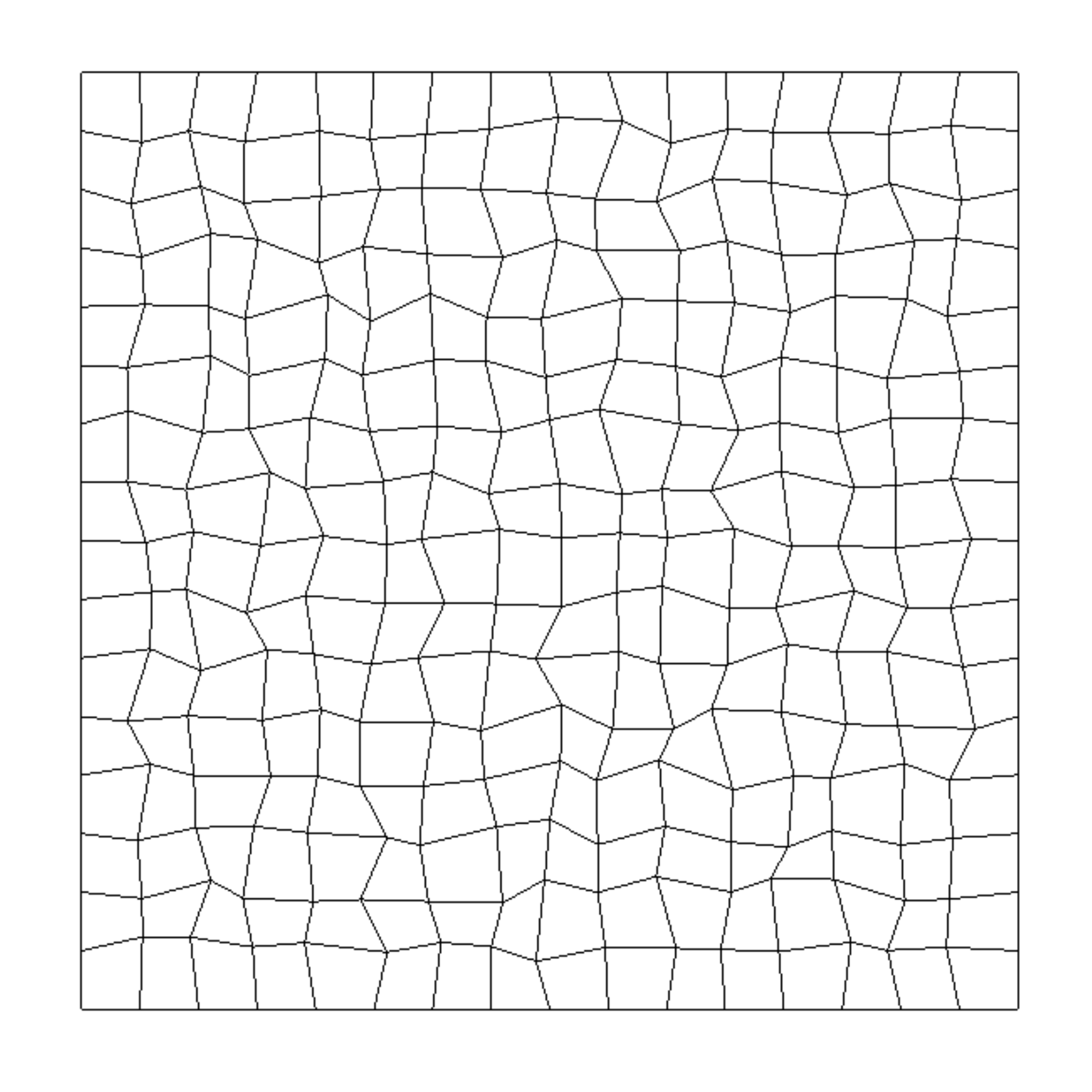}
			\caption{$q_h=0.7016$}
			\label{fig:mesh3}
		\end{subfigure}
		\caption{Example slices through a 3D  hexahedral mesh to illustrate the mesh quality metric.}
		\label{fig:mesh_plots}
	\end{figure}

	As has been stated, we will take the uniform periodic mesh on the domain $\mathbf{\Omega}\in[-\pi,\pi]^3$, and jitter the corner nodes of the elements that are interior to the domain. The amount of jitter is calculated using a time seeded random number shifted to be centred about zero and scaled by a global factor between zero and unity. The scaling factor is such that zero gives a uniform mesh and unity could lead to edges of zero length. After jittering, the solution points are then linearly positioned within the element using the thin plate spline radial basis function together with the mapping from the uniform to jittered corner nodes. This gives a linear mapping of uniform solution points to solution points within the jittered elements. Finally, a quality metric is needed to describe, in a single number, the relative quality of the meshes produced. We opted for a volume ratio shape factor, slightly redefined as:
	\begin{equation}
		q_h = \frac{6\sqrt{\pi}V_h}{S_h^{3/2}}
	\end{equation}
	where $S_h$ is the surface area of the hexahedral element and $V_h$ is the volume of  the hexahedral elements. The quality metric, $q_h$, is then defined as the ratio of the volume of the element to the volume of a sphere with the same surface area, with $q_h=\sqrt{\pi/6}$ for a perfect cube. To put this parameter into context, some example meshes are shown in Fig.~\ref{fig:mesh_plots}.
	
	The statistics that will be studied here are the decay of the kinetic energy and the enstrophy decay rate, which are defined respectively as:
	\begin{align}
		-\dx{E_k}{t} &= -\frac{1}{2\rho_0|\mathbf{\Omega}|}\dx{}{t}\int_{\mathbf{\Omega}}\rho(u^2+v^2+w^2)\mathrm{d}\mathbf{x} \\
		\epsilon &= \frac{\mu}{\rho_0^2|\mathbf{\Omega}|}\int_{\mathbf{\Omega}}\rho(\pmb{\omega}\cdot\pmb{\omega})\mathrm{d}\mathbf{x}
	\end{align}
	where $\pmb{\omega}=\nabla\times[u,v,w]^T$ is vorticity and $|\mathbf{\Omega}|$ is the domain volume.
	
	\begin{figure}
		\centering
		\begin{subfigure}[b]{0.45\linewidth}
			\centering
			\includegraphics[width=\linewidth]{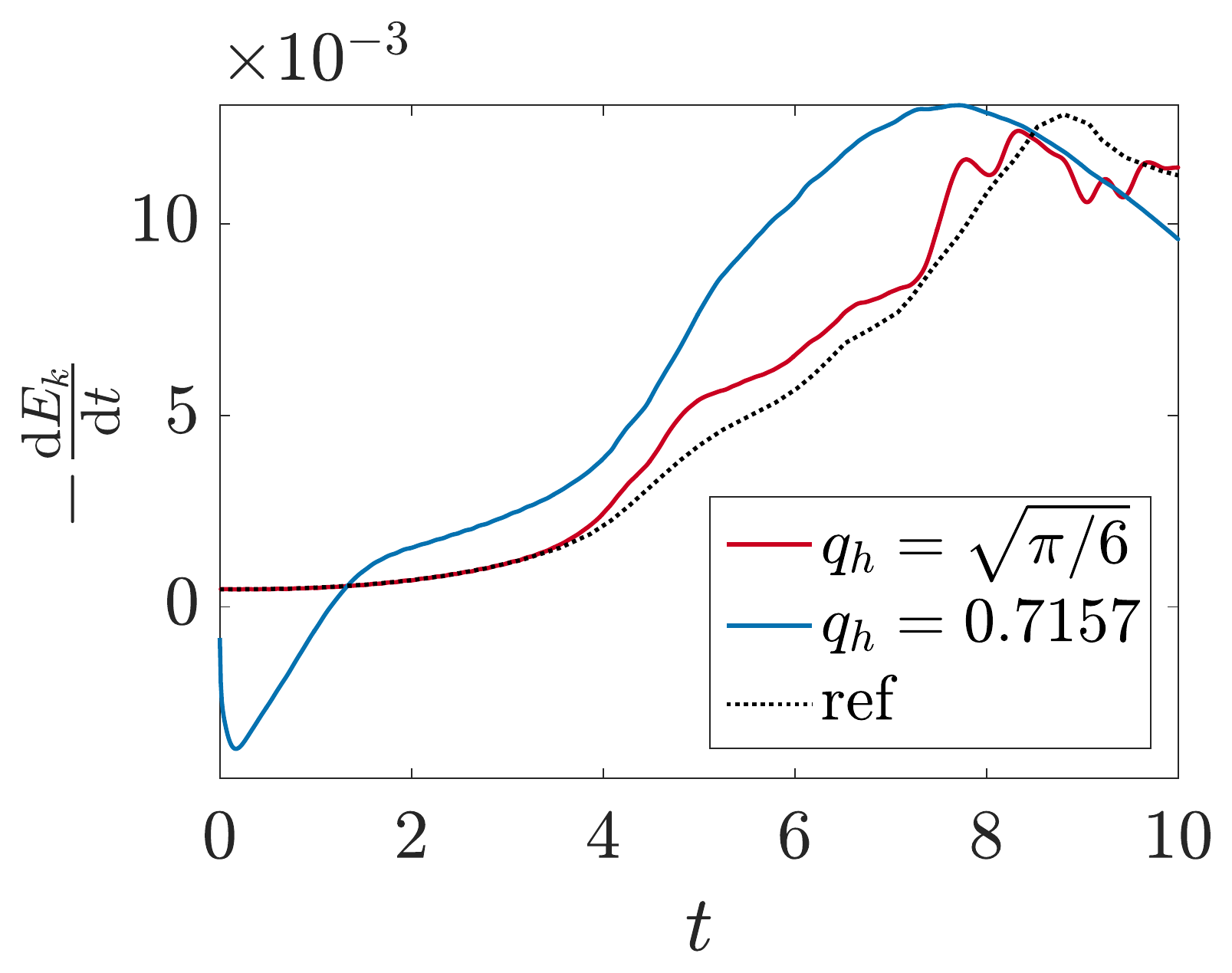}
			\caption{Selected turbulent kinetic energy dissipation.\\ \quad}
			\label{fig:FRHU_p2_tgv_j_diss}
		\end{subfigure}
		~
		\begin{subfigure}[b]{0.47\linewidth}
			\centering
			\includegraphics[width=\linewidth]{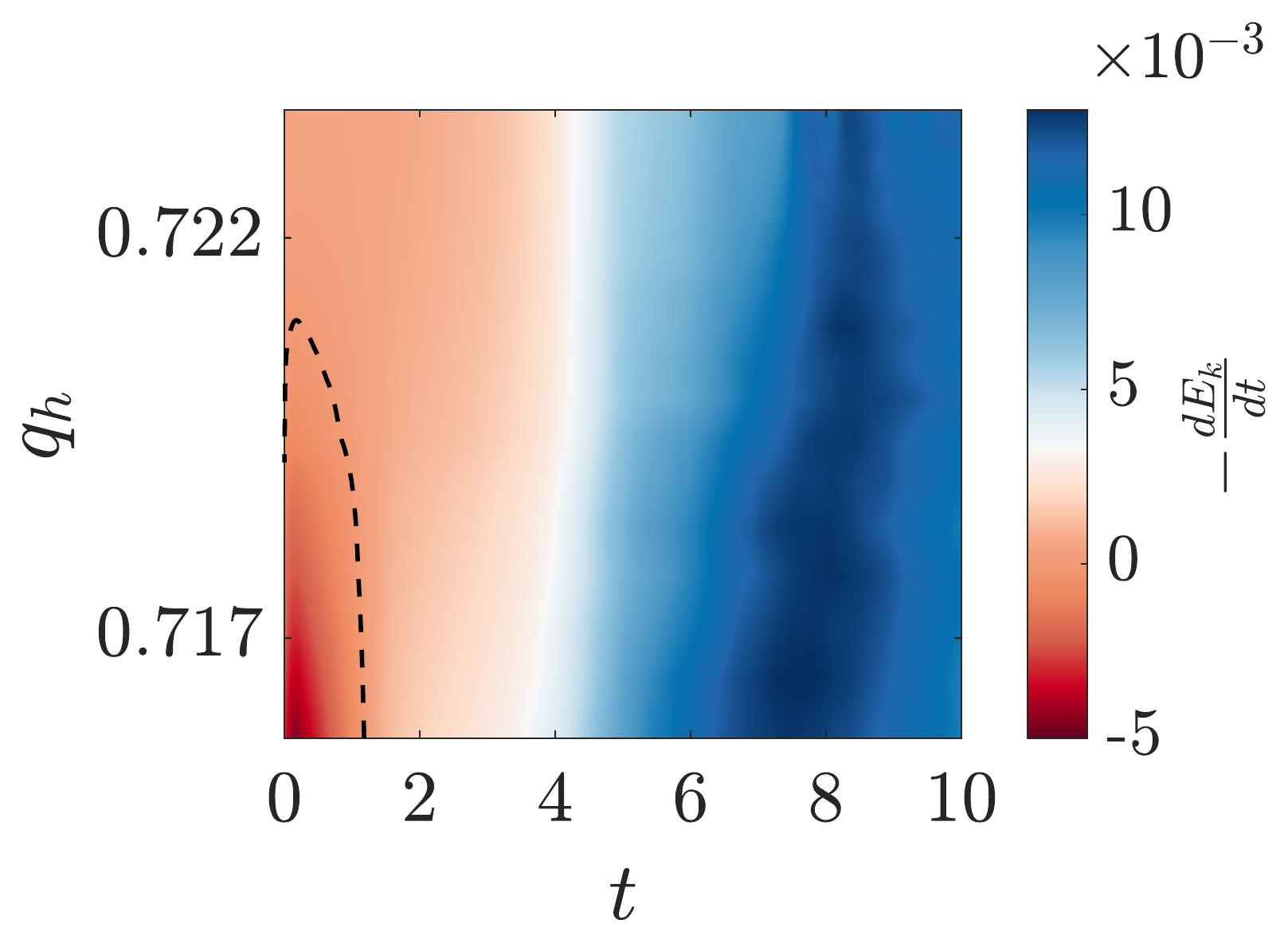}
			\caption{Variation of turbulent kinetic energy dissipation with jitter. Dashed contour at zero dissipation.}
			\label{fig:FRHU_p2_tgv_j}
		\end{subfigure}
		\caption{Effect of jitter on turbulent kinetic energy dissipation of the TGV ($R_e=1600$) for FR, $p=2$, with Huynh $g_2$ correction functions on a $40^3$ DoF mesh. Explicit time step size is $\Delta t = 1\times10^{-3}$.}
		\label{fig:FRHU_p2_TGV}
	\end{figure}
	
	\begin{figure}
		\centering
		\begin{subfigure}[b]{0.45\linewidth}
			\centering
			\includegraphics[width=\linewidth]{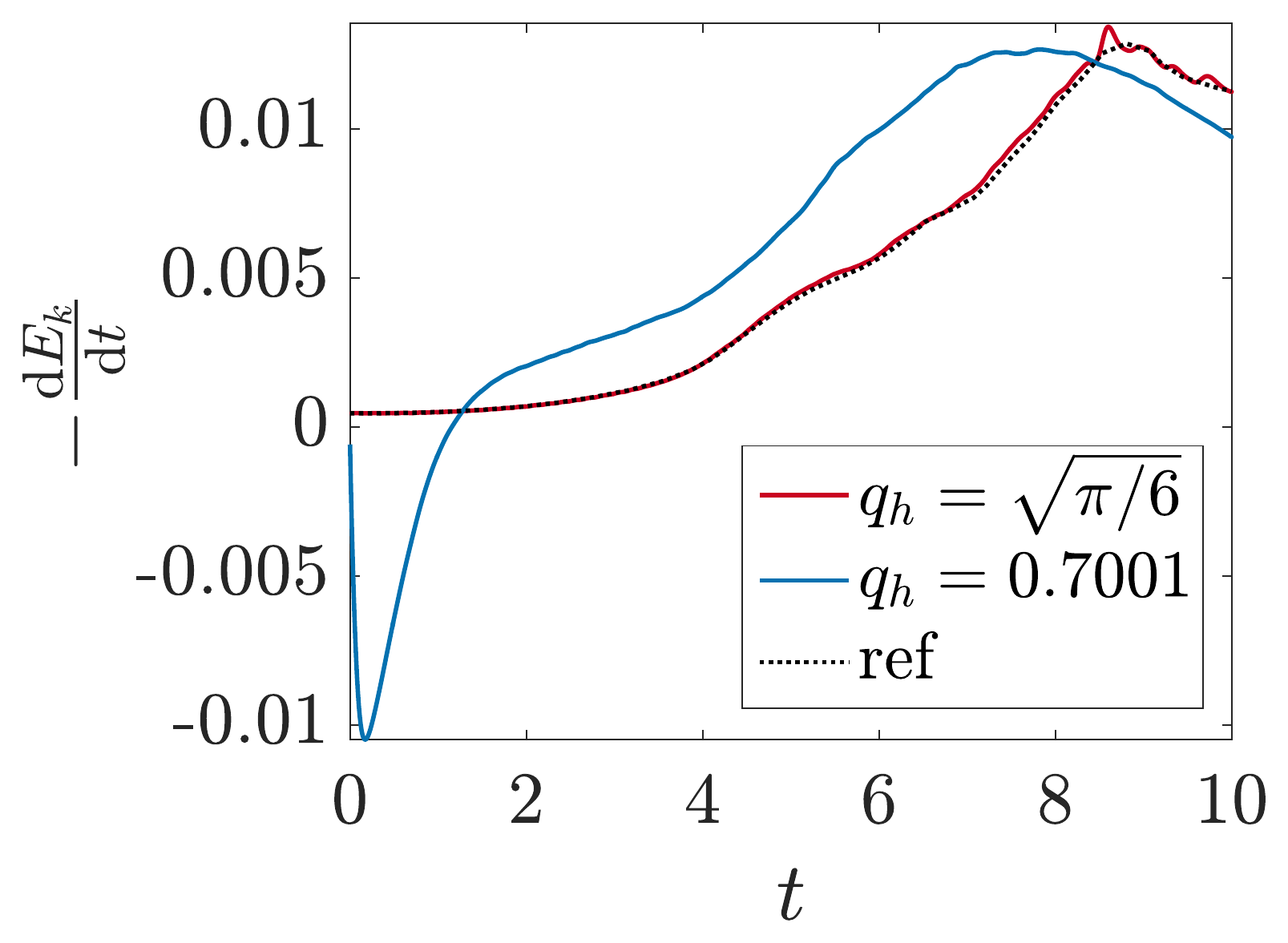}
			\caption{Selected turbulent kinetic energy dissipation.\\ \quad}
			\label{fig:FRHU_p4_tgv_j_diss}
		\end{subfigure}
		~
		\begin{subfigure}[b]{0.47\linewidth}
			\centering
			\includegraphics[width=\linewidth]{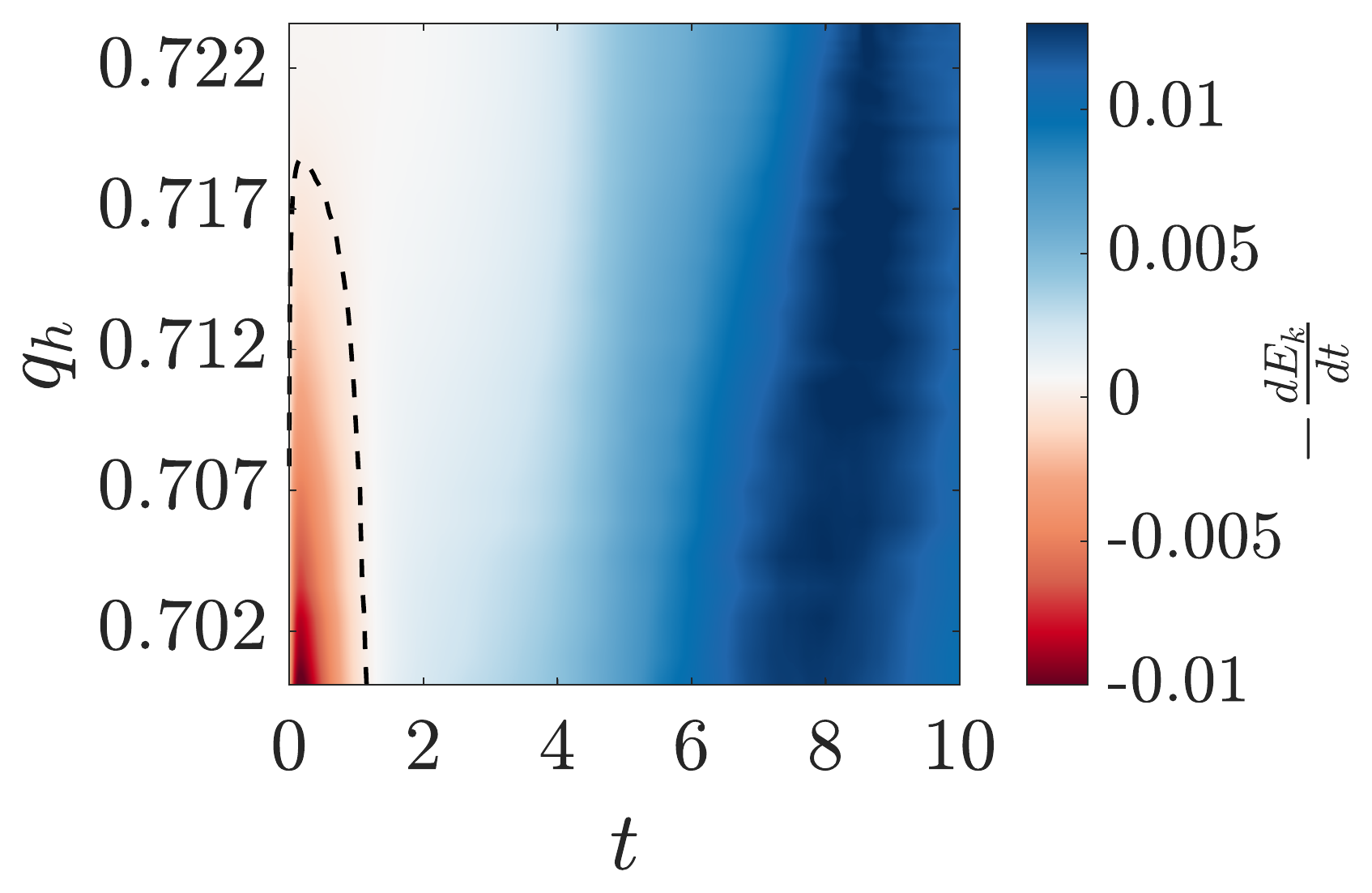}
			\caption{Variation of turbulent kinetic energy dissipation with jitter. Dashed contour at zero dissipation.}
			\label{fig:FRHU_p4_tgv_j}
		\end{subfigure}
		\caption{Effect of jitter on turbulent kinetic energy dissipation of the TGV ($R_e=1600$) for FR, $p=4$, with Huynh $g_2$ correction functions on a $120^3$ DoF mesh. Explicit time step is $\Delta t = 1\times10^{-3}$.}
		\label{fig:FRHU_p4_TGV}
	\end{figure}
	
	Fig.~\ref{fig:FRHU_p2_TGV}~\&~\ref{fig:FRHU_p4_TGV} shows the first of these results. First looking at Fig.~\ref{fig:FRHU_p2_tgv_j_diss}~\&~\ref{fig:FRHU_p4_tgv_j_diss} which shows two specific dissipation curves for an uniform and jittered mesh. At the beginning of the simulation, there is a clear time at which the global energy increases. 

    Extending these runs to cover multiple grid qualities, Fig.~\ref{fig:FRHU_p2_tgv_j}~\&~\ref{fig:FRHU_p4_tgv_j}, it is observed that as the grid quality decreases a region where turbulent kinetic energy increases soon emerges. As time progress, energy dissipation is again seen and the point of peak dissipation arrives early, moving from $t\approx 8.5$ to $t\approx 7.5$. The same behaviour is seen for both $p=2$ and $p=4$. From comparison of $p=2$ and $p=4$, it seems that $p=4$ is slightly more robust to grid deformation, as $p=4$ was able to run at $q_h\approx 0.7$. Whereas for $p=2$, $q_h$ could not be reduced much below $0.717$ for $120^3$ DoF without completely diverging.    

    The explanation of this is believed to be that initially the regions in the mesh that are locally expanding cause an increase in the energy due to the positive dissipation. This was discussed in Section~\ref{sec:dissdisp}, that for the linear advection equation, dissipation is positive at low wavenumbers for expanding grids and negative for contracting grids at the same wavenumber. It is thought that as the simulation progresses, the energy cascade of large scales to small scales then means that more of the solution lies in the more dissipative higher wavenumber region for both expanding and contracting grids. Therefore, the net dissipation at a later time is higher than the uniform case and hence the peak dissipation is earlier. This is consistent with a lower $R_e$ and hence higher global dissipation. We can conclude that the stability of this case is brought about by the physics of the Navier-Stokes and cascade of energy from low to high wavenumbers, which sidesteps the problem of positive dissipation of low wavenumbers on expanding meshes. This result is interesting as it is in slight contradiction to the result of Trojak~\etal~\cite{Trojak2017a}, where FR was found to be more robust to grid stretching than second order FV for Euler's equations. However, in that investigation, the Isentropic Convecting Vortex was considered where there was a large convective velocity. This may expose slightly different properties, although a full explanation is not known.
	
	\begin{figure}
		\centering
		\begin{subfigure}[b]{0.45\linewidth}
			\centering
			\includegraphics[width=\linewidth]{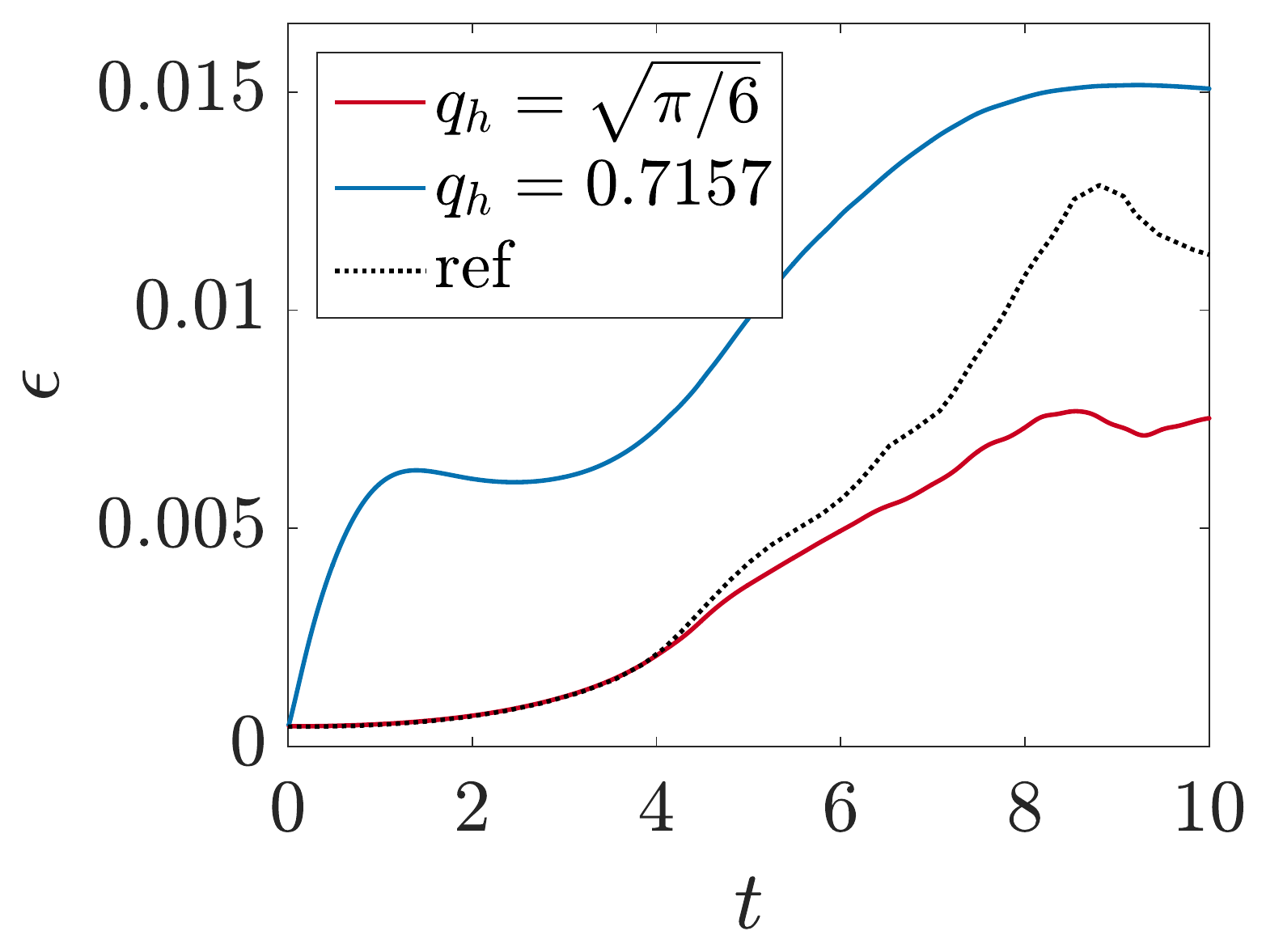}
			\caption{$p=2$}
			\label{fig:FRHU_p2_tgv_e2}
		\end{subfigure}
		~
		\begin{subfigure}[b]{0.47\linewidth}
			\centering
			\includegraphics[width=\linewidth]{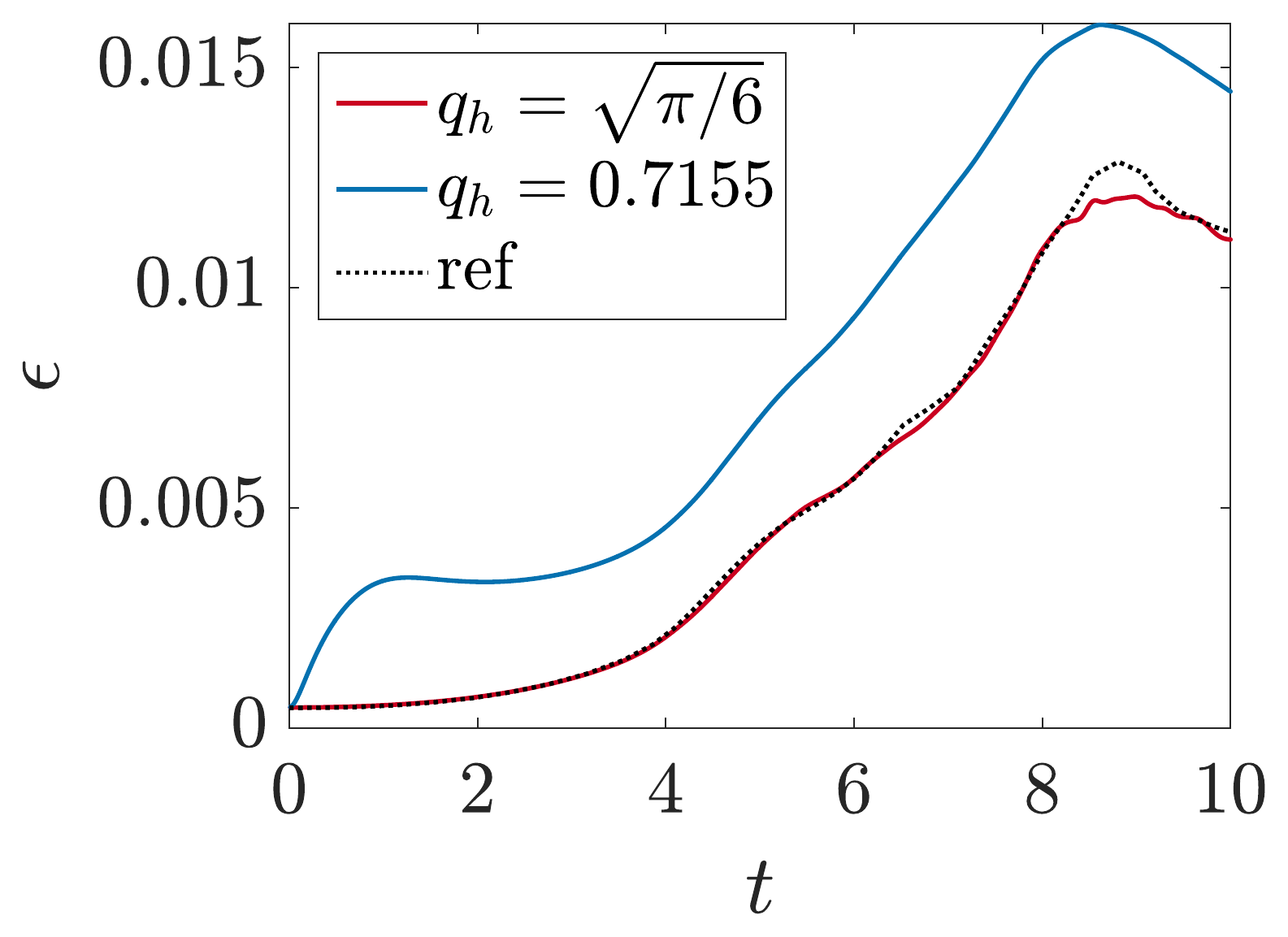}
			\caption{$p=4$}
			\label{fig:FRHU_p4_tgv_e2}
		\end{subfigure}
		\caption{Comparison of TGV enstrophy for $120^3$ degree of freedom grid with similar $q_h$.}
		\label{fig:FRHU_TGVe2}
	\end{figure}
	
	Studying the effect of jittered grids on enstrophy, shown in Fig.~\ref{fig:FRHU_TGVe2}, it is clear that as the grid is stretched the enstrophy increases. This is indicative of an increase in the vorticity, with the rise occurring within $t=0-1$. This is consistent with the assertion that the energy is being added in the small scales, as at this time there are only small scales present. Furthermore, if the energy were only added to the bulk flow, i.e. $k=0$, it would not be seen in the enstrophy. After the initial increase, the enstrophy returns to following the trend of uniform case. However, in the case of $p=2$, Fig.~\ref{fig:FRHU_p2_tgv_e2}, a larger initial increase is seen followed by a wider peak. The wider peak is similar in character to that of the uniform case and is due to the grid being mildly under-resolved in the $p=2$ case relative to the DNS.

\begin{figure}
		\centering
		\begin{subfigure}[b]{0.4\linewidth}
			\centering
			\includegraphics[width=\linewidth]{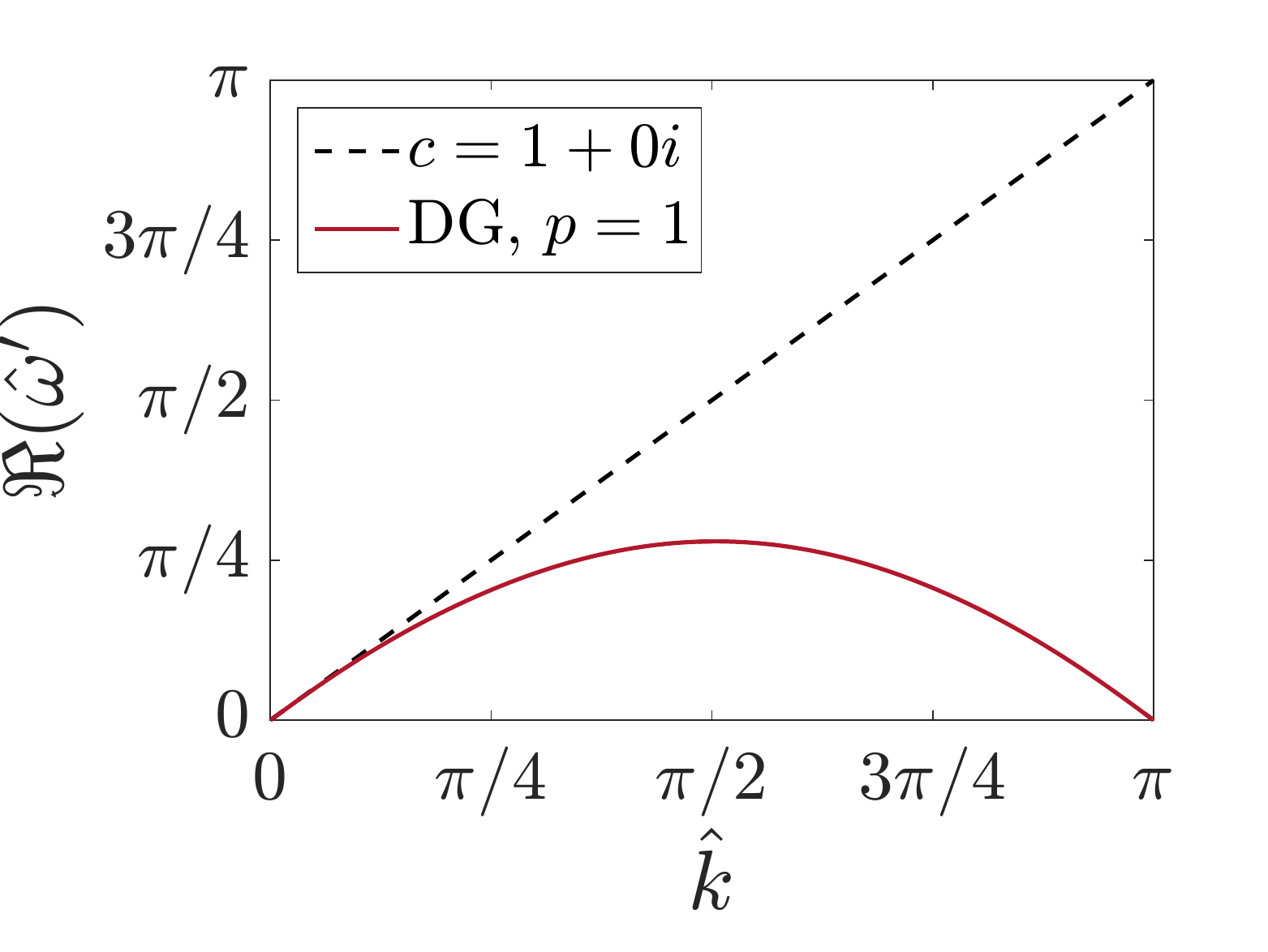}
			\caption{Dispersion \\ \quad}
			\label{fig:FR1_R}
		\end{subfigure}
		~
		\begin{subfigure}[b]{0.4\linewidth}
			\centering
			\includegraphics[width=\linewidth]{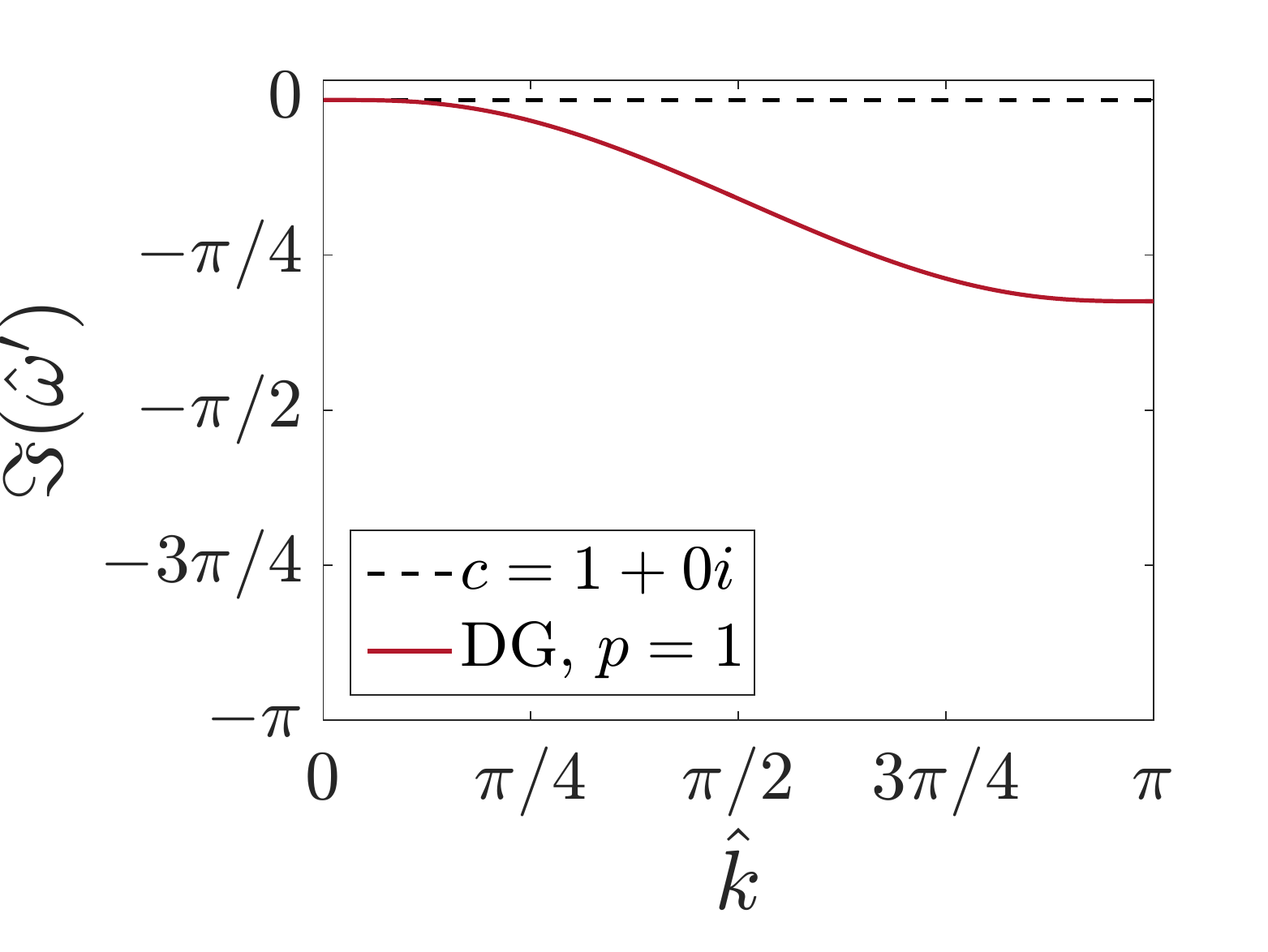}
			\caption{Dissipation}
			\label{fig:FR1_I}
		\end{subfigure}
		\caption{Dispersion and Dissipation relation for 1D upwinded FR, $p=1$, with DG correction function.}
		\label{fig:FR1}
	\end{figure}

    To provide some reference as to how FR performs relative to an established method we will use an edge-based Finite Volume (FV) method for comparison. The FV method is a standard central second order method with L2 Roe smoothing~\cite{Osswald2016} for stabilisation, which has been validated previously~\cite{Scillitoe2016}. The particular FR scheme used in this comparison is $p=1$, giving second order, the same as the FV scheme. However, this puts FR at a significant disadvantage as its numeric characteristics at low order are particularly poor. For example, consider the dispersion and dissipation relations in Fig.~\ref{fig:FR1}, which, by comparison to the result of Lele~\cite{Lele1992}, show that FR has noticeably lower resolving abilities when compared against a second order FD scheme. 

\begin{figure}
		\centering
		\begin{subfigure}[b]{0.45\linewidth}
			\centering
			\includegraphics[width=\linewidth]{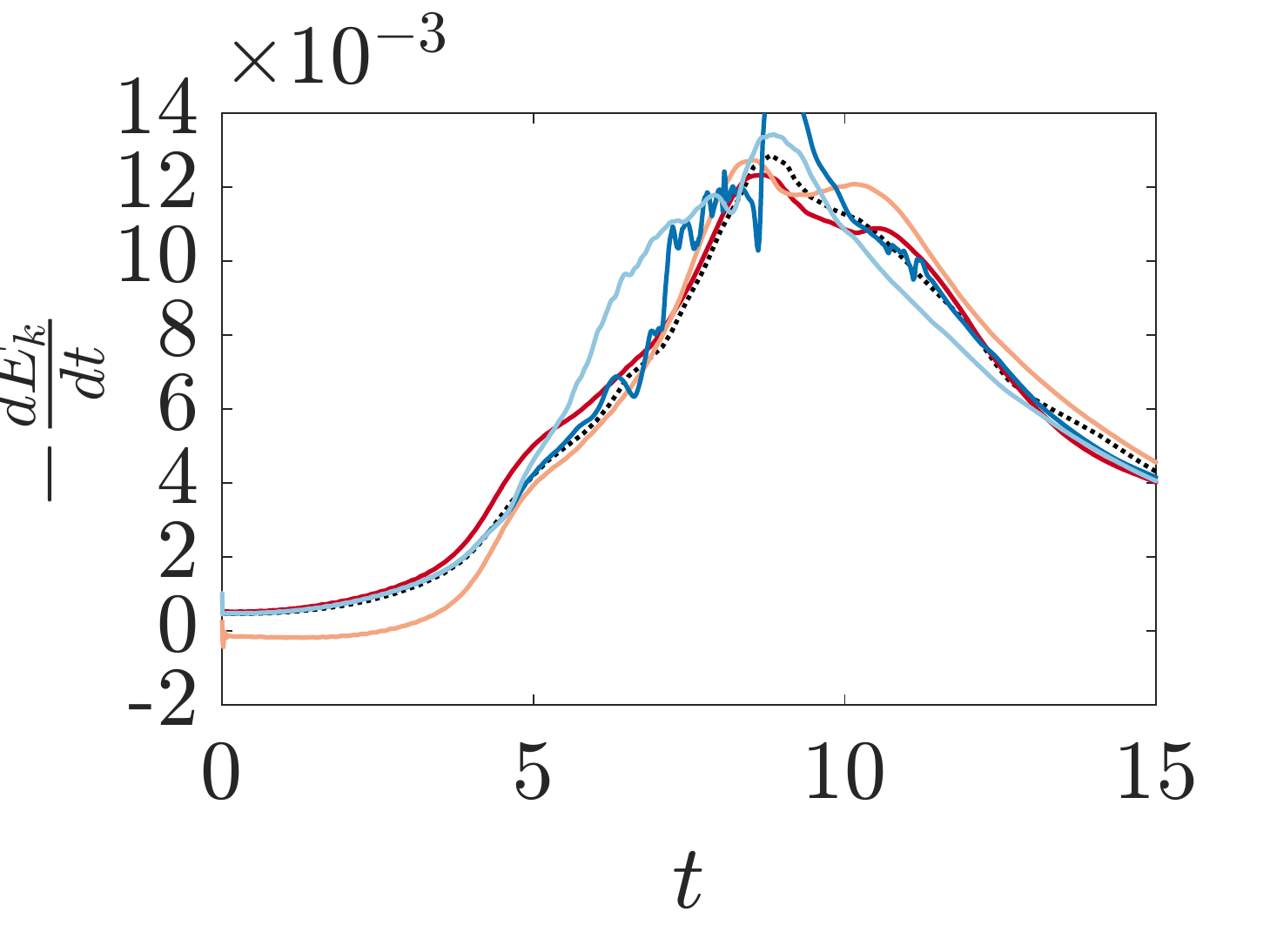}
			\caption{Kinetic energy dissipation}
			\label{fig:TGV_comp_e1}
		\end{subfigure}
		~
		\begin{subfigure}[b]{0.45\linewidth}
			\centering
			\includegraphics[width=\linewidth]{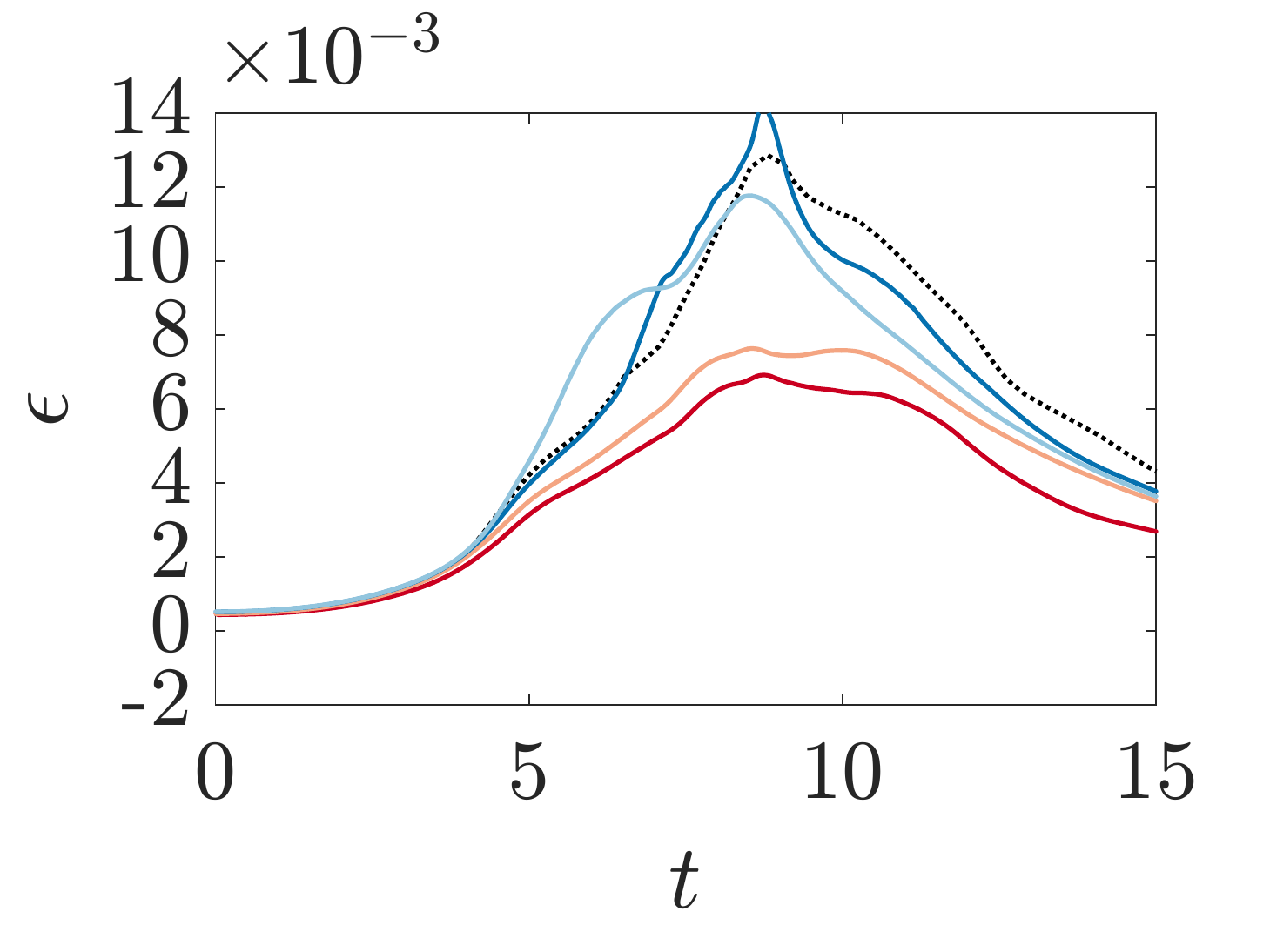}
			\caption{Enstrophy based dissipation}
			\label{fig:TGV_comp_e2}
		\end{subfigure}
		~
		\begin{subfigure}[b]{0.3\linewidth}
			\centering
			\includegraphics[width=\linewidth]{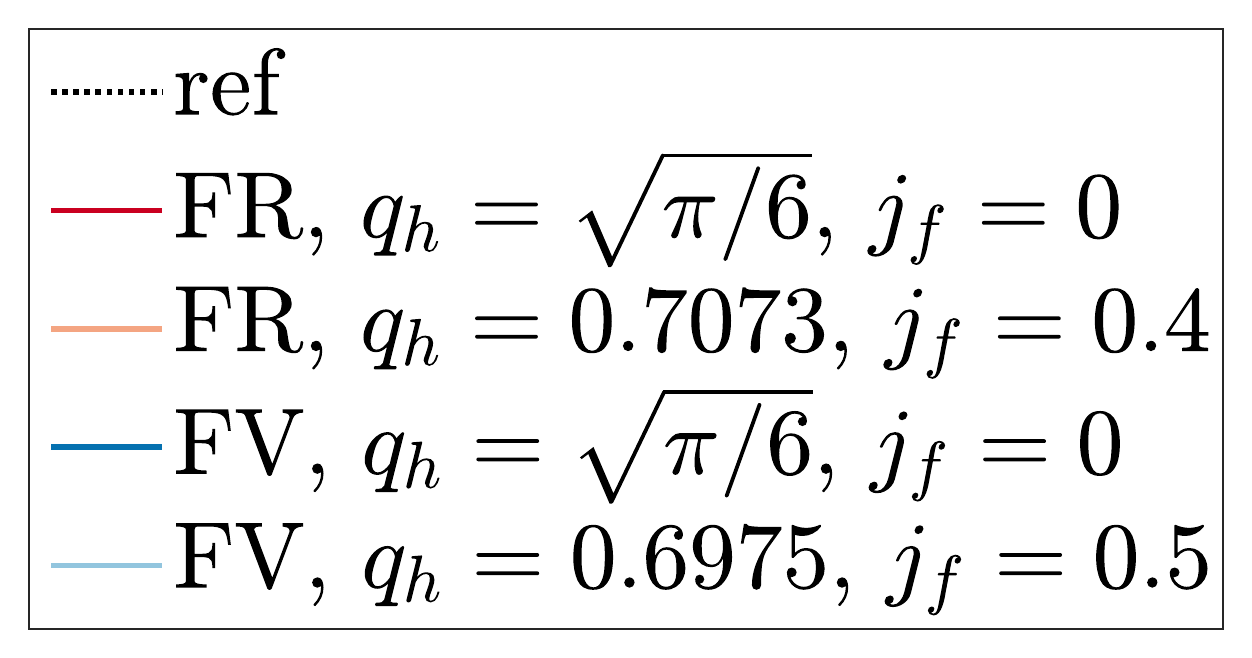}
		\end{subfigure}
		\caption{Comparison of FR, $p=1$ with DG correction functions with a second order central FV scheme with L2 Roe smoothing both with $170^3$ degrees of freedom and $\Delta t\approx 5\times10^{-4}$. A reference DNS solution is provided from \cite{Brachet1983}.}
		\label{fig:TGV_comp}
	\end{figure}

\begin{figure}
		\centering
		\begin{subfigure}[b]{0.45\linewidth}
			\centering
			\includegraphics[width=\linewidth]{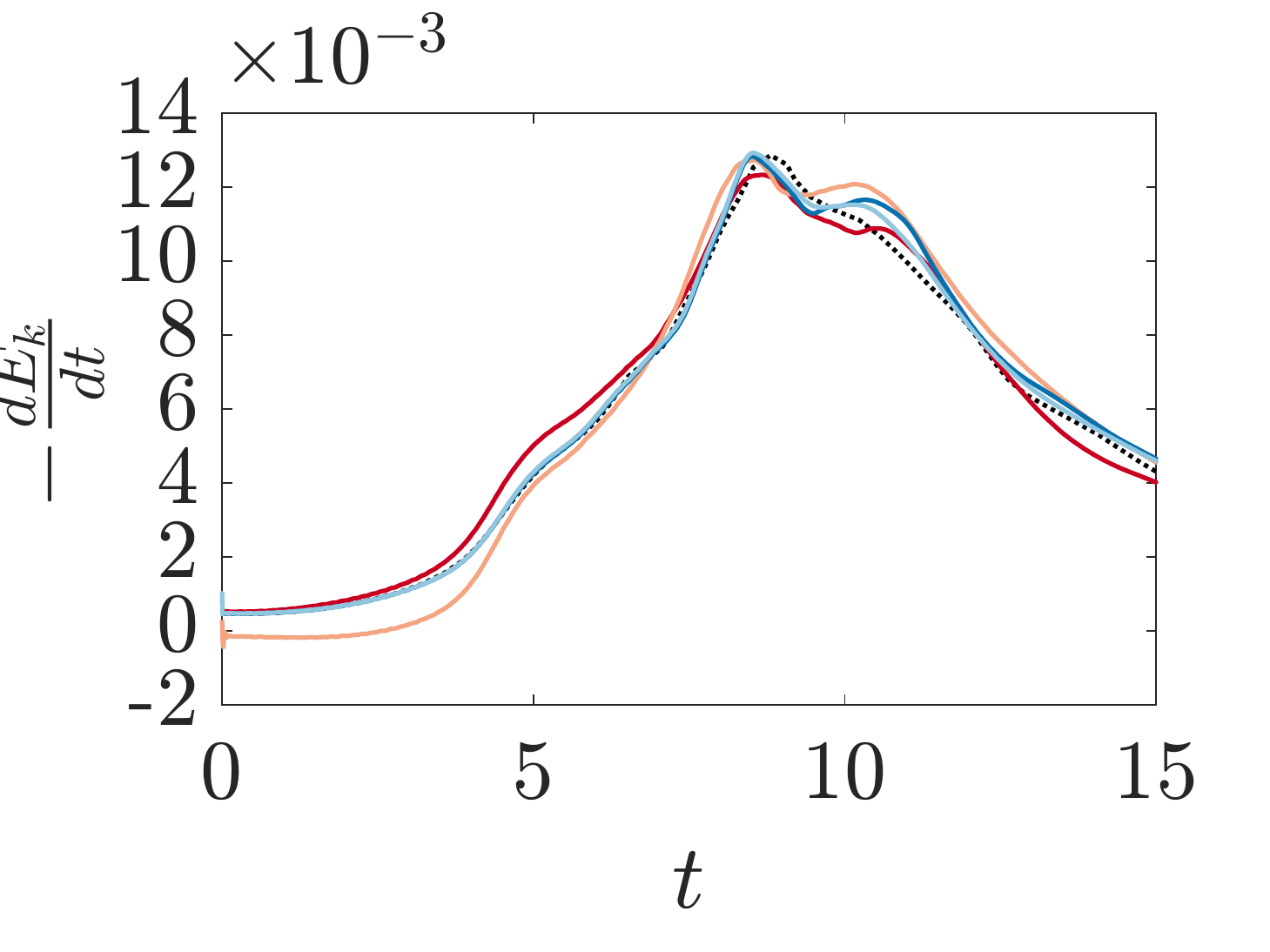}
			\caption{Kinetic energy dissipation}
			\label{fig:TGV_comp_kep_e1}
		\end{subfigure}
		~
		\begin{subfigure}[b]{0.45\linewidth}
			\centering
			\includegraphics[width=\linewidth]{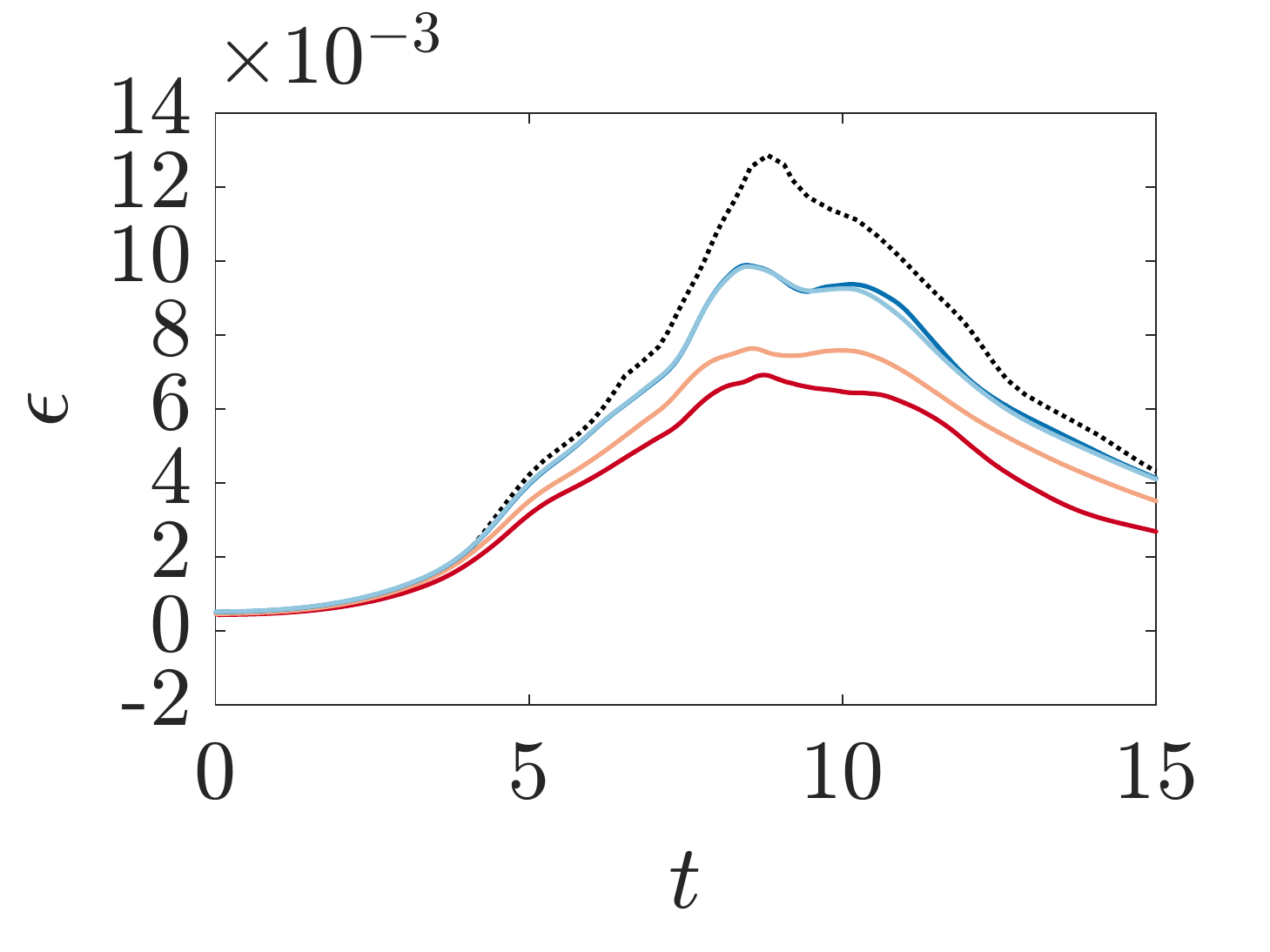}
			\caption{Enstrophy based dissipation}
			\label{fig:TGV_comp_kpe_e2}
		\end{subfigure}
		~
		\begin{subfigure}[b]{0.3\linewidth}
			\centering
			\includegraphics[width=\linewidth]{TGV_comp_legend.pdf}
		\end{subfigure}
		\caption{Comparison of FR, $p=1$ with DG correction functions with a second order KEP FV scheme with L2 Roe smoothing both with $170^3$ degrees of freedom and $\Delta t\approx 5\times10^{-4}$. A reference DNS solution is provide from \cite{Brachet1983}.}
		\label{fig:TGV_kep_comp}
	\end{figure}

	With this in mind, we present the results of tests on various jittered grids with a total of $170^3$ degrees of freedom in Fig.~\ref{fig:TGV_comp}. For the uniform case, the enstrophy clearly shows that FR is underresolved compared to FV, which is also shown by a slightly increased rate of dissipation earlier --- indicating that the implicit filter is too narrow. If we now consider the effect of jittering, several things may be concluded. Firstly, we were unable to run FR with $j_f=0.5$ as the simulation quickly became unstable. Secondly, for $-dE_k/dt$ it seems that the peak value is less sensitive with FR than with FV, with central FV seeing some very large amplitude oscillation in $-dE_k/dt$. This is likely to be rooted in the central differencing at the interfaces, as if we change to a kinetic energy preserving formulation~\cite{Jameson2008,WatsonTucker2015}, as is displayed in Fig.~\ref{fig:TGV_kep_comp}, these oscillations are removed and the sensitivity to jitter is reduced. The enstrophy (Fig.~\ref{fig:TGV_comp_kpe_e2}) seems to indicate that a large amount of what seemed to be resolved energy may have in fact been dispersion induced fluctuations. However, in both cases FV was able to run with grids up to $j_f=0.9$ and $q_h=0.6382$ --- not shown --- it appears that in these cases the added stability of the smoothing has greatly helped FV, especially in the central difference case where running without smoothing caused the case to fail even at low levels of jitter.
	
	Before moving on, it must be noted that in FR we once again see a dip in $dE_k/dt$, which is also present to a lesser extent in both versions of FV tested. This is again linked to the instability caused by locally expanding grids, but in the case of FV the dip is smaller and is aided by the use of smoothing. The conclusion for this comparison is that FV is somewhat resilient to degradation to the mesh quality, with the resilience coming from smoothing for central differenced FV. This allowed highly warped meshes to be run, but at the expense of accuracy with excess dissipation affecting the solution. KEP was found to be far more resilient and could even run without smoothing. FR, when run at a low order, was unsuited to this problem, but did see less degradation compared to central FV and there is the potential for $p=1$ FR to equally benefit on poorer quality meshes from smoothing via a different Riemann solver, such as that of Roe~\cite{Roe1981}. 

    The result shown in this section clearly highlight a key issue of FR, that it works on the strong conservative form of the equation. Consequently, it relies on the accurate interpolation and accurate calculation of gradients to make sure that the scheme does not stray far from conservation. However, this implementation does not strictly enforce conservation and hence when inaccuracy in the gradient or interpolation is introduced there is no remedial action taken. Some recent work that aims to fix this issue, see Abe~\etal~\cite{Abe2015,Abe2016,Abe2018}, where the common interface values are used to enforce conservation. Following from the results of the present work, these extensions to FR should be considered essential.
	
\section{Conclusions}
\label{sec:conclusions}
Through this work, we have presented a theoretical extension of the FR von Neumann analysis to higher dimensions. This allowed as to understand the character of the dispersion and dissipation relations of FR as the incident angle of a wave was varied. Differences were noted between the behaviours of FR and finite differencing methods, primarily that FR saw a lower variation in character with the angle of incidence. The effect of higher dimensionality on the CFL limit was also found, with higher dimensionality causing a reduction in the CFL limit. 

    Investigations were then performed on deformed meshes and theoretically, the same behaviour was seen in two dimensions as in one. Specifically, that expanding meshes cause instability and contracting meshes cause excess dissipation, however, when coupled together the effects can act to cancel each other out. Numerical experiments were subsequently performed using the Taylor-Green vortex case, but with the element corner node positions jittered. Tests showed $5^{\mathrm{th}}$ order to be more resilient to poor quality meshes than $3^{\mathrm{rd}}$ order, but in both cases, the effect of localised regions of expansion are thought to be responsible for an initial increase in the kinetic energy of the solution. The appearance of smaller turbulent scales within the TGV solution as time progressed then counteracted this effect, as high wavenumbers on locally contracting regions experience excess dissipation. Lastly, a comparison was made between FR and a second order FV method. It was found that FR was more resilient to mesh deformation that FV methods, however, FR is far from optimal when running at second order. In both cases, it recommended that kinetic energy presenting/conservation methods should be used as they will greatly increase resilience to mesh quality.

\section*{Acknowledgements}
\label{sec:ack}
	The support of the Engineering and Physical Sciences Research Council of the United Kingdom is gratefully acknowledged under the award reference 1750012. This work was performed using resources provided by the Cambridge Service for Data Driven Discovery (CSD3) operated by the University of Cambridge Research Computing Service (http://www.csd3.cam.ac.uk/), provided by Dell EMC and Intel using Tier-2 funding from the Engineering and Physical Sciences Research Council (capital grant EP/P020259/1), and DiRAC funding from the Science and Technology Facilities Council (www.dirac.ac.uk).

\bibliographystyle{spmpsci}      
\bibliography{library}


\clearpage
\begin{appendices}
\section{Nomenclature}
		\begin{tabbing}
		  XXXXXXXX \= \kill
        \textit{Roman}\\
        $a$ \> convective velocity in x \\
        $b$ \> convective velocity in y \\
        $c(k)$ \> wavespeed at wavenumber $k$ \\
        $\mathbf{C}_{0\xi}$ \> centre cell FR matrix in $\xi$ \\
        $\mathbf{C}_{0\eta}$ \> centre cell FR matrix in $\eta$ \\
        $\mathbf{C}_L$ \> left cell FR matrix \\
        $\mathbf{C}_R$ \> right cell FR matrix \\
        $\mathbf{C}_B$ \> bottom cell FR matrix \\
        $\mathbf{C}_T$ \> top cell FR matrix \\
        $\mathbf{D}_{\xi}$ \> $\xi$ first derivative matrix \\
        $\mathbf{D}_{\eta}$ \> $\eta$ first derivative matrix \\
        $h_L \:\mathrm{\&}\: h_R$ \> left and right correction functions\\
        $h_B \:\mathrm{\&}\: h_T$ \> bottom and top correction functions\\
        $k_{nq}$ \> solution point Nyquist wavenumber, $(p+1)/\delta_j$\\
        $\hat{k}$ \> $k_{nq}$ normalised wavenumber, $[0,\pi]$ \\
        $l_k$ \> $k^{\mathrm{th}}$ Lagrange polynomial \\
        $p$ \> solution polynomial order \\
        $q_h$ \> element shape factor \\
        $\mathbf{Q}$ \> spatial scheme matrix \\
        $\mathbf{R}$ \> update matrix \\
        $u$ \> primitive in real domain \\	
        
        \textit{Greek}\\
        $\gamma$ \> grid geometric expansion factor \\
        $\delta_j$ \> mesh spacing\\
        $\eta$ \> $2^{\mathrm{nd}}$ computational dimension variable\\
        $\iota$ \> VCJH scheme correction function variable \\
        $\iota_+$ \> variable $\iota$ for peak temporal stability \\
        $\kappa(\mathbf{A})$ \> condition number of matrix $\mathbf{A}$ \\
        $\xi$ \>  $1^{\mathrm{st}}$ computational dimension variable \\
        $\rho(\mathbf{A})$ \> spectral radius of matrix $\mathbf{A}$ \\ 
        $\tau$ \> time step \\
        $\pmb{\Omega}$ \> solution domain \\
        $\pmb{\Omega}_n$ \> $n^{\mathrm{th}}$ solution sub-domain \\
        $\hat{\pmb{\Omega}}$ \> standardised sub-domain\\
        
        \textit{Subscript}\\
        $\mathrm{\bullet}_B$ \> variable at bottom of cell\\
        $\mathrm{\bullet}_L$ \> variable at left of cell\\
        $\mathrm{\bullet}_R$ \> variable at right of cell\\
        $\mathrm{\bullet}_T$ \> variable at top of cell\\
        \textit{Superscript}\\
        $\mathrm{\bullet}^T$ \> transpose\\
        $\mathrm{\bullet}^{\delta}$ \> local polynomial fit of value\\
        $\mathrm{\bullet}^{\delta C}$ \> correction to value\\
        $\mathrm{\bullet}^{\delta D}$ \> localised discontinuous polynomial fit of value\\
        $\mathrm{\bullet}^{\delta I}$ \> common interface value based on local polynomial fit of value\\
        $\hat{\mathrm{\bullet}}$ \> transformed variable \\
        $\overline{\mathrm{\bullet}}$ \> locally fitted polynomial of variable \\
        \textit{Other}\\
        $\hnabla$ \> gradient operator in computational domain \\
        $\Im{(z)}$ \> imagine part of $z$ given $z\in\mathbb{C}$ \\
        $\Re{(z)}$ \> real part of $z$ given $z\in\mathbb{C}$ \\
        $\mathbb{C}$ \> set of complex number \\
        $\mathbb{R}$ \> set of real numbers \\
      \end{tabbing}

\end{appendices}


\end{document}